\documentclass[a4paper,12pt, reqno]{amsart}
\usepackage{amssymb,amsmath,latexsym}
\usepackage{color}
\usepackage{comment}
\allowdisplaybreaks

\newcommand{\normal}{\color{black}}

\def\1{{\bf 1}}

\def\nn{\nonumber}
\def\bee{\begin{equation}}
\def\eee{\end{equation}}
\def\sA {{\mathcal A}} \def\sB {{\mathcal B}}

 \def\bH {{\mathbb H}} 
  
 \def\bN {{\mathbb N}}

\def\R {{\mathbb R}}  \def\bH {{\mathbb H}}

\newtheorem{thm}{Theorem}[section]
\newtheorem{lemma}[thm]{Lemma}

\newtheorem{prop}[thm]{Proposition}
\newtheorem{corollary}[thm]{Corollary}
\newtheorem{remark}[thm]{Remark}

\numberwithin{equation}{section}

\def\qed{{\hfill $\Box$ \bigskip}}
\def\NN{{\mathcal N}}

\def\CC{{\mathcal C}}

\def\DD{{\mathcal D}}
\def\FF{{\mathcal F}}
\def\EE{{\mathcal E}}

\def\R{{\mathbb R}}

\def\E{{\mathbb E}}

\def\P{{\mathbb P}}
\def\N{{\mathbb N}}
\def\eps{\varepsilon}
\def\wh{\widehat}
\def\wt{\widetilde}
\def\pf{\noindent{\bf Proof.} }

\usepackage{hyperref}
\begin{document}
\title[Markov processes with 
jump kernels decaying at the boundary]{On potential theory of Markov processes with 
jump kernels decaying at the boundary}
\author{ Panki Kim \quad Renming Song \quad and \quad Zoran Vondra\v{c}ek}
\thanks{Panki Kim: This work was  supported by the National Research Foundation of
Korea(NRF) grant funded by the Korea government(MSIP) (No. NRF-2015R1A4A1041675)
}
\thanks{Renming Song:Research supported in part by a grant from
the Simons Foundation (\#429343, Renming Song)}
\thanks{Zoran Vondra\v{c}ek: Research supported in part by the Croatian Science Foundation under the project 4197.}

 \date{}

\begin{abstract}
Motivated by some recent potential theoretic results on subordinate killed L\'evy processes
in open subsets of the Euclidean space, we study 
processes in an open set $D\subset \R^d$ defined via Dirichlet forms
with jump kernels of the form $J^D(x,y)=j(|x-y|)\sB(x,y)$ and 
critical killing functions. 
Here $j(|x-y|)$ is the L\'evy density of an isotropic stable process 
(or more generally, a pure jump isotropic unimodal L\'evy process) in $\R^d$. 
The main novelty is that the term $\sB(x,y)$ tends to 0 when $x$ or $y$ approach the boundary of $D$.
Under some general assumptions on $\sB(x, y)$, we construct the corresponding process and prove that non-negative harmonic functions of the process satisfy the Harnack inequality and 
Carleson's estimate. 
We give several examples of boundary terms satisfying those assumptions.  
The examples depend on four parameters, $\beta_1, \beta_2, \beta_3$,  $\beta_4$,  
roughly governing the decay of the boundary term near the boundary of $D$.

In the second part of this paper, we specialize to the case of the half-space 
$D=\R_+^d=\{x=(\wt{x},x_d):\, x_d>0\}$,  
the $\alpha$-stable kernel $j(|x-y|)=|x-y|^{-d-\alpha}$ and the killing function
$\kappa(x)=c x_d^{-\alpha}$, $\alpha\in (0,2)$, 
where $c$ is a positive constant.
Our main result in this part is a boundary Harnack principle which says that, 
for any $p>(\alpha-1)_+$, there are values of the parameters 
$\beta_1, \beta_2, \beta_3$,  $\beta_4$, and the constant $c$
such that non-negative harmonic functions  of the process must decay at 
the rate $x_d^p$ if they vanish near a portion of the boundary. 
We further show that there are values of the parameters 
$\beta_1, \beta_2, \beta_3$,   $\beta_4$, 
for which the boundary Harnack principle fails despite the fact that Carleson's estimate is valid.
\end{abstract}
\maketitle



\smallskip

\section{Introduction}\label{s:intro}
When studying discontinuous Markov processes,
one of the common
assumptions is that the jump kernel, 
which describes the intensity of jumps from one point in the state space to another, 
is a function of the distance (or is comparable to a function of the distance)
between these two points. 
In the context of stable (or even stable-like) processes and L\'evy processes
in the Euclidean space such an assumption is quite natural. For example, the jump kernel of 
the isotropic $\alpha$-stable process in $\R^d$ is given by $c|x-y|^{-d-\alpha}$, $d\ge 1$, 
$\alpha\in (0,2)$. 
The part process of this process in an open subset $D\subset \R^d$, 
the so-called killed $\alpha$-stable process in $D$, 
has the same jump kernel -- the intensity of jumps again depends only on the distance between the points. By removing the killing part of the killed process one is led to the censored $\alpha$-stable process with the state space $D\subset \R^d$. Still, the jump kernel remains the same.

We note that both the killed stable process and the censored stable process, although living in an open subset $D$ of  $\R^d$, are closely related to (and in some sense derived from) the stable process in $\R^d$. On the other hand, subordinate killed Brownian motions in $D\subset \R^d$, or more generally, subordinate killed L\'evy processes, are processes which are intrinsically defined in $D$ and are not a part of some larger processes living in $\R^d$. Potential theory of those processes was recently studied in \cite{KSV18-a, KSV18-b}. One of the key features of subordinate killed processes is that their jump kernels depend not only on the distance between the points, but also on the distance of each point to the boundary $\partial D$. 
This can be easily seen from 
the fact that the jump kernel of the subordinate killed process is equal to the integral of the transition density of the killed process (which depends also on the distance to the boundary) against the L\'evy measure of the subordinator. 

To be more precise, let us describe some of the results from \cite{KSV18-a, KSV18-b} in the case of stable processes and  stable subordinators. 
Let $D$ be a $C^{1,1}$ domain (connected open set) in $\R^d$. 
For $\delta\in (0,2]$, let $Z^D$ be a killed $\delta$-stable process in $D$ and 
$p^D(t,x,y)$ be its transition density.
Further, let $(T_t)$ be an independent (of $Z^D$) $\gamma/2$-stable subordinator, $\gamma\in (0,2)$, and let $Y_t^D:=Z^D_{T_t}$ be the subordinate process. Away from the boundary of $D$, the process $Y^D$ 
behaves like the $\alpha$-stable process with $\alpha=\delta \gamma/2$. 
On the other hand, the boundary behavior of the jump kernel $J^D(x,y)$ of $Y^D$, which is studied in \cite{KSV18-a, KSV18-b}, exhibits some sort of phase transition. 
In case $\delta=2$ (then $Y^D$ is a subordinate killed Brownian motion), it holds that 
\begin{equation}\label{e:intro-skbm}
J^D(x,y)\asymp \left(\frac{\delta_D(x)}{|x-y|}\wedge 1\right)\left(\frac{\delta_D(y)}{|x-y|}\wedge 1\right)|x-y|^{-d-\alpha}\, .
\end{equation}
Here and below, $a\asymp b$ means that $c\le b/a \le c^{-1}$ for some $c\in (0,1)$, 
 $a\wedge b:=\min \{a, b\}$, $a\vee b:=\max\{a, b\}$, and  $\delta_D(x)$ denotes the distance 
between $x$ and the boundary $\partial D$.

In case $\delta\in (0,2)$ (then $Y^D$ is a subordinate killed $\delta$-stable process) 
the situation is more complicated and more interesting.
It holds that (recall that $2\alpha=\delta \gamma$), 
\begin{equation}\label{e:intro-skst}
J^D(x,y)\asymp \begin{cases} 
\left(\frac{\delta_D(x)\wedge \delta_D(y)}{|x-y|}\wedge 1\right)^{\delta(1-\gamma/2)}|x-y|^{-d-\alpha},\quad \gamma\in (1,2),\\
\left(\frac{\delta_D(x)\wedge \delta_D(y)}{|x-y|}\wedge 1\right)^{\delta/2}
\log\left(1+\frac{(\delta_D(x)\vee \delta_D(y))\wedge |x-y|}{\delta_D(x)\wedge \delta_D(y)\wedge |x-y|}\right)
|x-y|^{-d-\alpha}, \quad \gamma=1, \\
\left(\frac{\delta_D(x)\wedge \delta_D(y)}{|x-y|}\wedge 1\right)^{\delta/2}\left(\frac{\delta_D(x)\vee \delta_D(y)}{|x-y|}\wedge 1\right)^{ (\delta/2)(1-\gamma) }|x-y|^{-d-\alpha},\quad \gamma\in (0,1).
\end{cases}
\end{equation}
Define $\sB(x,y):=J^D(x,y)|x-y|^{d+\alpha}$ so that $J^D(x,y)=\sB(x,y)|x-y|^{-d-\alpha}$. Then we can think of $\sB(x,y)$ as the boundary term of the jump kernel of $Y^D$ 
which depends on $\delta_D(x), \delta_D(y)$ and $|x-y|$ (as opposed to the term $|x-y|^{-d-\alpha}$ which depends solely on $|x-y|$).
Sharp two-sided estimates for $\sB(x,y)$ can be seen from \eqref{e:intro-skbm} and \eqref{e:intro-skst}. 

We note that the process $Y^D$ is transient with 
its killing function $\kappa^{Y^D}$ satisfying the lower bound 
$\kappa^{Y^D}(x)\ge c \delta_D(x)^{-\alpha}$, cf. \cite[(2.4)]{KSV18-b}.
Note that the killing function is not in the usual Kato class, 
i.e., $\lim_{t \to 0}\sup_{x\in D}\E_x\int_0^t \kappa^{Y^D}(Y^D_s)ds$ is not zero.

One of the main goals of \cite{KSV18-a, KSV18-b} was to prove a (scale invariant) boundary Harnack principle (BHP) with exact decay rate for non-negative functions harmonic with respect to $Y^D$ in case $D$ is a $C^{1,1}$ domain. 
For the subordinate killed Brownian motion ($\delta=2$) 
it was shown in \cite[Theorem 1.2]{KSV18-a} that BHP holds with the decay rate $\delta_D(x)$. 
In case of the subordinate killed $\delta$-stable process via an independent $\gamma/2$-stable subordinator  ($\delta, \gamma\in (0,2)$), 
the scale invariant BHP holds in case $\gamma\in (1,2)$ with the decay rate $\delta_D(x)^{\delta/2}$, cf.~\cite[Theorem 7.2]{KSV18-b}. Surprisingly, it turns out that, despite the fact that Carleson's estimate is valid, the (non-scale invariant) BHP fails when $\gamma\in (0,1]$, see \cite[Section 9]{KSV18-b}.

The goal of this paper is to study processes in open subsets of $\R^d$ associated with Dirichlet
forms with jump kernels decaying at the boundary and critical killing functions,
and to build a general framework which includes the processes studied in \cite{KSV18-a, KSV18-b} as examples.
To be more precise, 
let $X=(X_t, \P_x)$ be a pure-jump isotropic  unimodal L\'evy process whose 
 jump kernel $j(x,y)=j(|x-y|)$  satisfies 
 \begin{equation}\label{e:assumption-on-j}
 j(r)\asymp r^{-d}\Phi(r)^{-1}, \quad \text{for all } 
r>0,
\end{equation}
where $\Phi$ is an increasing function satisfying  the following weak scaling condition: 
There exist constants $0<\delta_1\le \delta_2<1$ and $a_1, a_2>0$ such that 
\begin{equation}\label{e:H-infty}
   a_1(R/r)^{2\delta_1}\leq\frac{\Phi(R)}{\Phi(r)}\leq a_2(R/r)^{2\delta_2},\quad 
0<r<R <\infty.
\end{equation}
A prototype of such a process $X$ is the isotropic $\alpha$-stable process in which case $\Phi(r)=r^{\alpha}$. This particular case already contains all the essential features of our results.

For a given open set $D\subset \R^d$, we consider the process $Y^D$ on $D$
associated with a pure jump Dirichlet form whose jump kernel has the form $J^D(x,y)=\sB(x,y)j(|x-y|)$.
Here $\sB(x,y)$ -- the boundary term --
depends on $\delta_D(x), \delta_D(y)$ and $|x-y|$, and is allowed to approach 0 at the boundary.
This is in contrast with previous works (see, for instance, \cite{CK03, CK08}) where $\sB(x,y)$
is assumed to be bounded between two positive constants, which can be viewed as a uniform ellipticity condition for non-local operators. In this sense our paper is the first systematic attempt
to study the potential theory of degenerate non-local operators.
Throughout this paper we will assume that $\sB$ is symmetric and bounded, namely

\medskip
\noindent
\textbf{(B1)} $\sB(x,y)=\sB(y,x)$ for all $x,y\in D$.

\medskip
\noindent
\textbf{(B2)} There exists a constant $C_2>0$ such that $\sB(x,y)\le C_2$ for all $x,y\in D$.

\medskip
Under these two conditions, the jump kernel $J^D$ gives rise to a regular Dirichlet
form $(\FF^D, \EE^D)$ on $L^2(D, dx)$. Thus (see, for instance, 
\cite[Example 1.2.4]{FOT}) there exists a Hunt process $Y^D$ associated with $(\FF^D, \EE^D)$.
Such a process is defined up to an exceptional set $\mathcal N$. 
We next kill the process $Y^D$ by using a  killing function $\kappa:D\to [0,\infty)$ satisfying
\begin{equation}\label{e:bound-on-kappa}
\kappa(x)\le C_1 \frac{1}{\Phi(\delta_D(x))}\, ,\quad x\in D\, ,
\end{equation}
for some $C_1>0$, to obtain the Hunt process $Y^{D,\kappa}$. 
Once killed, the process is sent to the cemetery $\partial$. We will use the convention that 
every function is automatically extended to be zero at $\partial$.
We further impose the following general conditions on the boundary term $\sB(x,y)$. These conditions 
are satisfied by a quite a few classes of examples.

\medskip
\noindent
\textbf{(B3)} For every $a\in (0,1)$ there exists $C_3=C_3(a)>0$ such that $\sB(x,y)\ge C_3$ whenever $\delta_D(x)\wedge \delta_D(y)\ge a|x-y|$.

\medskip
Condition \textbf{(B3)} ensures that, away from the boundary, the jumping kernel $J^D(x,y)$ is comparable to $j(|x-y|)$. In Subsection \ref{ss:regularization} we use this condition to remove the exceptional set $\mathcal N$ so that the process $Y^{D,\kappa}$ can start from every point in $D$. The next condition is needed in the analysis of the generator of the process  in Subsection \ref{ss:dynkin}.  
\medskip

\medskip
\noindent
\textbf{(B4)}  
If $\delta_2 \ge 1/2$, then there exist $\theta>2\delta_2-1$ and $C_4>0$ such that 
$$
|\sB(x, x)-\sB(x,y)|\le C_4\left(\frac{|x-y|}{\delta_D(x)\wedge \delta_D(y)}\right)^{\theta}\,.
$$

\medskip
\noindent
\textbf{(B5)} For every  $\epsilon \in (0,1)$ there exists $C_5=C_5(\epsilon)\ge 1$ with the following property: For all $x_0\in D$ and $r>0$ with $B(x_0, (1+\epsilon)r)\subset D$, we have
$$
C_5^{-1}\sB(x_1,z)\le \sB(x_2,z) \le C_5 \sB(x_1, z)\, ,\quad \text{for all }x_1,x_2\in B(x_0,r), \, \,z\in D\setminus B(x_0, (1+\epsilon)r)\, .
$$

Roughly speaking, this condition says that 
jumping intensity 
from points close to each other to a faraway point is comparable.
We note that \textbf{(B3)} implies that $\sB(x, y)>0$ for all $x, y\in D$.

Under these conditions we 
prove two versions of Harnack inequality for non-negative harmonic functions  
with respect to $Y^{D,\kappa}$.
Recall that 
a non-negative Borel function $f$ defined on $D$ is said to be  \emph{harmonic} in an open set $V\subset D$ 
with respect to $Y^{D,\kappa}$ if for every bounded open set $U\subset\overline{U}\subset V$,
$$
f(x)= \E_x \left[ f(Y^{D,\kappa}_{\tau_{U}})\right] \text{ and } 
\E_x \left[f(Y^{D,\kappa}_{\tau_{U}})\right]  <\infty, \qquad
\hbox{for all } x\in U,
$$
where $\tau_U=\tau_U^{Y^{D,\kappa}}:=\inf\{t>0:\, Y^{D,\kappa}_t\notin U\}$ is 
the first exit time 
of $Y^{D,\kappa}$ from $U$.

\begin{thm}[Harnack inequality]\label{t:uhp}
Suppose $D\subset \R^d$ is an open set and assume that 
\textbf{(B1)-(B5)} and \eqref{e:assumption-on-j}-\eqref{e:bound-on-kappa} hold true.

\begin{itemize}
\item[(a)] 
There exists a constant $C>0$ such that 
for any $r\in (0,1]$, $B(x_0, r) \subset D$ and any non-negative function $f$ in $D$ which is 
harmonic in $B(x_0, r)$ with respect to $Y^{D,\kappa}$, we have
$$
f(x)\le C f(y), \qquad \text{ for all } x, y\in 
B(x_0, r/2).
$$

\item[(b)] 
There exists a constant $C>0$ 
such that for anyy  $L>0$, any $r \in (0, 1]$, 
all $x_1,x_2 \in D$ with $|x_1-x_2|<Lr$ and $B(x_1,r)\cup B(x_2,r) \subset D$ 
and any non-negative function $f$ in $D$ which is  harmonic in $B(x_1,r)\cup B(x_2,r)$ with respect to $Y^{D, \kappa}$ we have
$$
f(x_2)\le C C_3(1/(2(L+1)))^{-1} L^{d+2\delta_2} f(x_1)\, .
$$
\end{itemize}
\end{thm}

We then proceed towards Carleson's estimate. 
Here we need two additional assumptions: one on the killing function $\kappa$,  the other on the boundary term $\sB$.
 In addition to \eqref{e:bound-on-kappa} we assume the corresponding lower bound: There exists $C_8>0$ such that
\begin{equation}\label{e:l-bound-on-kappa}
\kappa(x)\ge C_8 \frac{1}{\Phi(\delta_D(x))}\, ,\quad x\in D\, .
\end{equation}
We note  in passing   that the killing function of a subordinate killed Brownian motion in a
Lipschitz 
domain satisfies \eqref{e:bound-on-kappa} and \eqref{e:l-bound-on-kappa}, cf. \cite[Proposition 3.2]{SV03}.
For the boundary term we assume

\medskip
\noindent
\textbf{(B6)} There exist $\wh{\beta}>0$ and $C_6>0$ such that if 
$x, y, z\in D$ satisfy $\delta_D(x)\le \delta_D(z)$ and $|y-z|\le M |y-x|$ with $M\ge 1$, then
$
\sB(x,y)\le C_6 M^{\wh{\beta}}\sB(z,y)\, .
$

An open set $D\subset \R^d$ is $\kappa$-fat if there are $\overline{\kappa}\in (0,1/2]$ and $R>0$ such that for all $x\in \overline{D}$ and all $r\in (0,R]$, there is a ball $B(A_r(x), \kappa r)\subset D\cap B(x,r)$. The pair $(R, \overline{\kappa})$ is called the characteristics of the $\kappa$-fat open set $D$.

\begin{thm}[Carleson estimate]\label{t:carleson}
Suppose that $D\subset \R^d$ is a  $\kappa$-fat open set with characteristics 
$(R, \overline{\kappa})$. Assume that \textbf{(B1)-(B6)} and 
\eqref{e:assumption-on-j}-\eqref{e:l-bound-on-kappa} hold true. There exists a constant $C=C(R,\overline{\kappa})>0$  such that for every 
$Q\in\partial D$, $0<r<(R \wedge 1)/2$, and every 
non-negative  Borel function $f$ in $D$ which is harmonic in $D \cap B(Q, r)$ with respect to
$Y^{D,\kappa}$ and vanishes continuously on $ \partial D \cap B(Q, r)$, we have
\begin{equation}\label{e:carleson}
f(x)\le C f(x_0) \qquad \hbox{for all }  
x\in D\cap B(Q,r/2),
\end{equation}
where $x_0\in D\cap B(Q,r)$ with $\delta_D(x_0)\ge \overline{\kappa}r/2$.
\end{thm}

These results comprise the first part of the paper and are foundations for the second part in which we study in more detail the boundary behavior of non-negative harmonic functions in case $D=\R^d_+=\{x=(\wt{x}, x_d): x_d>0\}$, $j(x,y)=|x-y|^{-d-\alpha}$ and $\kappa(x)=c x_d^{-\alpha}$ 
(where $c$ is any positive constant).
The main goal is to explore conditions on $\sB$ which ensure that the BHP holds true and find the rate of decay of non-negative harmonic functions. 
Motivated by the three examples of the boundary term in Section \ref{s:assumptions} that satisfy conditions \textbf{(B1)}-\textbf{(B6)}, in the second part we assume that the 
boundary term $\sB(x,y)$ satisfies \textbf{(B1)}, \textbf{(B4)}
and the following sharp two-sided estimates:

\medskip
\noindent
\textbf{(B7)}
There exist $C_7\ge 1$ and  $\beta_1, \beta_2, \beta_3, \beta_4 \ge 0$ 
with $\beta_1>0$ if $\beta_3 >0$,
 and $\beta_2>0$ if $\beta_4 >0$ 
 such that
\begin{equation}\label{e:B7}
C_7^{-1}\wt{B}(x,y)\le \sB(x,y)\le C_7 \wt{B}(x,y)\, ,\qquad x,y\in \R^d_+\, ,
\end{equation}
where $\wt{B}(x,y)=\wt{B}_{\beta_1, \beta_2, \beta_3,  \beta_4 }(x,y)$ is 
$$
\Big(\frac{x_d\wedge y_d}{|x-y|}\wedge 1\Big)^{\beta_1}\Big(\frac{x_d\vee y_d}{|x-y|}\wedge 1\Big)^{\beta_2} \left[ \log\Big(1+\frac{(x_d\vee y_d) \wedge |x-y|}{x_d\wedge y_d\wedge |x-y|}\Big)\right]^{\beta_3} 
\left[ \log\Big(1+\frac{ |x-y|}{(x_d\vee y_d)\wedge |x-y|}\Big)\right]^{\beta_4} . 
$$
The first three terms 
in the definition of  $\wt{B}(x,y)$ with the parameters $\beta_1, \beta_2$ and $ \beta_3$
 are motivated by \eqref{e:intro-skst}, and
 the last term  in the definition of  $\wt{B}(x,y)$ with
  the parameter $\beta_4$
is motivated by \cite{CKSV21} (see the paragraph above Lemma \ref{l:B-tilde-assumptions}).

We show in Section \ref{s:assumptions} that 
the boundary term $\wt{B}(x,y)$ satisfies   \textbf{(B1)}-\textbf{(B6)}. 
Thus $\sB$ also satisfies \textbf{(B2)}, \textbf{(B3)}, \textbf{(B5)} and \textbf{(B6)}.

The following additional assumption on $\sB(x,y)$ is natural for the state space $\R^d_+$. 

\medskip
\noindent
\textbf{(B8)} 
For all $x,y\in \R^d_+$ $a>0$, $\sB(ax,ay)=\sB(x,y)$. In case $d\ge 2$, for  all $x,y\in \R^d_+$ and $\wt{z}\in \R^{d-1}$, $\sB(x+(\wt{z},0), y+(\wt{z},0))=\sB(x, y)$.
\medskip

We now describe the BHP we establish in more detail. 
 Let $J^{\R_+^d}(x,y)=|x-y|^{-d-\alpha}\sB(x,y)$, $x,y\in \R_+^d$, where the boundary term $\sB(x,y)$ 
satisfies \textbf{(B1)}, \textbf{(B4)}, \textbf{(B7)} and \textbf{(B8)}.
To every parameter $p\in ((\alpha-1)_+, \alpha+\beta_1)$ we associate a constant 
$C(\alpha, p, \sB)\in (0, \infty)$
depending on $\alpha$, $p$ and $\sB$ 
defined, in case $d\ge 2$, as
\begin{equation}\label{e:explicit-C}
C(\alpha, p, \sB)=
\int_{\R^{d-1}}\frac{1}{(|\wt{u}|^2+1)^{(d+\alpha)/2}}\left( \int_0^1 \frac{(s^p-1)(1-s^{\alpha-p-1})}{(1-s)^{1+\alpha}}
\sB\big((1-s)\wt{u}, 1), s\mathbf{e}_d \big)\, ds\right)d\wt{u}\, 
\end{equation}
where $\mathbf{e}_d=(\tilde{0}, 1)$.
In case $d=1$, $  C(\alpha, p, \sB)$ is defined as
$$
C(\alpha, p, \sB)=\int_0^1 \frac{(s^p-1)(1-s^{\alpha-p-1})}{(1-s)^{1+\alpha}}
\sB\big(1, s \big)\, ds.
$$
Note that the function $p\mapsto C(\alpha, p, \sB)$ is non-decreasing,  $\lim_{p\downarrow (\alpha-1)_+}C(\alpha, p, \sB)=0$ and that
  $\lim_{p\uparrow \alpha+\beta_1}C(\alpha, p, \sB)=\infty$ 
(see Lemma \ref{l:LB-on-g} and Remark \ref{r:C-increasing}).

Let
\begin{align}
\label{e:kappap}
\kappa(x)
=C(p,\alpha, \sB)x_d^{-\alpha}, \qquad x\in \R^d,
\end{align} 
be the killing function.
Note that $\kappa$ depends on $p$, but we omit this dependence from the notation for simplicity.
We denote  by $Y^{\R^d_+, \kappa}$  the corresponding process killed by $\kappa$.

We are now ready to state the BHP.
We will only give the statement of the result for $d\ge 2$. The statement in the $d=1$ case is similar and simpler.
For any $a, b>0$ and $\wt{w} \in \R^{d-1}$, we define $D_{\wt{w}}(a,b):=\{x=(\wt{x}, x_d)\in \R^d:\, |\wt{x}-\wt{w}|<a, 0<x_d<b\}$.

\begin{thm}\label{t:BHP}
Assume that \textbf{(B1)}, \textbf{(B4)}, \textbf{(B7)}-\textbf{(B8)} and \eqref{e:kappap} hold true.
Suppose that $p\in ((\alpha-1)_+, \alpha+\beta_1)$
and either (a) $\beta_1=\beta_2=\beta\ge0$ and 
$\beta_3=\beta_4=0$,
 or (b) $p<\alpha$.
Then there exists $C>0$ such that for every $r>0$, $\wt{w} \in \R^{d-1}$, any non-negative function $f$ in $\R^d_+$ which is harmonic in $D_{\wt{w}}(2r, 2r)$ with respect to $Y^{\R_+^d, \kappa}$ and vanishes continuously on 
$B({\wt{w}}, 2r)\cap \partial \R^d_+$, we have 
\begin{equation}\label{e:TAMSe1.8}
\frac{f(x)}{x^p_d}\le C\frac{f(y)}{y^p_d}, 
\quad x, y\in D_{\wt{w}}(r/2, r/2) .
\end{equation}
\end{thm}

Note that the power $p$ in the decay rate comes from the constant $C(\alpha, p, \sB)$ in the killing function.
In \cite{BBC}, a BHP is established for censored $\alpha$-stable processes, $\alpha\in (1, 2)$, in smooth domains. The censored process 
has no killing,  $\wt{B}(x,y) \equiv 1$, 
and the explicit decay rate $p$ in the BHP of \cite{BBC} is $\alpha-1$.

In case $\beta_1=\beta_2=1$, $\beta_3= \beta_4= 0$ and $p=1$, 
this result covers the BHP (for $D=\R_+^d$) from \cite{KSV18-a}. In case 
$\frac{\alpha}2<p<\alpha$, $\beta_1=2p-\alpha$ and 
$\beta_2=\beta_3= \beta_4=0 $
it covers the BHP from \cite{KSV18-b}.
The main novelty of Theorem \ref{t:BHP} comes from  part (a)
where we can allow arbitrarily large $\beta$, hence arbitrarily large $p$, implying very fast decay rate of harmonic functions. For fixed $\beta$, this decay rate is a consequence of the constant $C(\alpha, p, \sB)$ in the killing term $\kappa$. By increasing $\beta$ we enlarge the possible range for parameter $p\in ((\alpha-1)_+, \alpha+\beta)$. 
Note that when $\sB(x, y)\equiv 1$ (corresponding to the case   
$\beta_1=\beta_2=\beta_3= \beta_4= 0 $),  
the possible range for the parameter $p$ is $((\alpha-1)_+, \alpha)$. This says
that the power $p$ in the decay rate of non-negative harmonic functions of $-(-\Delta)^{\alpha/2}$ with critical killing can not be higher than $\alpha$ (see also \cite{CKSV19}). This is in sharp contrast with the case of Laplacian with critical killing, where the power $p$ can be arbitrarily large. Indeed, for $f(x)=x_d^p$, it clearly holds that $\Delta f(x)-p(p-1)x_d^{-2}f(x)=0$ for all
$p\in (1, \infty)$.
The second novelty of Theorem \ref{t:BHP}
 is that in case $p<\alpha$ we can allow $\beta_2>0$ as well as 
 logarithmic terms in $\sB(x,y)$.

On the other hand, similarly as in \cite{KSV18-b}, for certain values of the parameters 
$\beta_1$, $\beta_2$, $\beta_3$ 
 and  $\beta_4$, 
we can show that BHP fails. We say that the non-scale-invariant BHP holds near the boundary of $\R^d_+$ if there is a constant $\wh{R}\in (0,1)$  such that for any $r \in (0, \wh{R}\, ]$, there exists  a constant $c=c(r)\ge 1$ such that for all $Q\in \partial \R^d_+$ and non-negative functions $f, g$ in $\R^d_+$ which are harmonic  in $\R^d_+ \cap B(Q,r)$ with respect to $Y^{\R^d_+, \kappa}$ and vanish continuously on $ \partial \R^d_+ \cap B(Q, r)$, we have
$$
\frac{f(x)}{f(y)}\,\le c\,\frac{g(x)}{g(y)} \qquad
\hbox{for all } x, y\in  
(\R^d_+ \cap B(Q, r/2)).
$$ 
\begin{thm}\label{t:counterexample}
Assume that 
$d \ge 2$ and \textbf{(B1)}, \textbf{(B4)}, \textbf{(B7)}-\textbf{(B8)} and \eqref{e:kappap} hold true.
Suppose that $\alpha+\beta_2<p<\alpha+\beta_1$.
Then the non-scale-invariant boundary Harnack principle is not valid for $Y^{\R_+^d, \kappa}$. 
\end{thm}

Theorems \ref{t:carleson} and  \ref{t:counterexample}  imply that, when $\alpha+\beta_2<p<\alpha+\beta_1$, under the assumptions 
\textbf{(B1)}, \textbf{(B4)}, \textbf{(B7)}-\textbf{(B8)} and \eqref{e:kappap},  the Carleson estimate holds but the BHP fails.
This is quite interesting since it is known that, for diffusions, the Carleson estimate and BHP are 
equivalent. See \cite{A} and the references therein.

Now we describe the structure of the paper. We start with a section on examples of the boundary term and check that they satisfy all of the conditions \textbf{(B1)}-\textbf{(B8)}. The examples are motivated by and are extensions of the boundary terms of the processes studied in \cite{KSV18-a, KSV18-b}. 
This section is not essential for understanding of the main development later in the paper and the reader may wish to only glance through it first.

Section \ref{s:interior} is on the interior results and uses only assumptions \textbf{(B1)}-\textbf{(B5)} and 
\eqref{e:assumption-on-j}-\eqref{e:bound-on-kappa}. 
It consists of three subsections. In the first one we remove the exceptional set $\NN$. The main step towards this is the identification of the process $Y^{D, \kappa}$ killed upon exiting a relatively compact $C^{1,1}$-open subset $U$ of $D$ with a certain process which can start from every point of $U$. In Subsection \ref{ss:dynkin} we analyze the generator $L^\sB$ of the process $Y^{D, \kappa}$, cf.~\eqref{e:defn-LBD}. The difficulty comes from the boundary  term  which destroys the symmetry at the singularity. Assumption \textbf{(B4)} plays a central role in getting around this difficulty. An important result of this subsection is Proposition \ref{p:Dynkin} which gives a Dynkin formula for smooth functions of compact support. 
Note that our process may not be a Feller process.
In Subsection \ref{HI} we prove Theorem \ref{t:uhp} by following the approach in \cite{BL, SV04}.

In Section \ref{s:carleson} we add assumption \textbf{(B6)} and the lower bound \eqref{e:l-bound-on-kappa} on the killing function $\kappa$,  and prove Theorem \ref{t:carleson}. The proof follows the method already used before in different contexts. This section concludes the first part of the paper.

In the second part of this paper, we specialize to the case of the half-space $D=\R_+^d=\{x=(\wt{x},x_d):\, x_d>0\}$,  the $\alpha$-stable kernel $j(|x-y|)=|x-y|^{-d-\alpha}$ and the killing function $\kappa(x)=c x_d^{-\alpha}$, $\alpha\in (0,2)$. We assume that \textbf{(B1)}, \textbf{(B2)}, \textbf{(B4)}, \textbf{(B7)} and \textbf{(B8)} hold true.
The main purpose of Section \ref{s:estimates} is to establish estimates related to exit probabilities from small boxes near the boundary. To accomplish this, we first show that the action of the operator $L^\sB$ on the function $x\mapsto x_d^p$ is equal to 0. This indicates that if the BHP holds, the correct decay rate has to be $x_d^p$. Then we establish the following two-sided estimates:
\begin{equation}\label{e:intro-1}
\E_x \int_0^{\tau_U} (Y_t^d)^{\beta_1}|\log Y_t^d|^{\beta_3} \, dt \asymp x_d^p\, ,\quad x\in U\, ,
\end{equation}
where $U:=\{x=(\wt{x}, x_d)\in \R^d_+:\, |\wt{x}|< \frac12, x_d<\frac12\}$ and 
$Y=(Y^1, \dots, Y^d)$ stands for $Y^{\R^d_+, \kappa}$. 
One direction is relatively easy, while the other requires constructing suitable test functions and estimating the action of the operator $L^\sB$ on them, cf.~Lemma \ref{l:key-lemma}. This turns out to be a formidable task, and therefore we postpone the proof of this lemma to Section \ref{s:proof-of-key-lemma}. Another important ingredient in the proof of \eqref{e:intro-1} as well as some other results from Section \ref{s:estimates} are some extensions of Proposition \ref{p:Dynkin} to functions which are neither smooth nor of compact support. In order not to interrupt the flow of the presentation, we relegate the rigorous proofs of the extensions to Section \ref{s:10}. Estimates \eqref{e:intro-1} are used in a rather straightforward way to establish the already mentioned exit probability estimates from small boxes, cf.~Lemma \ref{l:exit-probability-estimate}. These exit probability estimates together with \eqref{e:intro-1}  are the main tools in proving the BHP in Section \ref{s:bhp}. Additionally, Section \ref{s:estimates} contains the following two-sided estimates when $\alpha<p<\alpha+\beta_1$ and $\alpha+\beta_2<p$:
\begin{equation}\label{e:intro-2}
\E_x \int_0^{\tau_U} (Y_t^d)^{\beta_2}\, dt \asymp x_d^{\alpha+\beta_2}\, ,\quad x\in U\, .
\end{equation}
These will be used in proving
the failure of the BHP for the above range of parameters in Section \ref{s:counter}. 

Sections \ref{s:bhp} and \ref{s:counter} contain proofs of Theorem \ref{t:uhp} and Theorem \ref{t:counterexample}, respectively. These are modeled after the corresponding proofs in 
\cite{KSV18-b} and 
use the box method (first adapted to the case jump processes in \cite{BB}), the Harnack inequality and Carleson's estimate.

Throughout this paper, the positive constants 
$\beta_1$, $\beta_2$, $\beta_3$, $\beta_4$,
$\theta$, $\delta_1$, $\delta_2$, $\delta_*$, $\wh{\beta}$  will remain the same.
We will use the following convention:
Capital letters $C, C_i, i=1,2,  \dots$ will denote constants
in the statements of results and assumptions. The labeling of these constants will remain the same. Lower case letters 
$c, c_i, i=1,2,  \dots$ are used to denote the constants in the proofs
and the labeling of these constants starts anew in each proof.
The notation $c_i=c_i(a,b,c,\ldots)$, $i=0,1,2,  \dots$ indicates  constants depending on $a, b, c, \ldots$.
We will use ``$:=$" to denote a
definition, which is read as ``is defined to be".
For any $x\in \R^d$, $r>0$ and $0<r_1<r_2$, we use  $B(x, r)$ to denote the
open ball of radius $r$ centered at $x$ and use $A(x, r_1, r_2)$ to denote the
annulus $\{y\in \R^d: r_1\le |y-x|<r_2\}$.
For a Borel subset $V$ in $\R^d$, $|V|$ denotes  the Lebesgue measure of $V$ in $\R^d$,  $\delta_U:=\mathrm{dist}(U, \partial D)$ and $d_U:=\mathrm{diam}(U)$.
For  a function space $\bH(U)$
on an open set $U$ in $\R^d$, we let   
$\bH_c(U):=\{f\in\bH(U): f \mbox{ has  compact support}\}$ and $\bH_b(U):=\{f\in\bH(U): f \mbox{ is bounded}\}$. 
We use the shorthand notation $\log(1+a)^{\gamma}$ for the more precise $(\log(1+a))^{\gamma}$.

For the convenience of the reader
we give below the table of important constants and operators used throughout the paper.

\begin{table}[h!]
\begin{tabular}{*{2}{c}}
\hline
Notation &    Description 
\\ \hline 
\medskip
$\delta_1, \delta_2$  & Indices of weak scaling condition of $\Phi$ in \eqref{e:H-infty}  \\
\medskip
$\beta_1, \beta_2, \beta_3, \beta_4$ & Parameters in the definition of $\wt{B}(x,y)$ in \textbf{(B7)} \\\medskip
$C(\alpha, p, \sB)$ & The constant in the killing function $\kappa$ defined in \eqref{e:explicit-C} \\
\medskip
$L^\sB$ &The  extended infinitesimal generator of the process $Y^{D, \kappa}$ defined in \eqref{e:defn-LBD}\\
\medskip
$L^\sB_\alpha$ &
The extended  infinitesimal generator of the process $Y^{\R^d_+, 0}$ defined in \eqref{d:Lalpha}. 
\\
\medskip
\end{tabular}
\caption{
List of of the important notations
}
\end{table}

\section{Examples of the boundary term }\label{s:assumptions}
In this section we give examples of the boundary term $\sB(x,y)$ that satisfy 
assumptions \textbf{(B1)-(B8)}. 
Our motivation comes from the boundary terms appearing in \eqref{e:intro-skst}. 

We start with three elementary lemmas.  
We recall that a function $f:I\to \R$, $I\subset \R$, is said to be almost increasing if there exists $c>0$ such that $f(x)\le c f(y)$ whenever $x,y\in I$ and $x<y$. An almost decreasing function is defined analogously.

Let $\beta_1, \beta_2, \beta_3,  \beta_4 \ge 0$ so that $\beta _1 >0$ if $ \beta_3 >0$,  and $\beta_2>0$ if $\beta_4>0$. 
For $s,t,u>0$, define
$$
\sA(s,t,u):=\left(\frac{s\wedge t}{u}\wedge 1\right)^{\beta_1}\left(\frac{s\vee t}{u}\wedge 1\right)^{\beta_2}\ {\log\left(1+\frac{(s\vee t)\wedge u}{s\wedge t\wedge u}\right)^{\beta_3}}  \ {\log\left(1+\frac{ u}{(s\vee t)\wedge u}\right)^{\beta_4}}\, . 
$$

\begin{lemma}\label{l:2-assumptions} 
\begin{itemize}
\item[(a)] 
 The function $u\mapsto \sA(s,t,u)$ is almost decreasing on $(0,\infty)$.
\item[(b)]  The function $s\mapsto \sA(s,t,u)$ is almost increasing on $(0,\infty)$.
\end{itemize}
\end{lemma}
\pf 
In the proof we will use the following fact (which can be proved by elementary calculus) several times:
\begin{align}
\label{l:1-assumptions}
\text{For $\beta>0$ and $\gamma\ge 0$,
$s\mapsto s^{\beta}\log\left(1+{s}^{-1}\right)^{\gamma}$
is almost increasing on $(0,1]$.}
\end{align}

\noindent
(a) Note that $\sA$ is symmetric in $s$ and $t$, hence we may assume that $s\ge t$. Note that
$$\sA(s,t,u)=
\begin{cases}
(\log 2)^{\beta_3  + \beta_4 \normal}, &\text{if } s\ge t\ge u;\\
\left(\frac{t}{u}\right)^{\beta_1}{\log\left(1+\frac{u}{t}\right)^{\beta_3}} (\log 2)^{\beta_4}, 
&\text{if }s\ge u\ge t;\\
\left(\frac{t}{u}\right)^{\beta_1}\left(\frac{s}{u}\right)^{\beta_2}{\log\left(1+\frac{s}{t}\right)^{\beta_3}}  \log\left(1+\frac{u}{s}\right)^{\beta_4}, \normal &\text{if }u\ge s\ge t.
\end{cases}
$$
 Thus,  by \eqref{l:1-assumptions}, we see that $u\mapsto \sA(s,t,u)$ is almost decreasing on $(0, \infty)$. 

\noindent (b)
Here we distinguish two cases:

\noindent
Case (i): $t\ge u$. Then
$$\sA(s,t,u)=
\begin{cases}
\left(\frac{s}{u}\right)^{\beta_1}{\log\left(1+\frac{u}{s}\right)^{\beta_3}} (\log 2)^{\beta_4},  \normal &\text{if }   s\le u\le t	;\\
(\log 2)^{\beta_3  +\beta_4  },&\text{if }	u\le s\le t \text{ or }u\le t\le s.
\end{cases}
$$
It follows from \eqref{l:1-assumptions} that $s\mapsto \sA(s,t,u)$ is almost increasing. 

\noindent
Case (ii): $t\le u$. Then
$$\sA(s,t,u)=
\begin{cases}
\left(\frac{s}{t}\right)^{\beta_1}\left(\frac{t}{u}\right)^{\beta_1+\beta_2}{\log\left(1+\frac{t}{s}\right)^{\beta_3}}  {\log\left(1+\frac{u}{t}\right)^{\beta_4}}, &\text{if }	s\le t\le u;\\
\left(\frac{t}{u}\right)^{\beta_1}\left(\frac{s}{u}\right)^{\beta_2}{\log\left(1+\frac{s}{t}\right)^{\beta_3}}  {\log\left(1+\frac{u}{s}\right)^{\beta_4}}, &\text{if }	t\le s\le u;\\
\left(\frac{t}{u}\right)^{\beta_1}{\log\left(1+\frac{u}{t}\right)^{\beta_3}}  (\log 2)^{\beta_4},   &\text{if }	t\le u\le s.
\end{cases}
$$
By \eqref{l:1-assumptions}, the function $s\mapsto \sA(s,t,u)$ is almost increasing on $(0,t]$  and almost increasing on $[t,u]$. It is clearly  constant on $[u,\infty)$. 
Thus, $s\mapsto \sA(s,t,u)$ is almost increasing on $(0,\infty)$. \qed

Let $D\subset \R^d$ be an open set, and recall that $\delta_D(x)$ denotes the distance  
between $x\in D$ and the boundary $\partial D$. 
For simplicity, let
$$
L(x,y):={\log\left(1+\frac{(\delta_D(x)\vee \delta_D(y))\wedge |x-y|}{\delta_D(x)\wedge \delta_D(y)\wedge |x-y|}\right)}
$$ 
and 
$$
K(x,y)={\log\left(1+\frac{|x-y|}{(\delta_D(x)\vee \delta_D(y))\wedge |x-y|}\right)}.
$$
Let $\beta_1, \beta_2, \beta_3,  \beta_4 \ge 0$ so that $\beta _1 >0$ if $ \beta_3 >0$,  and $\beta_2>0$ if $\beta_4>0$.  We define the function $\wt{B}:D\times D\to (0, \infty)$ by
\begin{equation}\label{e:def-of-B-tilde}
\wt{B}(x,y):=\left(\frac{\delta_D(x)\wedge\delta_D(y)}{|x-y|}\wedge 1\right)^{\beta_1}\left(\frac{\delta_D(x)\vee \delta_D(y)}{|x-y|}\wedge 1\right)^{\beta_2} L(x,y)^{\beta_3}  K(x,y)^{\beta_4}  \, .
\end{equation}
Observe that $\lim_{y\to x}\wt{B}(x,y)=(\log 2)^{\beta_3  +\beta_4 }$, so for $y=x$ we interpret the right-hand side above as $\wt{B}(x,x)=(\log 2)^{\beta_3  + \beta_4 }$.
Note that $\wt{B}(x,y)=\sA(\delta_D(x), \delta_D(y), |x-y|)$.

Note that when $\beta_1=\beta_2=\beta$ and $\beta_3  =\beta_4=0$, we have
\begin{equation}\label{e:beta1=beta2}
\wt{B}(x, y)= \left(\frac{\delta_D(x)}{|x-y|}\wedge 1\right)^{\beta}\left(\frac{\delta_D(y)}{|x-y|}\wedge 1\right)^{\beta}\, ,
\end{equation}
which is, in case $\beta=1$, comparable to the boundary term in \eqref{e:intro-skbm}. When $\beta_1=\delta(1-\gamma/2)$, $\beta_2=0$ and $\beta_3   =\beta_4=0$,  $\wt{B}(x,y)$ is comparable to the boundary term in the first line of \eqref{e:intro-skst}.  When $\beta_1=\delta/2$, $\beta_2=0$, $\beta_3=1$  and $\beta_4=0$,   $\wt{B}(x,y)$ is comparable to the boundary term in the second line of \eqref{e:intro-skst}. 
If $\beta_1=\delta/2$, 
$\beta_2=(\delta/2)(1-\gamma)$ 
($0<\gamma <1$) and $\beta_3   =\beta_4=0$, then $\wt{B}(x,y)$ is comparable to the boundary term in the third line of \eqref{e:intro-skst}. 

 Finally, we outline an example where $\beta_4=1$. For more details we refer the reader to \cite[Example 7.3]{CKSV21} (where different notation is used). Let $D=\{x\in \R^d: \, x_d>0\}$ be the upper half space in $\R^d$, $\delta\in (0,2)$ and $q\in [\delta-1, \delta)\cap (0, \delta)$. The underlying process $Z^D$ corresponds to a Feynman-Kac semigroup of the part process in $D$ of the reflected $\delta$-stable process in $\overline{D}$ via a multiplicative functional involving parameter $q$. The process $Z^D$ is subordinated by an independent $\gamma/2$ subordinator, $\delta \gamma=2\alpha$, to obtain the process $Y^D$. 
In the case $q=\frac12\delta(1-\frac\gamma2)(=\frac12(\delta-\alpha))$, the 
jump kernel $J(x,y)$ of $Y^D$ is comparable to 
$|x-y|^{-\alpha-d}\wt{B}(x,y)$ 
where in $\wt{B}(x,y)$ we have $\beta_1=\beta_2=q$, $\beta_3=0$ and $\beta_4=1$.

\begin{lemma}\label{l:B-tilde-assumptions} 
The boundary term $\wt{B}(x,y)$ satisfies Assumptions \textbf{(B1)-(B6)}.
\end{lemma}
\pf Assumptions \textbf{(B1)} and \textbf{(B2)} are clear.

\noindent \textbf{(B3)} Let $a\in (0,1)$ and assume that $\delta_D(x)\wedge \delta_D(y)> a|x-y|$. Then the quantities in the  first two parentheses in $\wt{B}(x,y)$ are larger than $a$,  and 
those 
in $L(x,y)$ and $K(x,y)$ are equal to $\log 2$. 
Therefore, $\wt{B}(x,y)\ge (\log 2)^{\beta_3  + \beta_4 }a^{\beta_1+\beta_2}=:C_3(a)$. 

\noindent \textbf{(B4)} Suppose $\delta_D(x)\wedge \delta_D(y)\ge |x-y|$. Then clearly 
 $\wt{B}(x,y)=(\log 2)^{\beta_3+\beta_4} $, so $(\log 2)^{\beta_3+\beta_4}-\wt{B}(x,y)=0$, 
and 
the assumption holds with any $\theta>0$, in particular, we can take $\theta=1$.

\noindent \textbf{(B5)} Let $\varepsilon\in (0,1)$, $x_0\in D$, and $r>0$ so that $B(x_0, (1+\varepsilon)r)\subset D$, and let $x_1,x_2\in B(x_0,r)$, $z\in D\setminus B(x_0, (1+\varepsilon)r)$.  It is easy to see that
$$
\frac{1}{3}\delta_D(x_2)\le \left(1+\frac{2}{1+\varepsilon}\right)^{-1}\delta_D(x_2)\le \delta_D(x_1)\le \left(1+\frac{2}{1+\varepsilon}\right)\delta_D(x_2) \le 3 \delta_D(x_2)\, ,
$$
and similarly
$
3^{-1} |x_2-z|\le |x_1-z|\le 3 |x_2-z|\, .
$
Now it is straightforward to obtain that 
$$
\left(\frac{\delta_D(x_1)\wedge\delta_D(z)}{|x_1-z|}\wedge 1\right)\le 9 \left(\frac{\delta_D(x_2)\wedge\delta_D(z)}{|x_2-z|}\wedge 1\right)
$$
 and the analogous estimate with the minimum replaced by the maximum, 
$L(x_1,z)\le 9 L(x_2,z)$  and $K(x_1,z)\le 9 K(x_2,z)$.  Hence
$
\wt{B}(x_1, z)\le 9^{\beta_1+\beta_2+\beta_3  +\beta_4 } \wt{B}(x_2, z)\, .
$

\noindent \textbf{(B6)} Let $M\ge 1$ and let $x,y,z\in D$ satisfy $\delta_D(x)\le \delta_D(z)$ and $|y-z|\le M|y-x|$. By Lemma \ref{l:2-assumptions}, there is a 
constant $c=c(\beta_1, \beta_2, \beta_3, \beta_4)>0$ such that 
\begin{align*}
&\wt{B}(x,y)=\sA(\delta_D(x), \delta_D(y), |x-y|)\le c \sA(\delta_D(z), \delta_D(y), M^{-1}|z-y|)\\
&= c
\Big(\frac{\delta_D(z)\wedge \delta_D(y)}{M^{-1}|z-y|}\wedge 1\Big)^{\beta_1}\Big(\frac{\delta_D(z)\vee \delta_D(y)}{M^{-1}|z-y|}\wedge 1\Big)^{\beta_2}
\log\Big(1+\frac{(\delta_D(z)\vee \delta_D(y))\wedge (M^{-1}|z-y|)}{\delta_D(z)\wedge \delta_D(y)\wedge (M^{-1}|z-y|)}\Big)^{\beta_3} \\
& \quad  \times \log\Big(1+\frac{M^{-1}|z-y|)}{(\delta_D(z)\vee \delta_D(y))\wedge M^{-1}|z-y|)}\Big)^{\beta_4} \\ 
&\le  c M^{\beta_1+\beta_2}\Big(\frac{\delta_D(z)\wedge \delta_D(y)}{|z-y|}\wedge 1\Big)^{\beta_1}\Big(\frac{\delta_D(z)\vee \delta_D(y)}{|z-y|}\wedge 1\Big)^{\beta_2}     {\log\Big(1+\frac{(\delta_D(z)\vee \delta_D(y))\wedge (|z-y|)}{\delta_D(z)\wedge \delta_D(y)\wedge (|z-y|)}\Big)^{\beta_3}} \\
& \quad  \times \log\Big(1+\frac{|z-y|)}{(\delta_D(z)\vee \delta_D(y))\wedge |z-y|)}\Big)^{\beta_4}  \\ 
&=c M^{\beta_1+\beta_2}\wt{B}(z,y)\, .
\end{align*}
The second inequality above follows from the fact that, when $s \ge u$,   $t \to  \frac{s\wedge t}{u\wedge t}$ 
 and $t\mapsto \frac{t}{s\wedge t} $ are  increasing on $(0, \infty)$. 
\qed

Let 
\begin{equation}\label{e:def-of-B-hat}
\wh{B}(x,y):=\left(\frac{(\delta_D(x)\wedge \delta_D(y))^{\beta_1}(\delta_D(x)\vee \delta_D(y))^{\beta_2}}{|x-y|^{\beta_1+\beta_2}}\wedge 1\right) L(x,y)^{\beta_3}  K(x,y)^{\beta_4} \, .
\end{equation}
The next result shows that $\wt{B}$ and $\wh{B}$ are comparable.
\begin{lemma}\label{l:comp} For all $x,y\in D$, 
$
\wh{B}(x,y)\ge \wt{B}(x,y) \ge 2^{-(\beta_1 \vee \beta_2)}\wh{B}(x,y)\, .
$
\end{lemma}
\pf The first inequality is trivial, so we only prove the second. Without loss of generality, we assume that $\delta_D(x)\le \delta_D(y)$. It is clear that the two sides are equal 
when $\delta_D(x)\wedge\delta_D(y) \ge {|x-y|}$ or when 
$\delta_D(x)\vee\delta_D(y )\le {|x-y|}$.
So we assume that
$
{\delta_D(x)}< {|x-y|}<{\delta_D(y)}.
$
Then $\delta_D(y)\le \delta_D(x)+|y-x|\le 2|x-y|$. Thus
\begin{align*}
&\frac{\delta_D(x)^{\beta_1}\delta_D(y)^{\beta_2}}{|x-y|^{\beta_1+\beta_2}}\wedge 1 \le 2^{\beta_2}\left(\frac{\delta_D(x)^{\beta_1}}{|x-y|^{\beta_1}}\wedge 1\right)=
2^{\beta_2}\left(\frac{\delta_D(x)^{\beta_1}}{|x-y|^{\beta_1}}\wedge 1\right)\left(\frac{\delta_D(y)^{\beta_2}}{|x-y|^{\beta_2}}\wedge 1\right).
\end{align*}
\qed

Since $\wh{B}$ is comparable to $\wt{B}$, we immediately obtain that $\wh{B}(x,y)$ also satisfies \textbf{(B1)-(B3}), \textbf{(B5)} and \textbf{(B6)}. If $\delta_D(x)\wedge \delta_D(y)\ge |x-y|$, 
we have that $\wh{B}(x,y)=(\log 2)^{\beta_3  +\beta_4  }$, 
so \textbf{(B4)} also holds with any $\theta>0$. In particular, we can take $\theta=1$. 

We now introduce  yet  another boundary 
term. Let
\begin{equation}\label{e:def-of-B}
\bar{B}(x,y):=\frac{(\delta_D(x)\wedge \delta_D(y))^{\beta_1}(\delta_D(x)\vee \delta_D(y))^{\beta_2}}{|x-y|^{\beta_1+\beta_2}+(\delta_D(x)\wedge \delta_D(y))^{\beta_1}(\delta_D(x)\vee \delta_D(y))^{\beta_2}} L(x,y)^{\beta_3}  K(x,y)^{\beta_4} .
\end{equation}
The elementary inequalities
$
\frac12 \left(\frac{a}{b}\wedge 1\right) \le  \frac{a}{a+b} \le \, \frac{a}{b}\wedge 1, \quad a,b>0\, ,
$ 
imply
\begin{lemma}\label{l:comp-2}
For all $x,y\in D$,
$
2^{-1} \wh{B}(x,y)\le\bar{B}  (x,y)\le \wh{B}(x,y)\, .
$
\end{lemma}
Thus $\bar{B} $ is comparable to $\wt{B}$ and $\wh{B}$. Since Assumptions \textbf{(B1)-(B3)}, \textbf{(B5)} and \textbf{(B6)} depend only on comparability, we conclude that they are also satisfied by the boundary term 
$\bar{B}(x,y)$. Let $x,y\in D$ such that $\delta_D(x)\wedge \delta_D(y)\ge |x-y|$. Then 
 $L(x,y)=K(x,y)=\log 2$,  so
$$
(\log 2)^{\beta_3 + \beta_4 }-\bar{B}  (x,y)=\frac{|x-y|^{\beta_1+\beta_2}}{|x-y|^{\beta_1+\beta_2}+(\delta_D(x)\wedge \delta_D(y))^{\beta_1}(\delta_D(x)\vee \delta_D(y))^{\beta_2}}\le \frac{|x-y|^{\beta_1+\beta_2}}{(\delta_D(x)\wedge \delta_D(y))^{\beta_1+\beta_2}}.
$$
Hence \textbf{(B4)} is satisfied with $\theta=\beta_1+\beta_2$. 

Note that when $\beta_1=\beta_2=\beta_3 = \beta_4  =0$, $\wt B(x, y)$, $\wh B(x, y)$ and 
$\bar B(x, y)$ 
reduce to 1, so the trivial case $\sB(x, y)\equiv 1$ satisfies all our assumptions.

Now we consider a subordinate killed L\'evy process in an arbitrary open subset $D\subset \R^d$ and show that, under some mild assumptions,
 its boundary term satisfies \textbf{(B1)}-\textbf{(B6)}. 
 
Let $Z$ be   a subordinate Brownian motion $\R^d$, $d \ge 2$, 
via a subordinator whose Laplace exponent $\varphi$  is a complete Bernstein function 
satisfying the following global weak scaling condition: 
There exist $a_1, a_2>0$ and $0<\delta_3\le \delta_4 <1$ such that
\begin{equation}\label{e:weak-scaling-phi}
a_1 \left(\frac{R}{r}\right)^{\delta_3}\le \frac{\varphi(R)}{\varphi(r)} \le a_2 \left(\frac{R}{r}\right)^{\delta_4}\, , 
\qquad 0<r\le R<\infty\, .
\end{equation}
Then it has a continuous transition density $p^{Z}(t, |x-y|)$.
Let $Z^{D}$ be the subprocess  of $Z$ killed upon leaving $D$.  
Then  the transition density  $p^{Z, D}(t, x, y)$ of $Z^{D}$  has the form 
\begin{align*}p^{Z, D}(s, x, y)=p^{Z, D}(s, |x -y|)-\E_x \left[ p^{Z}\big(s - \tau^{Z}_D,|Z_{\tau^{Z}_D}-y|\big) : \tau^{Z}_D < s\right],\end{align*}
where $\tau^Z_D:=\inf\{t>0:\, Z_t\notin D\}$.
 
Let $T=(T_t)_{t\ge 0}$ be another subordinator, independent of $Z$, with Laplace exponent $\psi$ and 
L\'evy measure $\nu$. We assume that $\psi$ is a complete Bernstein function (so that  $\nu$ has complete monotone density $\nu(x)$)
satisfying the following global weak scaling condition:
There exist $b_1, b_2>0$ and $0<\gamma_1\le \gamma_2 <1$ such that
\begin{equation}\label{e:weak-scaling-psi}
b_1 \left(\frac{R}{r}\right)^{\gamma_1}\le \frac{\psi(R)}{\psi(r)} \le b_2 \left(\frac{R}{r}\right)^{\gamma_2}\, , 
\qquad 0<r\le R<\infty\, .
\end{equation}
So $\psi$ also has  no drift.

Let $Y_t$ and $X_t$ be the processes defined by $Y_t:=Z^{D}_{T_t}$  and $X_t:=Z_{T_t}$, respectively.  
 The jump kernel $J^{ D}(x, y)$ of $Y$ and the jump kernel $J(|x-y|)$ of $X$ 
 are given respectively by
\begin{align*}
J^{ D}(x, y)=\int^\infty_0p^{Z, D}(s, x, y) \nu(s) ds \text{ and } J(|x-y|)=\int^\infty_0p^{Z}(s, |x- y|) \nu(s) ds\, .
\end{align*}
It is clear that $J^{ D}(x, y)\le J(|x-y|)$. 
By \cite[Theorem 3.4]{KSV14}, there exists $c\ge 1$ such that
\begin{equation}\label{e:JD-estimate}
c^{-1} \frac{(\psi \circ \varphi)(|x-y|^{-2})}{|x-y|^{d}}\le 
J(|x-y|)\le 
c \frac{(\psi \circ \varphi)(|x-y|^{-2})}{|x-y|^{d}}\, ,\qquad x\neq y.
\end{equation}
Let $\sB(x, y)= J^D(x,y)/j(|x-y|)$ for $x\neq y$ and $\sB(x, x)=1$. Then $\sB(x, y)\le 1$.

We first check that $\sB$ satisfies  \textbf{(B4)}.
Assume that $\delta_D(y) \ge  \delta_D(x) \ge  |x-y|$.  Then
\begin{align*}
&j(|x- y|)-J^{ D}(x, y)
=\int^\infty_0\left(p^{Z}(s, |x- y|)-p^{Z, D}(s, x, y)\right)
\nu(s)ds\\
&=\int^\infty_0\E_x \left[ p^{Z}\big(s - \tau^{Z}_D, |Z_{\tau^{Z}_D}-y|\big) : \tau^{Z}_D < s\right]\nu(s)ds=
\E_x \int^\infty_{ \tau^{Z}_D} 
p^{Z}\big(s - \tau^{Z}_D,|Z_{\tau^{Z}_D}-y|\big) 
\nu(s)ds
\nn\\
&=
\E_x \int^\infty_{0} 
p^{Z}\big(v,|Z_{\tau^{Z}_D}-y|\big) 
\nu(v+\tau^{Z}_D)ds
\le
\E_x \int^\infty_{0} 
p^{Z}\big(v,|Z_{\tau^{Z}_D}-y|\big) 
\nu(v)ds=
\E_x \left [ J\big(|Z_{\tau^{Z}_D}- y|\big)\right].
\end{align*}
Since $|Z_{\tau^{Z}_D}- y| \ge \delta_D(y)$ and $j$ is monotone, we have 
$j(|x- y|)-J^{ D}(x, y) \le  j(\delta_D(y))$.
Hence by  \eqref{e:weak-scaling-phi}, \eqref{e:weak-scaling-psi} and \eqref{e:JD-estimate},
\begin{align}
&\sB(x, x)-\sB(x, y)=
1-\sB(x, y)\le \frac{j(\delta_D(y))}{j(|x-y|)}\nonumber\\
&\le c \left(\frac{|x-y|}{ \delta_D(y)}\right)^{d}
\frac{(\psi \circ \varphi)(\delta_D(y)^{-2})}{(\psi \circ \varphi)(|x-y|^{-2})}
\le c \left(\frac{|x-y|}{ \delta_D(y)}\right)^{d+2 \gamma_1\delta_3}.
\label{e:1-Bnew}
\end{align}
Thus \textbf{(B4)} holds with $\theta=d+2 \gamma_1\delta_3>1$.
If $D$ is a bounded $\kappa$-fat open set, then it follows from \cite[Propositions 4.5 and 5.3, and the proof of Theorem 5.5]{KSV18-b} that $\sB(x,y)$  satisfies assumptions \textbf{(B1)}-\textbf{(B3)}, 
\textbf{(B5)} and \textbf{(B6)} as well. 

If $Z$ is instead a Brownian motion in $\R^d$ with $d \ge 3$, then by following the same argument  one can show that 
the corresponding boundary term $\sB(x,y)$  satisfies  \textbf{(B4)} with $\theta=d+2 \gamma_1$.  
If $D$ is either a bounded $C^{1,1}$ domain, or a domain consisting of all the points above the graph of a bounded globally $C^{1,1}$ function, or
$D$ is a $C^{1,1}$
 domain with compact complement, it follows from \cite{KSV18-a} that 
the corresponding boundary term $\sB(x,y)$  satisfies assumptions \textbf{(B1)}-\textbf{(B3)}, \textbf{(B5)} and \textbf{(B6)} as well.

In the remaining part of this section we assume that $D=\R^d_+$, 
$Z=Z^{(\delta)}$ is a $\delta$-stable process in $\R^d$ and $T_t$ is an independent $(\gamma/2)$-stable subordinator where $\delta\in (0, 2]$ and $\gamma\in (0, 2)$. Let $\alpha=\delta\gamma/2$ and denote by $J^{(\alpha), \R^d_+}(x, y)$ the jump kernel of $Y_t:=Z^{(\delta), \R^d_+}_{T_t}$
It follows from \eqref{e:intro-skbm} and \eqref{e:intro-skst} that 
the boundary term $\sB(x, y)= J^D(x,y)/j(|x-y|)$ satisfies the assumption 
\textbf{(B7)} with $\beta_4=0$.
Let $p^{(\delta), \R^d_+}(t, x, y)$ be the transition density of $Z^{(\delta), \R^d_+}$.
By using the scaling property in the third equality and the change of variables $s=a^{\delta}t$ in the fourth equality below, we get
\begin{align*}
J^{(\alpha), \R^d_+}(ax, ay)
&=c(\gamma)\int^\infty_0p^{(\delta), \R^d_+}(a^{\delta}(a^{-\delta}s), ax, ay)s^{-1-\frac{\gamma}2}ds\\
&=c(\gamma)\int^\infty_0a^{-d}p^{(\delta), \R^d_+}(a^{-\delta}s, x, y)s^{-1-\frac{\gamma}2}ds\\
&=a^{-(d+\alpha)}c(\gamma)\int^\infty_0p^{(\delta), \R^d_+}(t, x, y)t^{-1-\frac{\gamma}2}dt=a^{-(d+\alpha)}J^{(\alpha), \R^d_+}(x, y).
\end{align*}
This proves that $\sB(x,y)$ is homogeneous. The second part of  \textbf{(B8)} follows from translation  
invariance property of the transition density: $p^{(\delta), \R^d_+}(t, x+(\wt{z},0), y+(\wt{z},0))=p^{(\delta), \R^d_+}(t, x, y)$ for all $x,y\in \R^d_+$ and all $\wt{z}\in \R^{d-1}$.

\section{Interior results}\label{s:interior}

Let $D\subset \R^d$ be an open set, $D\neq \R^d$.
From now on we always assume that 
\textbf{(B1)}-\textbf{(B3)} and 
\eqref{e:assumption-on-j}--\eqref{e:bound-on-kappa}   hold.
Recall that $\delta_D(x):=\mathrm{dist}(x, \partial D)$ and that we have defined the kernel $J^D(x,y)$ as 
$
J^D(x,y)=\sB(x,y) j(|x-y|).
$
Define
$$
\EE^D(u,v):=\frac12 \int_D\int_D (u(x)-u(y))(v(x)-v(y))J^D(x,y)\, dy\, dx.
$$
By Fatou's lemma, $(C_c^{\infty}(D), \EE^D)$ is closable in $L^2(D, dx)$. Let 
$\FF^D$ 
be the closure of $C_c^{\infty}(D)$ under $\EE^D_1:=\EE^D+(\cdot, \cdot)_{L^2(D,dx)}$. Then $(\FF^D, \EE^D)$ is a regular Dirichlet form on $L^2(D, dx)$.

Recall that  $\kappa:D\to [0,\infty)$ is a Borel function satisfying \eqref{e:bound-on-kappa}, namely
$
\kappa(x)\le C_1 {\Phi(\delta_D(x))^{-1}}$, $ x\in D$,
for some $C_1>0$. 
This implies that $\kappa$ is locally bounded in $D$.
Set
$$
\EE^{D, \kappa}(u,v):=\EE^D(u,v)+\int_D u(x)v(x)\kappa(x)\, dx\, .
$$
Since $\kappa$ is locally bounded, the measure $\kappa(x)dx$ is a positive Radon measure  charging no set of zero capacity. Let $\FF^{D,\kappa}:=\wt{\FF^D}\cap L^2(D, \kappa(x)dx)$, where $\wt{\FF^D}$ is the family of all quasi-continuous functions on $\FF^D$. By \cite[Theorems 6.1.1 and 6.1.2]{FOT}, $(\FF^{D,\kappa}, \EE^{D, \kappa})$ is a regular Dirichlet form on $L^2(D, dx)$ having $C_c^{\infty}(D)$ as a special standard core.
Let $(Y_t^{D,\kappa}, \P_x)$, $t\ge 0$, $x\in D\setminus \NN$, be the associated Hunt 
process with lifetime $\zeta$, where $\NN$ is an exceptional set.
We add a cemetery point $\partial$ to the state space $D$ and define $Y_t^{D,\kappa}=\partial$
for $t\ge \zeta$. We will write $D_{\partial}=D\cup \{\partial\}$.

For an open set $U\subset \overline{U}\subset D$, let $Y^{D,\kappa, U}$ be the process $Y^{D,\kappa}$ killed upon exiting $U$, that is, the part of the process $Y^{D,\kappa}$ in $U$. The Dirichlet form of $Y^{D,\kappa, U}$ is equal to $(\FF^{D,\kappa}_U, \EE^{D, \kappa})$, where $\FF^{D,\kappa}_U=\{u\in \FF^{D,\kappa}: \, u=0 \textrm{ q.e.~on  } D\setminus U\}$. 
For $u,v\in \FF^{D,\kappa}_U$,
\begin{equation}\label{e:EDkappa}
\EE^{D, \kappa}(u,v)=\frac{1}{2}\int_U\int_U (u(x)-u(y))(v(x)-v(y))J^D(x,y)\, dy\, dx +\int_U u(x)u(y)\kappa_U (x)\, dx\, ,
\end{equation}
where
\begin{equation}\label{e:kappaU}
\kappa_U(x)=\int_{D\setminus U} J^D(x,y)\, dy +\kappa(x)\, , \quad x\in U\, .
\end{equation}
Moreover, by \cite[Theorem 3.3.9]{CF}, since $C_c^{\infty}(D)$ is a special standard core of $(\FF^{D,\kappa}, \EE^{D, \kappa})$, the set $\{u\in C_c^{\infty}(D):\, \mathrm{supp}(u)\subset U\}=C_c^{\infty}(U)$ is a core of $(\FF^{D,\kappa}_U, \EE^{D, \kappa})$.

\subsection{Regularization of the process}\label{ss:regularization} 
In this subsection we show that one can remove the exceptional set $\NN$ and start the process $Y^{D, \kappa}$ from every point $x\in D$.
For this purpose, we 
will use the process $Z$ on $\R^d$, with jump kernel $J_\gamma$
defined below, as a tool. The process $Z$ can start from every point in $\R^d$.

Let $U$ be a relatively compact $C^{1,1}$ open subset of $D$. 
Note that for $x,y\in U$, $\delta_D(x)\wedge \delta_D(y)\ge \mathrm{dist}(U, \partial D)$, so with $a:=\mathrm{dist}(U, \partial D)/\mathrm{diam}(U)$ we have that $\delta_D(x)\wedge \delta_D(y)\ge a|x-y|$.  Therefore, it follows from \textbf{(B2)} and \textbf{(B3)} that there exists a constant $c_1=c_1(U)\ge 1$ such that $c^{-1}\le \sB(x,y)\le c$ for all $x,y\in U$. Together with \eqref{e:assumption-on-j} this implies that there 
exist  $c_2>0$  and $c_3>0$ such that
\begin{equation}\label{e:JD-in-U}
\frac{c_2}{|x-y|^d\Phi(|x-y|)}\le J^D(x,y)\le \frac{c_3}{|x-y|^d\Phi(|x-y|)}\, ,\quad x,y\in U\, .
\end{equation} 
For $\gamma>0$ define a kernel $J_{\gamma}(x,y)$ on $\R^d\times \R^d$ by $J_{\gamma}(x,y)=J^D(x,y)$ for $x,y\in U$, and $J_{\gamma}(x,y)= \gamma j(|x-y|)$ otherwise. 
By\eqref{e:assumption-on-j} 
and \eqref{e:JD-in-U},
there exist $c_4>0$ and $c_5>0$  such that
$$
\frac{c_4}{|x-y|^d\Phi(|x-y|)}\le J_{\gamma}(x,y)\le \frac{c_5}{|x-y|^d\Phi(|x-y|)}\, ,\quad x,y\in \R^d\, .
$$
It is now straightforward to check that all conditions of \cite[Theorem 1.2]{CK08} (as well as the geometric condition of the paper) are satisfied.  
For $u\in L^2(\R^d, dx)$, define
$$
\CC(u,u):=\frac12 \int_{\R^d}\int_{\R^d} (u(x)-u(y))^2 J_{\gamma}(x,y)\, dx\, dy \text{ and }
\DD(\CC):=\{u\in L^2(\R^d):\, \CC(u,u)<\infty\}\, .
$$
Note that $C_c^{\infty}(\R^d)$ is a special standard core of $ \DD(\CC)$.

Let
$$
\wt{q}(t,x,y):=\Phi(t)^{-d}\wedge \frac{t}{|x-y|^d\Phi(|x-y|)}\, , \quad t>0, \ x,y\in \R^d.
$$
It follows from 
\cite{CK08} that there exists  a conservative Feller and strongly Feller process $Z$ associated with $( \DD(\CC), \CC)$ that can start from every point in $\R^d$. Moreover, the process $Z$ has a continuous transition density function $p(t,x,y)$  on $(0, \infty)  \times \R^d\times \R^d$ (with respect to the Lebesgue measure) that satisfies the following estimate: There exists  $c_6\ge 1$ such that
$$
c_6^{-1}\wt{q}(t,x,y)\le p(t,x,y)\le c_6 \wt{q}(t,x,y)\, ,\quad t>0, \ x,y\in \R^d\, .
$$

Denote the part of the process $Z$ killed upon exiting $U$ by $Z^U$. The Dirichlet form of $Z^U$ is $(\DD_U(\CC), \CC)$, where for $u,v\in \DD_U(\CC)$,
\begin{eqnarray*}
\CC(u,v)&=&\frac{1}{2}\int_U \int_U (u(x)-u(y))(v(x)-v(y))J_{\gamma}(x,y)\, dy\, dx +\int_U u(x)v(x)\kappa^{Z}_U(x)\, dx\\
&=&\frac{1}{2}\int_U \int_U (u(x)-u(y))(v(x)-v(y))J^D(x,y)\, dy\, dx +\int_U u(x)v(x)\kappa^{Z}_U(x)\, dx
\end{eqnarray*}
with
\begin{equation}\label{e:kappa-Z-new}
\kappa^{Z}_U(x)=\int_{\R^d\setminus U}J_{\gamma}(x,y)\, dy=\gamma \int_{\R^d\setminus U}j(|x-y|)\, dy\, ,\quad x\in U\, ,
\end{equation}
and $\DD_U(\CC)=\{u\in \DD(\CC):\, u=0 \textrm{ q.e.~on  } \R^d\setminus U\}$. Again, by \cite[Theorem 3.3.9]{CF}, $C_c^{\infty}(U)$ is a core of $(\DD_U(\CC), \CC)$.

Recall $\delta_U=\mathrm{dist}(U, \partial D)$ and $d_U=\mathrm{diam}(U)$. For all $x,y\in U$ we have that $\delta_D(x)\wedge \delta_D(y)\ge (\delta_U/d_U)|x-y|$, hence by \textbf{(B3)} we have $\sB(x,y)\ge c_7$ 
(where $c_7=C_3(2^{-1} \wedge (\delta_U/d_U))$)
implying together with \textbf{(B2)}  that
\begin{equation}\label{e:equal-domains-1}
c_7 j(|x-y|)\le J^D(x,y)=\sB(x,y)j(|x-y|)\le C_2 j(|x-y|)\, ,\quad x,y\in U\, .
\end{equation}
Next, let $V$ be the $\delta_U/2$-neighborhood of $U$, that is $V:=\{x\in D:\, \mathrm{dist}(x, U)<\delta_U/2\}$. Then
$$
\kappa_U(x)=\int_{D\setminus V}J^D(x,y)\, dy+\int_{V\setminus U}J^D(x,y)\, dy +\kappa(x)\, .
$$
Similarly as above we conclude that $c_8 j(|x-y|)\le J^D(x,y)\le C_2 j(|x-y|)$ for all $x,y\in V$ with $c_8:=C_3((\delta_U/2)/(d_U+\delta_U/2))>0$.
Moreover, $\sup_{x\in U}\int_{D\setminus V}J^D(x,y)\, dy =:c_9\ < \infty$ (this can be shown by splitting the integral into two parts -- over 
$(D\setminus V)\cap B(0,R_0)$ and $(D\setminus V)\cap B(0,R_0)^c$ for appropriate $R_0$)
and  $\|\kappa_{|U}\|=:c_{10}<\infty$. Therefore
$$
c_8\int_{V\setminus U}j(|x-y|)\, dy \le \kappa_U(x)\le c_9 +C_2 \int_{V\setminus U}j(|x-y|)dy+c_{10}\, ,\quad x\in U\, .
$$
Since 
$$
\inf_{x\in U}\int_{V\setminus U}j(|x-y|)\, dy \ge |V\setminus U|j(\mathrm{diam}(V))=:c_{11}>0\, ,
$$
we conclude that
$$
c_8\int_{V\setminus U}j(|x-y|)\, dy \le \kappa_U(x)\le c_{12} \int_{V\setminus U}j(|x-y|)\, dy\, ,\quad x\in U\, .
$$
Further, since
$$
\kappa_U^Z(x)=\gamma \left(\int_{\R^d\setminus V}  j(|x-y|)\, dy+\int_{V\setminus U} j(|x-y|)\, dy\right)\, ,\quad x\in U\, 
$$
and $\sup_{x\in U}\int_{\R^d\setminus V}j(|x-y|)\, dy=:c_{13}<\infty$, we see that there is a constant $c_{14}>0$ such 
$$
\int_{V\setminus U}j(|x-y|)\, dy \le \gamma^{-1}\kappa_U^Z(x)\le c_{14} \int_{V\setminus U}j(|x-y|)\, dy\, ,\quad x\in U\, .
$$
It follows that
$
c_{12}^{-1}\kappa_U(x)\le \gamma^{-1}\kappa_U^Z(x)\le c_{14}c_8^{-1}\kappa_U(x)$ for all $x\in U$
with positive constants $c_8, c_{12}, c_{14}$ not depending on $\gamma$. Now we choose $\gamma>0$ so small that $\gamma c_{14}c_8^{-1}\le 1$. With this choice we get that $\kappa_U^Z(x)\le \kappa_U(x)$ for all $x\in U$. In particular, with $c_{15}:=\gamma c_{12}^{-1}$ we see that
\begin{equation}\label{e:kappa-U-kappa-Z-comparable}
c_{15}\kappa_U(x)\le \kappa_U^Z(x)\le \kappa_U(x)\, , \qquad x\in U\, .
\end{equation}

It follows that for $u\in C_c^{\infty}(U)$,
\begin{eqnarray*}
\EE^{D,\kappa}_1(u,u)&=&\EE^{D,\kappa}(u,u)+\int_U u(x)^2 \, dx\\
&=&\frac{1}{2}\int_U \int_U (u(x)-u(y))^2 J^D(x,y)\, dy\, dx +\int_U u(x)^2 \kappa_U(x)\, dx+\int_U u(x)^2 dx \\
&\asymp & \frac{1}{2}\int_U \int_U (u(x)-u(y))^2 J_{\gamma}(x,y)\, dy\, dx +\int_U u(x)^2 \kappa_U^Z(x)\, dx+\int_U u(x)^2 dx\\
& =&\CC(u,u) +\int_U u(x)^2 dx=\CC_1(u,u)\, .
\end{eqnarray*}
Since $C_c^{\infty}(U)$ is a core of both $(\FF^{D,\kappa}_U, \EE^{D, \kappa})$ and $(\CC, \DD_U(\CC))$, we conclude that $\FF^{D,\kappa}_U=\DD_U(\CC)$.

We now define $\wt{\kappa}:U\to \R$ by 
$
\wt{\kappa}(x):=\kappa_U(x)-\kappa_U^Z(x)$,  $x\in U.
$
By the choice of $\gamma$ we have that $\wt{\kappa}\ge 0$. On the other hand, since
 $\kappa$ is bounded in $U$ and $U$ is $C^{1,1}$, 
it follows from \eqref{e:kappa-Z-new} and \eqref{e:kappa-U-kappa-Z-comparable} that there is a constant $c_{16}>0$ such that
\begin{equation}\label{e:estimate-kappaU}
\wt{\kappa}(x)\le \kappa_U(x)\le  \frac{c_{16}}{\Phi( \delta_U(x))}\, , \quad x\in U.
\end{equation}
Let $\mu(dx)=\wt{\kappa}(x)\, dx$ be a measure on $U$. For $t>0$ and $a\ge 0$ define
$$
N^{U, \mu}_a(t):=\sup_{x\in \R^d}\int_0^t \int_{z\in U: \delta_U(z)>a\Phi^{-1}(t)}\wt{q}(s,x,z)\mu(dz)\, ds\, .
$$
By the definition of $\wt{q}$ and 
\eqref{e:estimate-kappaU} one can  check that $\sup_{t<1}N^{U, \mu}_a(t)<\infty$ and $\lim_{t\to 0} N^{V,\mu}_0(t)=0$ for every relatively compact open 
set $V\subset U$, that is, $\mu\in \textbf{K}_1(U)$ in the notation of \cite[Definition 2.12]{CKSV19}. 

Let $A_t:=\int_0^t \wt{\kappa} (Z^U_s)\, ds$. Then $(A_t)_{t\ge 0}$ is a positive continuous additive functional of $Z^U$ in the strict sense (i.e.~without an exceptional set) with Revuz measure $\wt{\kappa} (x)dx$. For any non-negative Borel function $f$ on $U$ let
$$
T^{U, \wt{\kappa} }_t f(x):=\E_x[\exp(-A_t)f(Z^U_t)]\, ,\qquad t>0, x\in U\, ,
$$
be the Feynman-Kac semigroup of $Z^U$ associated with $\wt{\kappa} (x)dx$. By \cite[Propostion 2.14]{CKSV19}, the Hunt process $Z^{U, \wt{\kappa} }$ on $U$ corresponding to the transition semigroup $(T^{U,\wt{\kappa} }_t)_{t\ge 0}$ has a transition density $q^U(t,x,y)$ (with respect to the Lebesgue measure) such that $q^U(t,x,y)\le c_{17} \wt{q}(t,x,y)$ for $t<1$. Further, $(t,y)\mapsto q^U(t,x,y)$ is continuous for each $x\in U$.

According to \cite[Theorem 6.1.2]{FOT}, the Dirichlet form $\CC^{U, \wt{\kappa} }$ corresponding to $T^{U, \wt{\kappa} }_t$ is regular and is given by
\begin{align*}
\CC^{U,\wt{\kappa} }(u,v)
=\frac{1}{2}\int_U\int_U (u(x)-u(y))(v(x)-v(y)) J^D(x,y) dy  dx+\int_U u(x)v(x)\kappa_U(x)\, dx
\end{align*}
with the domain
$
\DD^{\wt{\kappa}}_U=\DD_U(\CC)\cap L^2(U, \wt{\kappa}(x)dx)\, .
$
Since $(\DD^{\wt{\kappa}}_U, \CC^{U,\wt{\kappa}})$ is regular, the set $\DD^{\wt{\kappa}}_U\cap C_c(U)=\DD_U(\CC)\cap C_c(U)$ is its core.
By comparing with \eqref{e:EDkappa} we see that 
$$
\EE^{D, \kappa}(u,v)=\CC^{U, \wt{\kappa}}(u,v)\, , \qquad u,v\in C_c^{\infty}(U)\, .
$$
 Now we have to argue that the Dirichlet spaces $(\FF^{D,\kappa}_U, \EE^{D, \kappa})$ and $(\DD^{\wt{\kappa}}_U, \CC^{U, \wt{\kappa}})$ are equal. We know that $C_c^{\infty}(U)$ is a core for $\EE^{D, \kappa}$. 
One can easily check that this should be also true for $\CC^{U, \wt{\kappa}}$.
 Further,  $C_c^{\infty}(U)\subset  C_c(U)\cap \{u\in L^2(U, dx):\, \CC^U(u,u)<\infty\}$ (which is a core). Clearly, $C_c^{\infty}(U)$ is dense in $C_c(U)$ with uniform norm. 
It is easy to see that  $C_c^{\infty}(U)$ is dense in $C_c(U)\cap \{u\in L^2(U, dx):\, \CC^U(u,u)<\infty\}$ with $\CC^{U,\kappa_U}_1$ norm.
Thus the process $Z^{U, \wt{\kappa}}$ coincides with $Y^{D, \kappa, U}$. We summarize the above discussion in the following lemma.

\begin{lemma}\label{l:YDkappaU}
The process $Y^{D, \kappa, U}$ can be refined to start from every point in $U$. Moreover, it is strongly Feller.
\end{lemma}

\begin{prop}\label{p:YDregularized}
The process $Y^{D, \kappa}$ can be refined to start from every point in $D$.
\end{prop}
\pf Let $(U_n)_{n\ge 1}$ be a sequence of bounded smooth open sets such that
$U_n\subset \overline{U_n}\subset U_{n+1}$ and $\cup_{n\ge 1}U_n=D$. For each $n\ge 1$, let $Y^{(n)}=Y^{D,\kappa, U_n}$ and let $\zeta^{(n)}$ be the
lifetime of $Y^{(n)}$. The Dirichlet form of $Y^{(n)}$ is
$(\FF^{D,\kappa}_{U_n}, \EE^{D, \kappa})$, where $\FF^{D,\kappa}_{U_n}=\{u\in \FF^{D,\kappa}: \, u=0 \textrm{ q.e.~on  } D\setminus U_n\}$. 
It follows from  Lemma \ref{l:YDkappaU} that 
each process $Y^{(n)}$ can start from every point in $U_n$.
These processes are consistent in the sense that $Y^{(n+1)}_t=Y^{(n)}_t$ for $t\in [0, \zeta^{(n)})$.
We define a new process $\widetilde{Y}$ by
$\widetilde{Y}_t=Y^{(1)}_t$ for $t\in [0, \zeta^{(1)})$, and for each $n> 1$, 
$\widetilde{Y}_t=Y^{(n)}_t$ for $t\in [\zeta^{(n-1)}, \zeta^{(n)})$. 
The process $\widetilde{Y}$ can start from every point in $D$.
It is easy to check  that the Dirichlet form of $\widetilde{Y}$ is $(\FF^{D,\kappa}, \EE^{D, \kappa})$. Thus we can assume the process $Y_t^{D,\kappa}$ can start from  every point in $D$. \qed

\subsection{Analysis of the generator}\label{ss:dynkin}
In this subsection we assume that, in addition to \textbf{(B1)}-\textbf{(B3)}, \textbf{(B4)} also holds.
Let 
$$
C_c^2(D; \R^d)=\{f:D\to  \R : \text{ there exists } u \in  C_c^2(\R^d) \text{ such that } u=f \text{ on } D\}
$$
be the space of functions on $D$ that are restrictions of $C_c^2(\R^d)$ functions. Clearly, if $f\in C_c^2(D; \R^d)$ then $f\in C^2_b(D) \cap L^2(D)$. 

We introduce the operator
\begin{eqnarray}\label{e:defn-LBD}
L^\sB f(x)&:=&\textrm{p.v.}\int_{D}(f(y)-f(x))J^{D}(x,y)\, dy- \kappa (x)f(x)\\
&=&\lim_{\varepsilon \to 0} \int_{D, |y-x|> \varepsilon}(f(y)-f(x))J^{D}(x,y)\, dy- \kappa (x)f(x)\, ,\quad x\in D\, ,\nonumber 
\end{eqnarray}
defined for all functions $f:D\to \R$ for which the principal value integral makes sense. 
We will show that this is the case when $f\in C_c^2(D; \R^d)$. Let us start with the following auxiliary result.

\begin{lemma}\label{l:operator-interpretation-B}
There exists a constant $C_{9}>0$ such that for any 
bounded Lipschitz function with Lipschitz constant $L$, 
we have for all  $x\in D$ and $r\in (0,\delta_D(x)]$,
\begin{equation}\label{e:Lemma 3.3}
\int_D \left|f(y)-f(x)\right|j(|y-x|) |\sB(x,x)-\sB(x,y)|\, dy \le C_{9}(\|f\|_{\infty}+r
L)\Phi(r)^{-1}\, .
\end{equation}
\end{lemma}
\pf
Since the proof for the case $\delta_2 <1/2$ is simpler, we give the proof for the case $\delta_2  \ge 1/2$.

Note that
\begin{align*}
&\int_D \left|f(y)-f(x)\right|j(|y- x |) |\sB(x,x)-\sB(x,y)|\, dy \\
&\le 
L\int_{D, |y-x|<r/2} |y-x|j(|y-x|) |\sB(x,x)-\sB(x,y)|\, dy\\
& \quad  +4C_2\|f\|_{\infty}\int_{D, |y-x|\ge r/2}j(|y-x|) \, dy =:I_1+I_2\, ,
\end{align*}
where we used \textbf{(B2)} in the  inequality.
It follows from $\delta_D(x)\ge r$ that, if $|y-x|<r/2$, then $\delta_D(y)>r/2$ and thus
$\delta_D(y)\wedge \delta_D(x)>r/2>|y-x|$. Thus by \textbf{(B4)}, \eqref{e:assumption-on-j} and \eqref{e:H-infty},
\begin{eqnarray*}
I_1&\le  &c_1C_4  
L2^{\theta} r^{-\theta} \int_{|y-x|<r/2}|y-x|^{1+\theta}|y-x|^{-d}\Phi(|y-x|)^{-1}\, dy \\
&\le & c_2 C_4  
L2^{\theta} r^{-\theta}\frac{1}{\Phi(r/2)}\int_0^{r/2}\frac{s^{\theta}\Phi(r/2)}{\Phi(s)}\, ds \\
&\le & a_2 c_2 C_4  
L2^{\theta} r^{-\theta}\frac{(r/2)^{2\delta_2}}{\Phi(r/2)}\int_0^{r/2}s^{\theta-2\delta_2}\, ds \le   
c_3
L r\Phi(r)^{-1}\, .
\end{eqnarray*}
Similarly, using \eqref{e:assumption-on-j} and \eqref{e:H-infty},
$
I_2\le c_4\|f\|_{\infty} \int_{r/2}^{\infty}t^{d-1}j(t)\, dt \le c_5
\|f\|_{\infty}\Phi(r)^{-1}\, .
$
Combining the estimates for $I_1$ and $I_2$ we get the statement of the lemma. \qed

\begin{prop}\label{p:operator-interpretation}
\begin{itemize}
\item[(a)] 
 Let $f\in C_c^2(D; \R^d)$. Then $L^\sB f$ is well defined 
and for all  $x \in D$ and $r>0$,
\begin{eqnarray}\label{e:operator-interpretation}
L^\sB f(x)&=&\sB(x,x) \int_{\R^d}{(u(y)-u(x)-\nabla u(x) {\bf 1}_{\{ |y-x|<r\}}\cdot (y-x))}j(|x-y|)dy\nn\\
& &+\sB(x,x)\int_{\R^d \setminus D}(u(x)-u(y))j(|x-y|)dy\nn\\
& & +\int_{D}{(u(y)-u(x))}j(|x-y|)(\sB(y, x)-\sB(x, x))dy
- \kappa(x)u(x),
\end{eqnarray}
where $u\in C_c^2(\R^d)$ is any function such that $u=f$ on $D$.
\item[(b)] There exists a constant $C_{10}>0$ such that for any 
$f\in C_c^2(D; \R^d)$, any $x\in D$ and any $r\in (0,\delta_D(x)]$ we have
\begin{equation}\label{e:LB-on-u-estimate}
|L^\sB f (x)|\le C_{10}\left(r^2\|\partial^2 u\|_{\infty}+r\|\nabla u\|_{\infty}+\|u\|_{\infty}\right) \Phi(r)^{-1}\, ,
\end{equation}
where $u\in C_c^2(\R^d)$ is any function such that $u=f$ on $D$.

\noindent
\item[(c)]  Let $f\in C_c^2(D)$. Then for all $0<\epsilon_0 \le \mathrm{dist}(\partial D, \mathrm{supp}(f))/2$, 
$$
\|L^\sB f \|_{\infty}\le C_{10}\left(\epsilon^2_0\|\partial^2 f\|_{\infty}+\epsilon_0\|\nabla f\|_{\infty}+\|f\|_{\infty}\right) \Phi(\epsilon_0)^{-1}\, .
$$
\end{itemize}
\end{prop}
\pf (a)
 Let $u\in C_c^2(\R^d)$ be such that $u=f$ on $D$. Fix $x\in D$ and let 
  $\varepsilon<r \wedge (\delta_D(x)/2)$. Then
\begin{align}
&\int_{D, \, |x-y|>\varepsilon} {(f(y)-f(x))}j(|x-y|)\sB(x,y) dy \nn\\
&=\sB(x, x)\int_{D, \, |x-y|>\varepsilon} {(u(y)-u(x))}j(|x-y|)dy\nn\\
&\qquad +\int_{D, \, |x-y|>\varepsilon} {(u(y)-u(x))}j(|x-y|)(\sB(x,y)-\sB(x, x))dy\nn\\
&=\sB(x, x)\int_{|x-y|>\varepsilon} {(u(y)-u(x))}j(|x-y|)dy+
\sB(x, x)\int_{\R^d \setminus D, \,  |x-y|>\varepsilon}{(u(x)-u(y))}j(|x-y|)dy\nn\\
&\qquad +\int_{D, \, |x-y|>\varepsilon}{(u(y)-u(x))}j(|x-y|)(\sB(x,y)-\sB(x, x))dy\nn\\
&=\sB(x, x)\int_{|x-y|>\varepsilon}{(u(y)-u(x)-\nabla u(x) {\bf 1}_{\{ |y-x|<r\}}\cdot (y-x))}j(|x-y|)dy\nn\\
&\qquad +\sB(x, x)\int_{\R^d \setminus D, \, |x-y|>\varepsilon}(u(x)-u(y))j(|x-y|)dy\nn\\
&\qquad +\int_{D, \, |x-y|>\varepsilon}{(u(y)-u(x))}j(|x-y|)(\sB(y, x)-\sB(x, x))dy. \nn
\end{align}
In the last integral above, we have used {\bf(B1)}.
By letting $\varepsilon \to 0$ and using 
Lemma \ref{l:operator-interpretation-B} (with $r$ there being $\delta_D(x)$)
for the third integral, we obtain \eqref{e:operator-interpretation}.

\noindent 
(b) 
Let $u\in C_c^2(\R^d)$ 
be any function such that $u=f$ on $D$. 
Fix $x\in D$ and let $r\in (0, \delta_D(x)]$. Then by part (a), 
\begin{align*}
L^\sB f(x)&= \sB(x,x)\int_{\R^d}\left(u(y)-u(x)-\nabla u(x){\bf 1}_{|y-x|<r}\cdot (y-x)\right)j(|y-x|)\, dy \\
& \quad +\sB(x,x)\int_{\R^d\setminus D} (u(x)-u(y))j(|y- x |)\, dy\\
&\quad  +\int_D(f(y)-f(x))j(|y-x|)(\sB(y,x)-\sB(x,x))\, dy -\kappa(x)f(x)\\
&=:I+II+III+IV\, .
\end{align*}
For $I$, we use
$$
\left|u(y)-u(x)-\nabla u(x){\bf 1}_{|y-x|<r}\cdot (y-x)\right|\le \|\partial^2 u\|_{\infty}|y-x|^2{\bf 1}_{|y-x|\le r}+2\|u\|_{\infty}{\bf 1}_{|y- x |\ge r}
$$
to get 
\begin{align*}
I&\le \sB(x,x)\int_{\R^d}\left(\|\partial^2 u\|_{\infty}|y-x|^2 {\bf 1}_{|y-x|\le r}+2\|u\|_{\infty}{\bf 1}_{|y-x|\ge r}\right)j(|y-x|)\, dy\\
&\le \sB(x,x)c(d) \left(\|\partial^2 u\|_{\infty}\int_0^r t^{d-1}t^2 t^{-d}\Phi(t)^{-1}\, dt+\int_r^{\infty}2\|u\|_{\infty}t^{d-1}t^{-d}\Phi(t)^{-1}\, dt\right)\\
&\le \sB(x,x) c(d) (\|\partial^2 u\|_{\infty}r^2+ 2\|u\|_{\infty})\Phi(r)^{-1}\, .
\end{align*}
For $II$ we use that $\delta_D(x)\ge r$ 
to get 
\begin{align*}
&II\le 2\sB(x,x)\|u\|_{\infty}\int_{B(x, \delta_D(x))}j(|y-x|)\, dy\\
& \le c(d) \sB(x,x)\|u\|_{\infty} \Phi(\delta_D(x))^{-1}\le c(d) \sB(x,x)c(d)\|u\|_{\infty} \Phi(r)^{-1}\, .
\end{align*}
$III$ is estimated in  Lemma \ref{l:operator-interpretation-B}, 
while for $IV$ we use \eqref{e:bound-on-kappa} to get
$$
IV\le C_1\|f\|_{\infty}\Phi(\delta_D(x))^{-1}\le C_1\|f\|_{\infty} \Phi(r)^{-1}\, .
$$

\noindent
(c) Take $x\in D$ such that $\delta_D(x)<\epsilon_0$. Then $B(x,\epsilon_0)\cap D\subset D\setminus \mathrm{supp}(f)$. 
Since $f(x)=0$ and also $f(y)=0$ for $y\in D\cap B(x,\epsilon_0)$, we have that 
$$
| L^\sB f(x) |\le  \int_{|y-x|\ge \epsilon_0} |f(y)|J^D(y,x)\, dy\le \|f\|_{\infty} \int_{|y-x|\ge \epsilon_0}j(|y-x|)\, dy \le c(d)\|f\|_{\infty}\Phi(\epsilon_0)^{-1}\, .
$$
 On the other hand, if $\delta_D(x)\ge \epsilon_0$, then we can take  $r=\epsilon_0$ and a function 
 $u\in C_c^2(D; \R^d)$ such that $u=f$ on $D$ with $\|u\|_\infty=\|f\|_\infty$, 
 $\|\nabla u\|_{\infty}=\|\nabla  f \|_{\infty}$ and  $\|\partial^2 u\|_{\infty}=\|\partial^2  f  \|_{\infty}$  
 in \eqref{e:LB-on-u-estimate} to get
$$
|L^\sB f (x)|\le C_{10}\left(\epsilon^2_0\|\partial^2 f\|_{\infty}+\epsilon_0\|\nabla f\|_{\infty}+\|f\|_{\infty}\right) \Phi(\epsilon_0)^{-1}\, .
$$
\qed

We mention in passing that the right-hand sides of \eqref{e:Lemma 3.3}
and \eqref{e:LB-on-u-estimate} tend to $\infty$ as $r\to 0$. We will apply them
with fixed $r>0$.

\begin{corollary}\label{c:L2-generator}
Let $(A, D(A))$ be the $L^2$-generator of the semigroup $(T_t)_{t\ge0}$ corresponding to ${\mathcal E}^{D, \kappa}$.
Then  $C_c^2(D; \R^d)\subset D(A)$ and 
$A|_{C_c^2(D; \R^d)}=L^\sB|_{C_c^2(D; \R^d)}$. 
\end{corollary}
\pf Since $\kappa$ is locally bounded, it suffices to show that, for $u\in  C_c^2(D; \R^d)$ and  $v \in C_c^2(D)$, 
\begin{align}
\label{e:DFY_1}
& \int_{D}\int_{D}{(u(y)-u(x))(v(y)-v(x))}J^D(x,y)(x,y) \, dy\, dx=
2\int_{D}\big(L^\sB u(x)-\kappa(x)u(x)\big) v(x)\, dx.
\end{align}
Note that, by the Lebesgue dominated convergence theorem and the symmetry of $\sB$, 
\begin{align}
\label{e:DFY_2}
& \int_{D}\int_{D} {(u(y)-u(x))(v(y)-v(x))}j(|x-y|)\sB(x,y) dydx
\nn\\
&=
\lim_{\varepsilon \downarrow
0}  
\int_{D}
\int_{y \in D:|x-y|>\varepsilon} {(u(y)-u(x))(v(y)-v(x))}j(|x-y|)\sB(x,y) dydx\nn\\
&=
2\lim_{\varepsilon \downarrow
0} 
 \int_{D}
\int_{y \in D:|x-y|>\varepsilon} {(u(y)-u(x))}j(|x-y|)\sB(x,y) dy\,  v(x)dx \nn\\
&=2\lim_{\varepsilon \downarrow 0} \int_{\text{supp}(v)}\int_{y \in D:|x-y|>\varepsilon} {(u(y)-u(x))}j(|x-y|)\sB(x,y) dy\,  v(x)dx.
\end{align}
Let $\varepsilon < \eps_0:=\mathrm{dist}(\partial D,\mathrm{supp}(v))/2$. Then 
by the estimates of $I$, $II$ and $III$ in the proof of  Proposition \ref{p:operator-interpretation} (b), we have that
\begin{align*}
&\sup_{x \in \mathrm{supp}(v),\varepsilon < \eps_0}
\left|\int_{y \in D:|x-y|>\varepsilon} {(u(y)-u(x))}j(|x-y|)\sB(x,y) dy  \right| \\
&\qquad \le c_1 \left(\eps_0^2\|\partial^2 u\|_{\infty}+\eps_0\|\nabla u\|_{\infty}+\|u\|_{\infty}\right) \Phi(\eps_0)^{-1}<\infty\, .
\end{align*}
Thus we can use the Lebesgue dominated convergence theorem to conclude that \eqref{e:DFY_1} holds. \qed

Corollary \ref{c:L2-generator} says that $L^\sB$ is the extended generator of the semigroup $(T_t)_{t\ge0}$ corresponding to ${\mathcal E}^{D, \kappa}$.

Consider now the process $Y^{D,\kappa, U}$ -- the process $Y^{D,\kappa}$ killed upon exiting $U\subset \overline{U}\subset D$. 
Denote $L^\sB_U u:=L^{\sB,U} u-\kappa^U(\cdot) u$, where
$$
L^{\sB,U} u(z):=\mathrm{p.v.} \int_U (u(y)-u(z))J^D(y,z)\, dy -\kappa(z)u(z)\, ,\qquad u\in U\, ,
$$
and
$$
\kappa^U(z):=\int_{D\setminus U}J^D(z,y)\, dy\, ,\qquad z\in U\, .
$$
We can write
$$
L^{\sB}_{U} u(z)=\mathrm{p.v.} \int_U (u(y)-u(z))J^D(y,z)\, dy -\kappa_U(z)u(z)\, ,\qquad u\in U\, ,
$$
where $\kappa_U=\kappa+\kappa^U$. If $U$ is a $C^{1,1}$-open set, then $\kappa^U(z)\le c_1\frac{1}{\Phi(\delta_U(z))}$. Since $\kappa$ is bounded in $U$ we have that
$\kappa_U(z)\le c_2\frac{1}{\Phi(\delta_U(z))}$. 
Thus $L^\sB_U$ has the same form in $U$ as $L^\sB$ in $D$. Hence, the analog of Proposition \ref{p:operator-interpretation} is valid. In particular, if $u\in C_c^2(U)$, then $\|L^\sB_U u\|_{\infty}<\infty$.

Further, if $u\in C_c^2(U)\subset C_c^2(D)$, then for $z\in U$, 
\begin{align*}
&L^\sB u(z)=\lim_{\epsilon \to 0}\int_{D, |y-z|<\epsilon}(u(z)-u(y))J^D(y,z)dy -\kappa(z)u(z)\\
&=\lim_{\epsilon \to 0}\int_{U, |y-z|<\epsilon}(u(z)-u(y))J^D(y,z)dy+\lim_{\epsilon \to 0}\int_{D\setminus U, |y-z|<\epsilon}(u(z)-u(y))J^D(y,z)dy-\kappa(z)u(z)\\
&=L^{\sB,U} u(z)-\kappa^U(z)u(z)-\kappa(z)u(z)= L^\sB_U u(z)\, .
\end{align*}

Recall that our basic process is $Y^{D, \kappa}=(Y^{D, \kappa}_t, \P_x, x\in D)$. 
The part of this process in $U\subset \overline{U}\subset D$ is $Y^{D, \kappa, U}=(Y^{D, \kappa, U}_t, \P_x, x\in U)$. 
\begin{lemma}\label{l:martingale-U}
Suppose that $U$ is a $C^{1,1}$-open set and $U\subset \overline{U}\subset D$.
For any $u\in C_c^2(U)$ and any $x\in U$,
\begin{equation}\label{e:martingale-U}
M_t^u:=u(Y^{D, \kappa, U}_t)-u(Y^{D, \kappa, U}_0)-\int_0^t L^\sB_U u(Y^{D, \kappa, U}_s)\, ds
\end{equation}
is a $\P_x$-martingale  with respect to the filtration of $Y^{D, \kappa, U}$.
\end{lemma}
\pf We follow the proof of \cite[Lemma 2.2]{GKK}. If $(A, \DD(A))$ denotes the $L^2$-generator of the semigroup $T_t$ of $Y^{D, \kappa, U}$, 
then using an argument similar to that used in Proposition \ref{p:operator-interpretation} and Corollary \ref{c:L2-generator} we get that 
$C_c^2(U)\subset \DD(A)$ and $A_{|C_c^2(U)}=(L^\sB_U)_{|C_c^2(U)}$. 
Then in the same way as in the proof of \cite[Lemma 2.2]{GKK} we get that
$$
T_t u(x)-u(x)=\int_0^t T_s L^\sB_U u(x)\, ds\,  \qquad \textrm{a.e. } x\in U\, . 
$$
Moreover, by the display before this lemma we know that $L^\sB_U u$ is bounded in $U$. 
Also, by Lemma \ref{l:YDkappaU}, $Y^{D, \kappa, U}$ is strongly Feller. Now we can follow word-by-word the second part of the proof of 
\cite[Lemma 2.2]{GKK} to get the desired conclusion.
\qed

Recall that, for an open set $U\subset D$, $\tau_U=\tau_U^{Y^{D,\kappa}}=\inf\{t>0:\, Y^{D,\kappa}_t\notin U\}$.

Note that $L^\sB_Uu(Y^{D, \kappa, U}_s)=L^\sB u(Y^{D, \kappa}_s)\mathbf{1}_{s<\tau_U\wedge \zeta}
=L^\sB u(Y^{D, \kappa}_s)\mathbf{1}_{s<\tau_U}$ and that $u(Y^{D, \kappa, U}_t)=u(Y^{D, \kappa}_t)\mathbf{1}_{t<\tau_U}$. 
Thus we can rewrite \eqref{e:martingale-U} as
\begin{equation}\label{e:martingale_U2}
M_t^u=u(Y^{D, \kappa}_t)\mathbf{1}_{t<\tau_U}-u(Y^{D, \kappa}_0)-\int_0^{t\wedge \tau_U}L^\sB u(Y^{D,\kappa}_s)ds\, .
\end{equation}

\begin{prop}\label{p:Dynkin}
For $u\in C_c^2(D)$ and $x\in D$, 
$$
M^u_t:=
u(Y^{D, \kappa}_t)-u(Y^{D, \kappa}_0)-\int_0^{t}L^\sB u(Y^{D,\kappa}_s)ds
$$
is a $\P_x$-martingale.
\end{prop}

\pf 
Let $(U_n)$ be an increasing sequence of relatively compact 
$C^{1,1}$ open subsets of $D$ such that $x\in U_1$, $\mathrm{supp}(u)\subset U_1$ and $U_n\uparrow D$. 
Let $M^{(n)}$ denote the martingale in \eqref{e:martingale_U2} with $U=U_n$, that is,
$$
M^{(n)}_t=u(Y^{D, \kappa}_t)\mathbf{1}_{t<\tau_{U_n}}-u(Y^{D, \kappa}_0)-\int_0^{t\wedge \tau_{U_n}}L^\sB u(Y^{D,\kappa}_s)ds\, .
$$
By using that both $u$ and 
$L^\sB u$ are bounded 
(by Proposition \ref{p:operator-interpretation} (c)), 
we conclude that
$M^u_t$
is a $\P_x$-martingale.
\qed

\medskip
Let $u\in C_c^2(D)$. Consider $y\in D$ and $r >0$ such that $B(y,2r)\subset D$. Let $\tau=\tau_{B(y,r)}$. 
Since $u$ and $L^\sB u$ are bounded, applying the optional stopping theorem to $M^u_t$ we get
$$
\E_y [u(Y^{D, \kappa}_{t\wedge \tau})]-u(y)=\E_y \left [\int_0^{t\wedge \tau}L^\sB u(Y^{D,\kappa}_s)ds\right]\, .
$$
Note that for $s <\tau$ it holds that $\delta_D(Y^{D, \kappa}_s)>r$ $\P_y$-a.s.  Now Proposition \ref{p:operator-interpretation} implies that $|L^\sB u(Y^{D, \kappa}_s)|\le C \Phi(r)^{-1}$ for $s<\tau$. 
Therefore, 
\begin{equation}\label{e:aaaa}
\E_y [u(Y^{D, \kappa}_{t\wedge \tau_{B(y,r)}})]-u(y)\le C \Phi(r)^{-1}t,  \quad \text{for all } r\in (0,1] \text{ and } B(y,2r)\subset D.
\end{equation}

\subsection{Harnack's inequality}\label{HI}
In this subsection we assume that \textbf{(B1)}-\textbf{(B5)} and 
 \eqref{e:assumption-on-j}-\eqref{e:bound-on-kappa} hold true.

For any $x\in D$ and Borel subset $A$ of $D_{\partial}$, we define
$N(x, A)=\int_{A\cap D}J^D(x, y)dy+\kappa(x){\bf 1}_A(\partial)$.  
Then it is known that $(N, t)$ is a L\'evy system for $Y^{D,\kappa}$ (cf.~\cite[Theorem 5.3.1]{FOT} and the argument in \cite[p.40]{CK03}),
that is, for any
non-negative Borel function $f$ on $D\times D_{\partial}$ vanishing on the
diagonal and any stopping time $T$,
$$
\E_x\sum_{s\le T}f(Y^{D,\kappa}_{s-}, Y^{D,\kappa}_s)=
\E_x\left(\int^T_0\int_{D_{\partial}}f(Y^{D,\kappa}_s, y)N(Y^{D,\kappa}_s, dy)ds\right), \quad x\in D.
$$

The following lemma is a simple consequence of  Assumptions \textbf{(B3)} and \textbf{(B5)}.
\begin{lemma}\label{l:two-estimates-J}
\begin{itemize}
\item[(a)] 
For all  $\varepsilon\in (0, 1)$, $x_0\in D$ and 
$r>0$
 such that $B(x_0, (1+\varepsilon)r)\subset D$,  we have
\begin{equation}\label{e:two-estiamtes-J-a}
C_3(\varepsilon/2)
 j(|x-y|)
\le J^{D}(x, y), \qquad x, y\in B(x_0, r).
\end{equation}
\item[(b)] 
For any $\varepsilon\in (0, 1]$, there exists a constant $=C_{11}(\varepsilon)\ge 1$  
such that for all $x_0\in D$ and all 
$r>0$ satisfying $B(x_0, (1+\varepsilon)r)\subset D$, it holds that
\begin{equation}\label{e:estimate-of-J-away}
J^{D}(z,x_1)\le C_{11} J^{D}(z,x_2)\, ,\quad x_1,x_2\in B(x_0,r), \ \ z\in D\setminus B(x_0, (1+\varepsilon)r)\, .
\end{equation}
\end{itemize}
\end{lemma}
\pf (a) For $x,y\in B(x_0,r)$, we have $\delta_D(x)\wedge \delta_D(y)\ge \varepsilon r$ and $|x-y|\le 2r$. Hence, $\delta_D(x)\wedge \delta_D(y)\ge (\varepsilon_0/2)|x-y|$. Thus by \textbf{(B3)}, 
$
J^D(x,y)=\sB(x,y)j(|x-y|)\ge C_3 j(|x-y|),
$
where $C_3=C_3(\varepsilon/2)$.

\noindent 
(b) By using \textbf{(B5)} we have $ J^{D}(z,x_1)=\sB(z,x_1)j(z,x_1)\le 
C_5 \sB(z,x_2) j(z,x_1)$. Since  $|x_1-z|\le |x_2-z|+|x_1-x_2|\le |x_2-z|+2r\le |x_2-z| +(2/\varepsilon)|x_2-z|=(1+2\varepsilon) |x_2-z|$, it follows from  \eqref{e:H-infty} that $j(|z-x_1|)\le c_1 j(|z-x_2|)$. This proves the estimate. \qed

\begin{lemma}\label{l: exit-time-probability}
There exists a constant $C_{12}>0$ such that for all 
$x\in D$ and $r>0$ with  $B(x,2r)\subset D$,
 it holds that
$$
\P_x(\tau_{B(x,r)} < t \wedge \zeta)\le C_{12}\Phi(r)^{-1} t\, .
$$
\end{lemma}
\pf Let $x\in D$ and $r >0$ be such that $B(x,2r)\subset D$. Let $f:\R^d \to [-1, 0]$ be a $C^2$ function such that $f(z)=-1$ for $|z|\le 1/2$, $f(y)=0$ for $|z|\ge 1$ and that $\|\nabla f\|_{\infty}+\|\partial^2 f\|_{\infty}=:c_1<\infty$. Define
$$
f_r(y):=f\left(\frac{y-x}{r}\right)\, .
$$
Then $f_r\in C_c^2(D)$, $f_r(y)=-1$ for $y\in B(x,r/2)$ and $f_r(y)=0$ for $y\in D\setminus B(x,r)$. Moreover, $\|\nabla f_r\|_{\infty}\le c_1/r$, $\|\partial^2 f_r\|_{\infty}\le c_1/r^2$.  
Using these and  applying Proposition \ref{p:operator-interpretation} (c), we get 
$
\|L^\sB f_r\|_{\infty}\le c_2 \Phi(r)^{-1}.
$

Since $f_r\in C_c^2(D)$, by Proposition \ref{p:Dynkin} that
$$
f_r(Y^{D, \kappa}_t)-f_r(Y^{D, \kappa}_0)-\int_0^t L^\sB f_r (Y^{D, \kappa}_s)\, ds
$$
is a $\P_y$-martingale for every $y\in D$. Thus for any $y\in B(x,r/2)$ and any stopping time $S$ with $\E_y[S]<\infty$,
\begin{align}\label{e:lemma3.7}
&\P_y(|Y^{D, \kappa}_S-x|\ge r, S<\zeta)=\E_y[1+f_r(Y^{D, \kappa}_S), |Y^{D, \kappa}_S-x|\ge r, S<\zeta) \nonumber \\
&\le \E_y[1+f_r(Y^{D, \kappa}_S)]=-f_r(y)+\E_y[f_r(Y^{D, \kappa}_S)]=\E_y\left[\int_0^S L^\sB f_r(Y^{D, \kappa}_s)\, ds\right] \nonumber\\
&\le \|L^\sB f_r\|_{\infty} \E_y[S] \le c_2 \Phi(r)^{-1} \E_y[S]\, .
\end{align}
The first inequality in the second line follows because $1+f_r\ge 0$. Note that here $f_r(Y^{D, \kappa}_S)$ makes sense regardless whether $S<\zeta$ or not (by definition $f_r(\partial)=0$). 

Now take $S=\tau_{B(x,r)}\wedge t$ and notice that $\{|Y^{D, \kappa}_{\tau_{B(x,r)}\wedge t}-x|\ge r\}=\{\tau_{B(x,r)}<t\}$. Indeed, $|Y^{D, \kappa}_{\tau_{B(x,r)}\wedge t}-x|\ge r$ implies that $\tau_{B(x,r)}\le t$. Consequently we have that
$$
\big\{|Y^{D, \kappa}_{\tau_{B(x,r)}\wedge t}-x|\ge r, \tau_{B(x,r)}\wedge t <\zeta\big\}=\big\{\tau_{B(x,r)}<t\wedge \zeta\big\}\, .
$$
Now it follows form \eqref{e:lemma3.7} that for $y\in B(x,r)$ we have
$$
\P_y\big(\tau_{B(x,r)}<t\wedge \zeta\big)=\P_y\big(|Y^{D, \kappa}_{\tau_{B(x,r)}\wedge t}-x|\ge r, \tau_{B(x,r)}\wedge t <\zeta\big)\le c_2 \Phi(r)^{-1}\E_y [\tau_{B(x,r)}\wedge t] \le c_2 \Phi(r)^{-1} t\, .
$$
\qed

\begin{lemma}\label{l:new}
For all $x\in D$ and all $r>0$ such that $B(x,2r)\subset D$, it holds that
$
\P_x(\tau_{B(x,r)}=\zeta < t)\le C_{13}\Phi(r)^{-1}t,
$
where $C_{13}=C_{1}.$
\end{lemma}

\pf
By the L\'evy system formula,
\begin{align*}
\P_x(\tau_{B(x,r)}=\zeta < t)=\E_x\sum_{s<t}{\bf 1}_{B(x, r)\times\{\partial\}}(Y^{D, \kappa}_{s-}, Y^{D, \kappa}_s)
=\E_x\int^t_0{\bf 1}_{B(x, r)}(Y^{D, \kappa}_s)\kappa(Y^{D, \kappa}_s)ds.
\end{align*}
Since $\kappa(y)\le C_1/\Phi(r)$ for $y\in B(x, r)$ by \eqref{e:bound-on-kappa}, we immediately get
$
\P_x(\tau_{B(x,r)}=\zeta < t)\le  C_{1}\Phi(r)^{-1}t.
$
\qed

Recall that $A(x, r_1, r_2)$ denotes the
annulus $\{y\in \R^d: r_1\le |y-x|<r_2\}$.
\begin{prop}\label{p:exit-time-estimate}
(a) 
 There exists a constant 
$C_{14}>0$ depending on $  a_1, a_2, \delta_1, \delta_2$
such that for all$x_0\in D$ and   $r>0$
 with $B(x_0,r) \subset D$, it holds that
$$
\E_x \tau_{B(x_0,r)}\ge C_{14} \Phi(r)\, ,
\quad x\in B(x_0,r/3).
$$

\noindent
(b)
 For every $\varepsilon >0$,  there exists 
$C_{15}>0$  depending on $\varepsilon$, $a_1$, $a_2$, $\delta_1$, $\delta_2$, 
$C_3$, $C_5$
such that for all $x_0\in D$ and all 
$r >0$
 satisfying $B(x_0, (1+\eps)r)\subset D$, it holds that
$$
\E_x \tau_{B(x_0,r)}\le C_{15}\Phi(r)\, ,\quad x\in B(x_0,r)\, .
$$
\end{prop}
\pf 
(a) 
Let $x\in D$ and $r>0$ 
be such that $B(x,r)\subset D$. 
It follows from Lemmas \ref{l: exit-time-probability}--\ref{l:new} and \eqref{e:H-infty}
that
$$
\P_x(\tau_{B(x, r/3)}< t) \le c_{1}\Phi(r)^{-1}t\, .
$$
Therefore,
$$
\E_x \tau_{B(x, r/3)}\ge t \P_x( 
\tau_{B(x, r/3)}
\ge t)\ge t(1-c_{1}\Phi(r)^{-1}t)\
$$
for all $t>0$. Choose $t=\Phi(r)/(2c_{1})$, so that $1-c_{1}\Phi(r)^{-1}t=1/2$. Then
$$
\E_x \tau_{B(x, r/3)} \ge \frac12 \Phi(r)/(2c_{1})=c_{2}\Phi(r)\, .
$$
Now let $B(x_0,r)\subset D$ and $x\in B(x_0,r/3)$. Then 
$B(x,r/2)\subset B(x_0,r)\subset D$. 
By what was proven above,
$$
\E_x \tau_{B(x_0,r)}\ge \E_x \tau_{B(x,r/6)}\ge c_{2}\Phi(r/2)\ge c_3 \Phi(r)\, .
$$

\noindent
(b)
Let $\varepsilon_0 :=\varepsilon/3$, $x_0\in D$ and $r>0$
so that $B(x_0, (1+3\eps_0)r)\subset D$.
 For $x\in B(x_0,r)$, by using Lemma \ref{l:two-estimates-J} (b)  in the third line,
\begin{eqnarray}\label{e:exit-time-estimate-1}
1&\ge &\P_x\big(Y^{D, \kappa}_{\tau_{B(x_0,r)}} \in A(x_0, (1+\varepsilon_0)r, (1+2\varepsilon_0)r)\big)\nonumber\\
&= &\E_x\int_0^{\tau_{B(x_0,r)}}\int_{A(x_0, (1+\varepsilon_0)r, (1+2\varepsilon_0)r)}J^D(u,Y^{D,\kappa}_s)\, du\, ds\nonumber \\
&\ge &C_{11}^{-1}c_1\E_x\tau_{B(x_0,r)}\int_{A(x_0, (1+\varepsilon_0)r, (1+2\varepsilon_0)r)}J^D(u,x_0)\, du.
\end{eqnarray}

For $u\in A(x_0,r, (1+2\varepsilon_0)r)$, we have by Lemma \ref{l:two-estimates-J} (a) that $J^D(u,x_0)\ge C_{11} 
j(|u-x_0|)$.
Therefore, 
by \eqref{e:assumption-on-j} and \eqref{e:H-infty},
\begin{eqnarray*}
\int_{A(x_0, (1+\varepsilon_0)r, (1+2\varepsilon_0)r)}J^D(u,x_0)\, du\ &\ge& C_{11}\int_{A(x_0, (1+\varepsilon_0)r, (1+2\varepsilon_0)r)}j(|u-x_0|)\, du\\
&\ge & c_3 \int_{(1+\varepsilon_0)r}^{(1+2\varepsilon_0)r}\frac{1}{t\Phi(t)}\, dt \ge c_4 \frac{1}{\Phi(r)}\, .
\end{eqnarray*}
By inserting this in \eqref{e:exit-time-estimate-1} we obtain
$
1\ge c_5{\E_x\tau_{B(x_0,r)}}/{\Phi(r)},
$
which is the required inequality. \qed

\begin{lemma}\label{l:A2}
There exists $C_{16}>0$ such that for 
all $x\in D$ and $r>0$ with $B(x,5r)\subset D$,
 and any Borel $A\subset B(x,r)$,
$$
\P_y(T_A<\tau_{B(x,3r)})\ge C_{16}\frac{|A|}{|B(x,r)|}\, , \qquad y\in B(x,2r).
$$
\end{lemma}
\pf Without loss of generality assume that ${\mathbb P}_y(T_A < {\tau}_{B(x,3r)})<1/4$.
Set $\tau={\tau}_{B(x,3r)}$. For $y\in B(x,2r)$ we have that $B(y,3r)\subset D$. Hence by Lemmas \ref{l: exit-time-probability} and \ref{l:new}, for every 
$y\in B(x,2r)$,
$
\P_y(\tau< t) \leq \P_y(\tau_{B(y,r)}< t) \leq c_1\Phi(r)^{-1} t\, .
$
 Choose $t_0= \Phi(r)/(4c_1)$, so that $\P_y(\tau< t_0) \leq 1/4$. Further, if $z\in B(x, 3r)$ and $u\in A \subset B(x,r)$, then $|u-z| \leq 4r$.  Since $j$ is decreasing,
$j(|u-z|) \geq j(4r)$. Moreover, $\delta_D(u)\wedge \delta_D(z)\ge r \ge (1/4)|u-z|$, implying by {\bf (B3)} that $\sB(u,z)\ge C_3(1/4)=C_3$. 
Thus,
\begin{align*}
&\P_y (T_A < \tau) \geq  \E_y \sum_{s\leq T_A \wedge \tau \wedge t_0}
{\mathbf 1}_{\{Y^{D, \kappa}_{s-}\neq Y_s, Y^{D, \kappa}_s\in A\}} 
\\
 &=  \E_y \int_0^{T_A \wedge \tau \wedge t_0}
\int_A j(|u-Y^{D, \kappa}_s|)\sB(u, Y^{D, \kappa}_s)\, du \, ds 
 \ge   C_3j(4r) |A| \E_y[T_A \wedge \tau \wedge t_0] \, ,
\end{align*}
where in the second line we used properties of the L\'evy system.
Next,
\begin{align*}
{\mathbb E}_y[T_A\wedge \tau \wedge t_0]   
 \ge  t_0 {\mathbb P}_y(T_A \ge \tau \ge t_0) 
 \ge  t_0[1-{\mathbb P}_y(T_A < \tau)-{\mathbb P}_y(\tau <t_0)] 
 \ge  \frac{t_0}{2} = \frac{\Phi(r)}{8 c_1}\, .
\end{align*}
The last two displays give that
$$
{\mathbb P}_y (T_A < \tau) \geq C_3 j(4r) |A| \frac{\Phi(r)}{8 c_1} \ge c_2 \frac{1}{r^d \Phi(r)}|A|\Phi(r)\ge c_3\frac{|A|}{|B(x,r)|}\, .
$$
 \qed

\begin{lemma}\label{l:A3}
 There exist $C_{17}>0$ and $C_{18}>0$
 with the property that if 
$r>0$, 
$x\in D$ are such that $B(x,2r)\subset D$,  and $H$ is a bounded
non-negative function with support in $D\setminus B(x,2r)$, then for every $z\in B(x,r)$,
$$
 C_{17} \E_z [{\tau}_{B(x,r)}] \int H(y)
J^D(x,y) \, dy \le \E_z H(Y^{D, \kappa}_{{\tau}_{B(x,r)}}) \le  C_{18} {\mathbb E}_z [{\tau}_{B(x,r)}] \int H(y)
J^D(x,y) \, dy \, .
$$
\end{lemma}

 \pf Let $y\in B(x,r)$ and $u\in D\setminus B(x,2r)$. By \textbf{(B5)} (or Lemma \ref{l:two-estimates-J} (b)), $J^D(u,y)\asymp J^D(x,y)$. Thus using the L\'evy system  we get
\begin{align*}
\E_z \left[ H(Y^{D, \kappa}_{\tau_{B(x,r)}}) \right]  &= \E_z \int_0^{\tau_{B(x,r)}} \int H(u) J^D(u, Y^{D, \kappa}_s)\, du\, ds \asymp \E_z \int_0^{\tau_{B(x,r)}} \int H(u)  J^D(u, x)\, du\, ds.
\end{align*}
 \qed

\noindent
\textbf{Proof of Theorem \ref{t:uhp}.} (a) Proposition \ref{p:exit-time-estimate} (a)-(b),  
Lemmas \ref{l:A2} and \ref{l:A3} imply that
conditions (A1), (A2) and (A3) of \cite{SV04} are satisfied for the  process $Y^{D, \kappa}$.
 Thus we can repeat the proofs of 
\cite[Theorems 2.2 and 2.4]{SV04} to finish the proof of part (a).
Note that conservativeness does not play any role. We omit the details.

\noindent
(b) 
By (a) we can and will assume that $L>2$ and $2r < |x_1-x_2|<Lr$. 
For simplicity, let $B_i=B(x_i,r)$, $i=1,2$. 
Then by using harmonicity in the first inequality, Theorem \ref{t:uhp}.(a) in the second inequality, and the L\'evy system formula in the second line, we have
\begin{eqnarray}\label{e:hp2-1}
f(x_1)&\ge &\E_{x_1}\big[f(Y^{D,\kappa}_{\tau_{B_1}}); Y^{D, \kappa}_{\tau_{B_1}}\in B(x_2, r/2)\big] \ge C_{11}^{-1} f(x_2) \P_{x_1}\big(Y^{D, \kappa}_{\tau_{B_1}}\in B(x_2, r/2)\big)\nonumber \\
&=& C_{11}^{-1} f(x_2)\E_{x_1} \int_0^{\tau_{B_1}}  \int_{B(x_2, r/2)}J^D(Y^{D, \kappa}_s,z)\, dz\, ds\, .
\end{eqnarray}
For $y\in B_1$ and $z\in B(x_2, r/2)$ we have by Lemma \ref{l:two-estimates-J}.(b) (with $\varepsilon=1$) that $J^D(y,z)\ge C_{11}^{-1}J^D(x_1, z)$. Further, $\delta_D(x_1)\wedge \delta_D(z)\ge r/2 \ge 
(2L+2)^{-1}
|x_1-z|$, hence by Lemma \ref{l:two-estimates-J} (a), $J^D(x_1,z)\ge C_{3} 
j(|x_1-z|)$, where
$C_{3}=C_3((2L+2)^{-1})$. 
By inserting this in \eqref{e:hp2-1}, 
and by using Proposition \ref{p:exit-time-estimate} (a) , we obtain
\begin{align*}
f(x_1)&\ge C_3 f(x_2) \E_x\tau_{B_1}\int_{B(x_2,r/2)}j(|x_1-z|)\, dz\\ 
&\ge c_1 C_3 f(x_2) \Phi(r) \frac{1}{((L+1)r)^d \Phi((L+1)r)}\, |B(x_2, r/2)|\\
& \ge C_3 f(x_2) L^{-d} \frac{\Phi(r)}{\Phi((L+1)r)}\ge  c_4 f(x_2) L^{-d} L^{-2\delta_2}\, .
\end{align*}
The last inequality follows from \eqref{e:H-infty}. \qed

By repeating the arguments of \cite[Theorem 4.9]{SV04}
and \cite[Theorem 4.1]{BL}, we immediately get the following result.

\begin{thm}\label{t:holder}
Let $f$ be a non-negative bounded function in $D$ which is harmonic in $B(x_0, r)$ with respect
to $Y^{D, \kappa}$.
Then there exist  $C_{19}>0$ and $\beta>0$ such that for any $r\in (0, 1]$ and $B(x_0, r)\subset D$,
$$
|f(x)-f(y)|\le C_{19}\|f\|_\infty 
(|x-y|/r)^\beta,  
\qquad x, y\in B(x_0, r/2).
$$
\end{thm}


\section{Carleson's estimate}\label{s:carleson}
In this section we prove Carleson's estimate. In addition to \textbf{(B1)-(B5)} and \eqref{e:bound-on-kappa}, we also assume that 
\textbf{(B6)}  and \eqref{e:l-bound-on-kappa} hold. 
\begin{lemma}\label{lower bound} 
There exists a constant $\delta_*>0$ such that for 
all $x\in D$,
$
\P_x\left(  \tau_{B(x,\delta_D(x)/2)}=\zeta     \right)\ge \delta_*
$ where   $\zeta$ is the lifetime of $Y^{D, \kappa}$.
\end{lemma}
\pf 
Let $x\in D$. 
By the L\'evy system formula,  \eqref{e:l-bound-on-kappa} and Proposition \ref{p:exit-time-estimate} (a),
\begin{align*}
&\P_x\left( \tau_{B(x,\delta_D(x)/2)}=\zeta     \right) \ge \P_x\left( \tau_{B(x,\delta_D(x)/3)}=\zeta     \right)
=
\E_x \int_0^\infty {\bf 1}_{B(x,\delta_D(x)/3)} (Y^{D, \kappa}_s) \kappa(Y^{D, \kappa}_s)ds\\
&\ge 
\E_x \int_0^{\tau_{B(x,\delta_D(x)/3)}} \kappa(Y^{D, \kappa}_s)ds
\ge \frac{c_3}{\Phi(\delta_D(x)/3)}\E_x\tau_{B(x,\delta_D(x)/3)}\ge c_4.
\end{align*}
\qed

\medskip
\noindent
\textbf{Proof of Theorem \ref{t:carleson}.}
In this proof, the constants $\delta_*, \nu, \gamma, \chi, \eta$ and $c_i$'s are always
independent of $r$. 
Let $f$ be a non-negative function on $D$ which is harmonic in $D\cap B(Q, r)$ and vanishes
continuously on $\partial D\cap B(Q, r)$. 
By Theorem \ref{t:uhp} (b)
and a standard chain argument, it suffices to prove (\ref{e:carleson}) for 
$x\in D\cap B(Q, \overline{\kappa} r/(24))$.

Let $\nu:=d+2\delta_2$. 
Recall that $\wh{\beta}$ is the constant from Assumption \textbf{(B6)}.
Choose $0<\gamma < 2\delta_1/(\nu+\wh{\beta})$. 
Fix an $x_0\in D\cap B(Q, r)$ with $\delta_D(x_0)\ge \overline{\kappa}r/2$.
For any $x\in D\cap B(Q,\overline{\kappa} r/(12))$, define
$$
B_0(x)=B(x,\delta_D(x)/2)\, ,
\qquad B_1(x)=B(x,r^{1-\gamma}\delta_D(x)^{\gamma})\,
$$
and
$$
B_2=B(x_0,\overline{\kappa}\delta_D(x_0)/3)\, ,\qquad 
B_3=B(x_0, 2\overline{\kappa}\delta_D(x_0)/3).
$$
Since $x\in B(Q,\overline{\kappa} r/(12))$, we have $\delta_D(x)<r/(12)$. 
In particular, we have that $B_0(x)\subset B(x,\delta_D(x)) \subset B_1(x)$.
By Lemma \ref{lower bound}, there exists $\delta_*=\delta_*(R,
\Lambda)>0$ such that
\begin{equation}\label{e:c:1}
\P_x(\tau^{Y^{D,\kappa}}_{B_0(x)}=\zeta)\ge \delta_*\, , \quad 
x\in D\cap B(Q,\overline{\kappa} r/(12))\, .
\end{equation}
By Theorem \ref{t:uhp}.(b) there exists $\chi>0$ such that
\begin{equation}\label{e:c:2}
f(x)<(\delta_D(x)/r)^{-\chi} f(x_0)\, ,\quad 
x\in D\cap B(Q,\overline{\kappa} r/(12))\, .
\end{equation}
Since $f$ is harmonic in $D\cap B(Q,r)$,  for every 
$x\in D\cap B(Q,\overline{\kappa} r/(12))$,
\begin{align}
&f(x)=\E_x\big[f\big(Y^{D,\kappa}(\tau_{B_0(x)})\big); Y^{D,\kappa}(\tau_{B_0(x)})\in B_1(x)\big]\nonumber\\
&\qquad+ \E_x\big[f\big(Y^{D,\kappa}(\tau_{B_0(x)})\big);Y^{D,\kappa}(\tau_{B_0(x)})\notin B_1(x)\big]. \label{e:c:3}
\end{align}
We first show that there exists $\eta>0$ such that for all 
$x\in D \cap B(Q, \overline{\kappa} r/(12)) $ with $\delta_D(x) < \eta r$,
\begin{equation}\label{e:c:4}
\E_x\big[f\big(Y^{D,\kappa}(\tau_{B_0(x)})\big); Y^{D,\kappa}(\tau_{B_0(x)})\notin
B_1(x)\big]\le f(x_0). 
\end{equation}

\noindent
Step 1: There exists $c_2>0$ such that
\begin{equation}\label{e:carleson-2}
f(x_0)\ge c_2 \Phi(r) \int_{D\setminus B_3}J^D(x_0,y) f(y)\, dy\, .
\end{equation}
Indeed, if $z\in B_2$ and $y\in D\setminus B_3$, then  by Lemma \ref{l:two-estimates-J} (b), 
 $J^{D}(z,y)\ge c_4 J^{D}(x_0,y)$. By using this estimate in the second inequality below, 
Proposition \ref{p:exit-time-estimate} (a)  in the third
 and $\delta_D(x_0) \ge \kappa r/2$ in the fourth, we get that
\begin{align*}
&f(x_0)\ge \E_{x_0}\left[f(Y^{D,\kappa}(\tau_{B_2})); Y^{D,\kappa}(\tau_{B_2})\notin B_3\right]=    \E_{x_0} \int_0^{\tau_{B_2}}\left( \int_{D\setminus B_3} J^{D}(Y^{D,\kappa}_t, y)f(y)\, dy \right) dt\nonumber \\
&\ge c_5  \E_{x_0}[\tau_{B_2}]\int_{D\setminus B_3}  J^{D}(x_0, y)f(y)\, dy \ge c_6 \Phi(\delta_D(x_0))\int_{D\setminus B_3}  J^{D}(x_0, y)f(y)\, dy \nonumber \\
&\ge c_7 \Phi(r)\int_{D\setminus B_3}  J^D(x_0, y)f(y)\, dy\, .
\end{align*}

\noindent
Step 2: For $z\in B_0(x)$ and 
$y\in D \setminus  B_1(x)$,
 by Lemma \ref{l:two-estimates-J} (b) it holds that $J^D(z,y)\le c_8 J^D(x,y)$. By using this in the second line below and 
Proposition \ref{p:exit-time-estimate} (b)
in the third, we get
\begin{align}\label{e:carleson-4}
&\E_x\big[f\big(Y^{D,\kappa}(\tau_{B_0(x)})\big); Y^{D,\kappa}(\tau_{B_0(x)}) \notin B_1(x)\big]
= \E_x \int_0^{\tau_{B_0(x)}}   \int_{D\setminus B_1(x)}J^{D}(Y^{D,\kappa}_t,y)f(y)\, dy\, dt  \nonumber\\
&\le  c_{9}\E_x\tau_{B_0(x)} \int_{D\setminus B_1(x)} J^{D}(x,y)f(y)\, dy \nonumber\\
&\le c_{10} \Phi( \delta_D(x)) \left(\int_{(D\setminus B_1(x))\cap B_3^c}J^{D}(x,y)f(y)\, dy 
+ \int_{(D\setminus B_1(x))\cap B_3}J^{D}(x,y)f(y)\, dy \right)\nonumber\\
&=:c_{10} \Phi(\delta_D(x))(I_1+I_2).
\end{align}

\noindent
Step 3:
Suppose now that $|y-x|\ge r^{1-\gamma}\delta_D(x)^{\gamma}$ and  $x\in B(Q, \overline{\kappa} r/(12))$. Then
$$
|y-x_0|\le |y-x|+2r\le |y-x|+2r^{\gamma}\delta_D(x)^{-\gamma}|y-x|\le 3r^{\gamma}\delta_D(x)^{-\gamma}|y-x|,
$$
and $\delta_D(x)\le \overline{\kappa} r/12< \overline{\kappa} r/2 \le \delta_D(x_0)$. 
Set $M=3r^{\gamma}/\delta_D(x)^{\gamma}\ge 1$ so that $|y-x_0|\le M|y-x|$.
Thus, by using \textbf{(B6)}  in the 
first inequality and 
\eqref{e:H-infty}
in the second inequality,
\begin{align}\label{e:gf1}  
{J^D(x, y)\le C_6 M^{\wh{\beta}}\sB(x_0,y) j(M^{-1}|y-x_0|)}\le c_{11} \left( \frac{\delta_D(x)}{r}\right)^{-\gamma(\wh{\beta}+\nu)} J^D(x_0, y)\, .
\end{align}
Now, using this and \eqref{e:carleson-2} , we get
\begin{align}\label{e:c:7}
&I_1\le c_{12} \left(\frac{\delta_D(x)}{r}\right)^{- \gamma (\nu+\wh{\beta})}\int_{(D\setminus B_1(x))\cap B_3^c }J^D(x_0, y) f(y) \, dy\nonumber\\
&\le  c_{13}\Phi(r)^{-1} \left(\frac{\delta_D(x)}{r}\right)^{- \gamma (\nu+\wh{\beta})} f(x_0).
\end{align}
 
\noindent
Step 4:
If $y\in B_3$, then $\delta_D(y)\le c_{12}r$ and 
$|y-x|\ge |x_0-Q|-|x-Q|-|y-x_0|>\delta_D(x_0)/2$.
By Theorem \ref{t:uhp}, there exists $C>0$ such that $f(y)\le C f(x_0)$ 
for all $y\in B_3$.
Thus,
\begin{align}\label{e:c:8}
&I_2 \le c_{14} f(x_0) \int_{(D\setminus B_1(x))\cap B_3} J^D(x, y) \, dy\le  c_{14} f(x_0)\int_{|y-x|>
\delta_D(x_0)/2}  J^D(x, y)\, dy\nonumber\\
 &\le   c_{15} f(x_0)\int_{|y-x|>
 \delta_D(x_0)/2
 } j(|x-y|)\, dy\le   c_{16} \Phi(r)^{-1}f(x_0)\, .
\end{align}

\noindent
Step 5:
Combining \eqref{e:carleson-4}, \eqref{e:c:7}, \eqref{e:c:8}, and using  \eqref{e:H-infty} and the definition of $\Phi$ in the last line, we obtain
\begin{align}\label{e:c:9}
&\E_x[f(Y^{D,\kappa}(\tau_{B_0(x)}));\, Y^{D,\kappa}(\tau_{B_0(x)})\notin B_1(x)]\nonumber\\
&\le c_{17} f(x_0)\left( \frac{\Phi(\delta_D(x))}{\Phi(r)} \left(\frac{\delta_D(x)}{r}\right)^{- \gamma (\nu+\wh{\beta})}+ \frac{\Phi(\delta_D(x))}{\Phi(r)} \right) \nonumber\\
&\le c_{18}f(x_0) \left(\left(\frac{\delta_D(x)}{r}\right)^{- \gamma (\nu+\wh{\beta})+2\delta_1}+ \left(\frac{\delta_D(x)}{r}\right)^{2\delta_1}\right).
\end{align}
  Since $2\delta_1-\gamma(\nu+\wh{\beta})>0$ (by the choice of $\gamma$),  we can choose 
  $\eta>0$ so that
$$
c_{18}\left(\eta^{2\delta_1-\gamma(\nu+\wh{\beta})} +\eta^{2\delta_1}\right)\,\le\, 1\, .
$$
Then  for $x\in  D \cap B(Q, \overline{\kappa} r/(12))$ with $\delta_D(x) < \eta r$, we
have by \eqref{e:c:9},
\begin{align*}
\E_x\left[f(Y^{D,\kappa}(\tau_{B_0(x)}));\, Y^{D,\kappa}(\tau_{B_0(x)})\notin
B_1(x)\right] \le  c_{19}\, f(x_0)
\left(\eta^{2\delta_1-\gamma(\nu+\wh{\beta})} +\eta^{2\delta_1}\right) \le
f(x_0)\, .
\end{align*}
This completes the proof of \eqref{e:c:4}.
We now prove Carleson's estimate \eqref{e:carleson} for 
$x\in D\cap B(Q, \overline{\kappa} r/(24))$ 
by a method of contradiction.
Without loss of generality, we may assume that $f(x_0)=1$.
Suppose that there exists 
$x_1\in D\cap B(Q, \overline{\kappa} r/(24))$ such
that $f(x_1)\ge K>\eta^{-\chi}\vee (1+\delta_*^{-1})$, where $K$ is a
constant to be specified later. By \eqref{e:c:2} and the assumption
$f(x_1)\ge K>\eta^{-\chi}$, we have
$(\delta_D(x_1)/r)^{-\chi}>f(x_1)\ge K>
 \eta^{-\chi}$, and hence
$\delta_D(x_1)<\eta r$.
By (\ref{e:c:3}) and  (\ref{e:c:4}),
$$
K\le f(x_1)\le \E_{x_1}\left[f(Y^{D}(\tau_{B_0(x_1)}));
Y^{D,\kappa}(\tau_{B_0(x_1)}) \in B_1(x_1) \right]+1\, ,
$$
and hence
$$
\E_{x_1}\left[f(Y^{D, \kappa}({\tau_{B_0(x_1)}})); Y^{D, \kappa}({\tau_{B_0(x_1)}})\in
B_1(x_1)\right] \ge f(x_1)-1 > \frac{1}{1+\delta_*}\, f(x_1)\, .
$$
In the last inequality of the display above we used the assumption
that  $f(x_1)\ge K>1+\delta_*^{-1}$. 
If $K \ge (24/\overline{\kappa})^{\chi /\gamma}$, then
$(\delta_D(x_1)/r)^\gamma <\overline{\kappa}/(24)$.
Thus $\overline{ B_1(x_1)}\subset B(Q, \overline{\kappa} r/(12))$. 
We now get from \eqref{e:c:1} that
\begin{align*}
&\E_{x_1}[f(Y^{D,\kappa}(\tau_{B_0(x_1)})), Y^{D,\kappa}(\tau_{B_0(x_1)})\in B_1(x_1)]\\
&=\E_{x_1}[
f(Y^{D,\kappa}(\tau_{B_0(x_1)})), Y^{D,\kappa}(\tau_{B_0(x_1)})\in B_1(x_1)\cap D]\\
& \le  \P_{x_1}
\left(Y^{D,\kappa}(\tau_{B_0(x_1)})\in D\right) \, \Big(\sup_{B_1(x_1)}f\Big)  \le (1-\delta_*) \, \Big(\sup_{B_1(x_1)}f\Big) \, .
\end{align*}
Therefore, 
$\sup_{B_1(x_1)}f> f(x_1)/(1-\delta_*^2)$, i.e.,
there exists $x_2\in D\cap B(Q, \overline{\kappa} r/(12))$ such that
$$
|x_1-x_2|\le r^{1-\gamma}\delta_D(x_1)^{\gamma} \quad \hbox{ and }
\quad
f(x_2)>\frac{1}{1-\delta_*^2}\, f(x_1)\ge \frac{1}{1-\delta_*^2}\, K\, .
$$
Similarly, if $x_k\in D\cap B(Q,  \overline{\kappa} r/(12))$
with $f(x_k)\geq
K/(1-\delta_*^2)^{k-1}$ for $k\ge 2$, then there exists 
$x_{k+1}\in D$
such that
\begin{equation}\label{e:c:10}
|x_k-x_{k+1}|\le r^{1-\gamma}\delta_D(x_k)^{\gamma}  \quad \hbox{ and
} \quad f(x_{k+1}) > \frac{1}{1-\delta_*^2}\, f(x_k)>
\frac{1}{(1-\delta_*^2)^k}\, K\, .
\end{equation}
From (\ref{e:c:2}) and (\ref{e:c:10}) it follows that
$(\delta_D(x_{k})/r)^{\gamma} \le (1-\delta_*^2)^{\gamma(k-1)/\chi}K^{-\gamma/\chi}$, for
every $k\ge 1$. 
Therefore by this and \eqref{e:c:10}, 
\begin{align*}
&|x_k-Q|\,\le\,|x_1-Q|+\sum_{j=1}^{k-1}|x_{j+1}-x_j|
\,\le\, \frac{\overline{\kappa} r}{24} +
r \sum_{j=1}^{\infty} 
(\delta_D(x_{k})/r)^{\gamma} 
\\
&\le
\frac{\overline{\kappa} r}{24}
+ r K^{-\gamma/\chi}\,
\frac{1}{1-(1-\delta_*^2)^{\gamma/\chi}}.
\end{align*}
Choose
$$
K=\eta^{-\chi}
\vee (1+\delta_*^{-1})\vee 
[(24/\overline{\kappa})^{\chi/\gamma}(1-(1-\delta_*^2)^{\gamma/\chi})^{-\chi/\gamma}]
$$
so that $K^{-\gamma/\chi}\, (1-(1-\delta_*^2)^{\gamma/\chi})^{-1}\le
\overline{\kappa}/(24)$. Hence 
$x_k\in D\cap B(Q, \overline{\kappa} r/(12))$
for every $k\ge 1$. 
Since $\lim_{k\to \infty}f(x_k)=\infty$ by \eqref{e:c:10},
 this contradicts the fact that
$f$ is bounded on $B(Q,r/2)$. This contradiction shows that $f(x)< K$
for every $x\in D\cap B(Q, \overline{\kappa} r/(24))$.
This completes the proof of the theorem. 
 \qed

\section{Estimates in the half-space}\label{s:estimates}
Starting from this section we assume that $D=\R^d_+$ and 
$j(|x-y|)=|x-y|^{-d-\alpha}$,  $d\ge 1$.  
 In the rest of this section, we will only deal with the case $d\ge 2$, the case $d=1$ being simpler.

For $x\in \R^d_+$, we write $x=(\wt{x}, x_d)$, $\wt{x}\in \R^{d-1}$, $x_d>0$, so that $\delta_D(x)=x_d$.
We have that  $J^{\R^d_+}(x,y)=\sB(x,y)|x-y|^{-d-\alpha}$ and note that with this assumption, 
$\Phi(\lambda)=\lambda^{\alpha}$ 
so that $\delta_1=\delta_2=\alpha/2$. 
Further, for $p\in ((\alpha-1)_+, \alpha+\beta_1)$, let $\kappa(x):=C(\alpha, p, \sB)x_d^{-\alpha}$, where 
$ C(\alpha,p,\sB)$ is the constant defined in \eqref{e:explicit-C}. 
For simplicity, we denote the process $Y^{\R^d_+, \kappa}$ by $Y$. Thus, $Y$ is the process associated with the Dirichlet form
\begin{align*}
\EE(u,v)=\frac12 \int_{\R^d_+}\int_{\R^d_+} \frac{(u(x)-u(y))(v(x)-v(y))}{|x-y|^{d+\alpha}}\sB(x,y) dydx
+C(\alpha, p, \sB)\int_{\R^d_+}u(x)v(x)x_d^{-\alpha}dx.
\end{align*}

We first record the scaling property of $Y$ here. For any $r>0$, define a process $Y^{(r)}$ by
$Y^{(r)}_t:=r Y_{r^{-\alpha} t}$.

\begin{lemma}\label{l:scaling-of-Y}
Assume that \textbf{(B1)},  \textbf{(B2)} and \textbf{(B8)} hold. 
Then $(Y^{(r)}, \P_{x/r})$ has the same law as $(Y, \P_x)$.
\end{lemma}
\pf
For any function $f$, define $f^{(r)}(x)= f(rx)$.  Let $P^{(r)}_t$ be the semigroup of $Y^{(r)}$. 
Then
$$
P^{(r)}_tf(x)=\E^{(r)}_xf(Y^{(r)}_t)=\E_{x/r}f(rY_{r^{-\alpha} t})=P_{r^{-\alpha} t}f^{(r)}(x/r).
$$
Thus for any nice functions $f$ and $g$,
\begin{align*}
\frac1t\int_{\R^d_+}(P^{(r)}_tf(x)-f(x))g(x)dx
=r^{d-\alpha}\frac1{r^{-\alpha} t}\int_{\R^d_+}(P_{r^{-\alpha} t}f^{(r)}(x )-f^{(r)}(x))g^{(r)}(x )dx.
\end{align*}
Since by \textbf{(B8)} $\sB(x/r, y/r)=\sB(x,y)$, 
by letting $t\downarrow 0$ we get
\begin{align*}
&\EE^{(r)}(f, g)=r^{d-\alpha}\EE(f^{(r)}, g^{(r)})\\
&=r^{d-\alpha}\left(\frac12 \int_{\R^d_+}\int_{\R^d_+} \frac{(f^{(r)}(x)-f^{(r)}(y))(g^{(r)}(x)-g^{(r)}(y))}{|x-y|^{d+\alpha}}\sB(x,y) dydx\right.\\
&\qquad\ \left. +\ C(\alpha, p, \sB)\int_{\R^d_+}f^{(r)}(x)g^{(r)}(x)x_d^{-\alpha}dx\right)\\
&=r^{d-\alpha}\left(\frac12 \int_{\R^d_+}\int_{\R^d_+} \frac{(f(rx)-f(ry))(g(rx)-g(ry))}{|x-y|^{d+\alpha}}\sB(x,y) dydx\right.\\
&\qquad\ \left.+\ C(\alpha, p, \sB)\int_{\R^d_+}f(rx)g(rx)x_d^{-\alpha}dx\right)\\
&=\frac12 \int_{\R^d_+}\int_{\R^d_+} \frac{(f(x)-f(y))(g(x)-g(y))}{|x-y|^{d+\alpha}}\sB(x,y) dydx+C(\alpha, p, \sB)\int_{\R^d_+}f(x)g(x)x_d^{-\alpha}dx=\EE(f, g).
\end{align*}
Thus $(Y^{(r)}, \P_{x/r})$ has the same law as $(Y, \P_x)$. \qed

An easy consequence of Lemma \ref{l:scaling-of-Y} is the following: Let $V$ be an open subset of $\R^d_+$, and for $r>0$ denote $rV:=\{rx:\, x\in V\}$. If $\tau_V=\inf\{t>0:\, Y_t\notin V\}$, then
\begin{equation}\label{e:exit-time-scaling}
\E_{rx}\tau_{rV}=r^{\alpha}\E_x \tau_V\, , \qquad x\in V\, .
\end{equation}

\emph{Throughout the remainder of this section and below we always assume that $\sB$ satisfies }
\textbf{(B1)},  \textbf{(B2)} \emph{and} \textbf{(B7)}-\textbf{(B8)}. 
Note that  \textbf{(B8)} implies that $x\mapsto \sB(x, x)$ is a constant on $D$. Without loss
of generality, will assume that $\sB(x, x)=1$.
Additional assumptions on $\sB$ will be imposed later when needed. 

The following observation will be used several times in the sequel. Recall the function $\wt{B}=\wt{B}_{\beta_1, \beta_2, \beta_3, \beta_4}$ appearing in \textbf{(B7)}. 
We note that if $\beta_4>0$, 
then, for any $\eps\in (0, \beta_2)$,  there exists $c_\eps>0$ such that  
\begin{align}\label{e:B7_2} 
&(\log 2)^{-\beta_4}\wt{B}_{\beta_1, \beta_2, \beta_3, 0}(x,y) \le \wt{B}_{\beta_1, \beta_2, \beta_3, \beta_4}(x,y) \le c_{\eps} \wt{B}_{\beta_1, \beta_2-\eps, \beta_3, 0}(x,y).
\end{align}

We first establish some auxiliary estimates that will be useful throughout the rest of the paper.  Recall $\log(1+a)^{\gamma}=(\log(1+a))^{\gamma}$. 

\begin{lemma}\label{l:key-estimate-for-B}
\begin{itemize}
\item[(a)] 
 There exists a constant $C_{20}>0$ such that for all $x,y\in \R^d_+$ satisfying $|x-y|\ge x_d$, it holds  that
\begin{equation}\label{e:key-estimate-for-B}
\sB(x,y)\le C_{20}x_d^{\beta_1}(|\log x_d|^{\beta_3}\vee 1)\big(1+{\bf 1}_{|y|\ge1}(\log|y|)^{\beta_3}\big)|x-y|^{-\beta_1}.
\end{equation}
\item[(b)]  For all $x,y\in \R^d_+$ satisfying $|x-y|\ge x_d$, $x_d\le e^{-1}$ and $|y|\le M$, it holds that
$$
\sB(x,y)\le C_{20}(1+(\log M)^{\beta_3})x_d^{\beta_1}|\log x_d|^{\beta_3}|x-y|^{-\beta_1}.
$$
\end{itemize}
\end{lemma}
\pf 
 By using \eqref{e:B7_2}, 
without loss of generality,  throughout the proof we assume that $\beta_4=0$.

(a) Since $|x-y|\ge x_d$ we get that
$$
\frac{(x_d\vee y_d)\wedge |x-y|}{x_d\wedge y_d \wedge |x-y|}\le  \frac{x_d\vee y_d}{x_d\wedge y_d}\, .
$$
Thus it follows from our assumptions on $\sB$  that
\begin{equation}\label{e:former-remark-3}
\sB(x,y)\le c_1 \left(\frac{x_d}{|x-y|}\right)^{\beta_1}\log\left(1+\frac{x_d\vee y_d}{x_d\wedge y_d}\right)^{\beta_3}.
\end{equation}
If $x_d/2\le y_d\le 2x_d$, then the log term above is estimated from above by a constant, implying \eqref{e:key-estimate-for-B}.

Note that for $a,b>0$ with $a\ge 2b$, it holds that $\log(1+a/b)\le 2\log(a/b)$. 
Using this and \eqref{e:former-remark-3} we see that if $y_d\ge 2x_d$ or $y_d\le x_d/2$, then 
\begin{align*}
\sB(x,y)\le  c_2\left(\frac{x_d \wedge y_d  }{|x-y|}\right)^{\beta_1} \left|\log\left(\frac{y_d}{x_d} \right) \right|^{\beta_3}
\le  c_3 \left(\frac{x_d}{|x-y|}\right)^{\beta_1}(|\log x_d|^{\beta_3}+|\log y_d|^{\beta_3})\, .
\end{align*}
We now give the  the proof  for the case $y_d\ge 2x_d$ and the case  $y_d\le x_d/2$ separately.

\noindent
\textbf{Case 1.} $y_d\ge 2x_d$:  
If $y_d<1$, then $|\log y_d|\le |\log x_d|$, so
$
\sB(x,y)\le 2c_3 x_d^{\beta_1}|\log x_d|^{\beta_3}|x-y|^{-\beta_1}\, .
$

If $y_d\ge 1$, we use the estimate $\log y_d\le \log |y|$ to get
$$
\sB(x,y)\le c_4x_d^{\beta_1}\left(|\log x_d|^{\beta_3}+(\log |y|)^{\beta_3}\right)|x-y|^{-\beta_1}\, .
$$
Therefore
\begin{eqnarray*}
\sB(x,y) &\le& c_5 x_d^{\beta_1}\left(|\log x_d|^{\beta_3}+{\bf 1}_{|y|\ge 1}(\log|y|)^{\beta_3}\right)|x-y|^{-\beta_1}\\
&\le &c_5 x_d^{\beta_1}(|\log x_d|^{\beta_3}\vee 1)\big(1+{\bf 1}_{|y|\ge1}(\log|y|)^{\beta_3}\big)|x-y|^{-\beta_1},
\end{eqnarray*}
which is precisely \eqref{e:key-estimate-for-B}.

\noindent 
\textbf{Case 2.} $y_d\le x_d/2$: It is straightforward to see that $f(t)=t^{\beta_1}((\log(1/t))^{\beta_3}\vee 1)$ is almost increasing on $(0,\infty)$. Therefore, there exists $c_6>0$ such that $y_d^{\beta_1}|\log y_d|^{\beta_3}\le c_6 x_d^{\beta_1}(|\log x_d|^{\beta_3}\vee 1)$. Hence,
$\sB(x,y)\le c_{7} x_d^{\beta_1}(|\log x_d|^{\beta_3}\vee 1) |x-y|^{-\beta_1}$, implying \eqref{e:key-estimate-for-B}.

\noindent (b) This follows immediately from part (a). \qed

\begin{lemma}\label{l:former-remark} 
Let $q<\alpha+\beta_1$.
\begin{itemize} 
\item[(a)] 
For all $x\in \R^d_+$, it holds that
\begin{align}\label{e:former-remark}
&\int_{\R^d_+, |y-x|\ge 2|x|}\frac{|y|^q}{|x-y|^{d+\alpha}}\sB(x,y)\, dy  
<  \infty\, .
\end{align}
\item[(b)]
For every $b>0$, there exists a constant $C_{21}=C_{21}(b)>0$ such that, for all $x\in \R^d_+$, it holds that
\begin{align}\label{e:former-remark-2}
&\int_{\R^d_+, |y|\ge (1+b)|x|}\frac{|y|^q}{|x-y|^{d+\alpha}}\sB(x,y)\, dy  \nonumber \\
&\le   C_{21}(b) x_d^{\beta_1}(|\log x_d|^{\beta_3}\vee 1)\int_{|y|\ge |x|}|y|^{-\beta_1+q-d-\alpha}\big(1+{\bf 1}_{|y|\ge1}(\log|y|)^{\beta_3}\big)\, dy
<  \infty\, .
\end{align}
\end{itemize}
\end{lemma}

\pf First note that since $q<\alpha+\beta_1$, it holds that
$$
\int_{|y|\ge |x|}|y|^{-\beta_1+q-d-\alpha}\big(1+{\bf 1}_{|y|\ge1}(\log|y|)^{\beta_3}\big)\, dy<  \infty\, .
$$
Further, if $|y-x|\ge |2|x|$ or $|y|\ge (1+b)|x|$, then $|y-x|\asymp |y|$. Now the conclusion for both (a) and (b) follows immediately from Lemma \ref{l:key-estimate-for-B}. \qed

For $f:\R^d_+\to \R$, we set for $x\in \R^d_+$,
\begin{align}\label{d:Lalpha}
L_{\alpha}^\sB f(x):=\textrm{p.v.} \int_{\R^d_+}\frac{f(y)-f(x)}{|y-x|^{d+\alpha}}\sB(x, y)\, dy=\textrm{p.v.}\int_{\R^d_+}(f(y)-f(x))J^{\R^d_+}(x,y)\, dy\, \, , 
\end{align}
whenever the principal value integral on the right hand side makes sense.

For $q>0$,
let 
$
g_q(y):=y_d^q=\delta_{\R^d_+}(y)^q\, .
$
Recall that $C(\alpha, p, \sB)$ is given in \eqref{e:explicit-C}.
The following is the analog of \cite[(5.4)]{BBC}.

\begin{lemma}\label{l:LB-on-g}
Let $p\in ((\alpha-1)_+, \alpha+\beta_1)$. Then
$C(\alpha, p, \sB)$ is finite and 
$$
L_{\alpha}^\sB g_p(x)= C(\alpha,p, \sB)x_d^{p-\alpha}, \quad x\in \R^d_+.
$$
\end{lemma}
\pf 
Recall  $\mathbf{e}_d=(\tilde{0}, 1)$. By  \textbf{(B8)},  we can for simplicity take $x=(\wt{0}, x_d)$.

By the change of variables $y=x_d z$,  and by using \textbf{(B8)}, we  have
\begin{align*}
L_{\alpha}^\sB g_p(x)=x_d^{p-\alpha} \textrm{p.v.} \int_{\R^{d-1}}\int_0^{\infty}\frac{z_d^p-1}{| (\wt{z}, z_d)-\mathbf{e}_d|^{d+\alpha}}\sB(\mathbf{e}_d,  (\wt{z}, z_d))\, 
dz_d d \wt{z}
=:x_d^{p-\alpha} I_1\, .
\end{align*}
Note that Lemma \ref{l:former-remark}(a)  takes care of integrability away of the point $x$.
Using the change of variables $\wt{z}=|z_d-1|\wt{u}$, we get
$
I_1= \int_{\R^{d-1}}{(|\wt{u}|^2+1)^{-(d+\alpha)/2}}I_2(\wt{u})\, d\wt{u},
$
where
 $$I_2(\wt{u})=I_2:=\lim_{\epsilon\to 0}\left(\int_0^{1-\epsilon}+\int_{1+\epsilon}^{\infty}\right)\frac{z_d^p-1}{|z_d-1|^{1+\alpha}} 
\sB\big(\mathbf{e}_d, (|z_d-1|\wt{u}, z_d)\big)dz_d .
$$
Using the change of variables $\wt{z}=|z_d-1|\wt{u}$, we get
$
I_1= \int_{\R^{d-1}}{(|\wt{u}|^2+1)^{-(d+\alpha)/2}}I_2(\wt{u})\, d\wt{u},
$
where
 $$I_2(\wt{u})=I_2:=\lim_{\epsilon\to 0}\left(\int_0^{1-\epsilon}+\int_{1+\epsilon}^{\infty}\right)\frac{z_d^p-1}{|z_d-1|^{1+\alpha}} 
\sB\big(\mathbf{e}_d, (|z_d-1|\wt{u}, z_d)\big)dz_d .
$$
Fix $\wt{u}$ and, by using the change of variables $s=1/z_d$, \textbf{(B1)} and \textbf{(B8)}, we get that  the second integral is equal to
\begin{align*}
&\int_0^{\frac{1}{1+\epsilon}}\frac{s^{\alpha-1-p}-s^{\alpha-1}}{(1-s)^{1+\alpha}}
\sB\big((1-s)\wt{u}, 1), s\mathbf{e}_d \big) ds\\
=&\int_0^{1-\epsilon}\frac{s^{\alpha-1-p}-s^{\alpha-1}}{(1-s)^{1+\alpha}}
\sB\big((1-s)\wt{u}, 1), s\mathbf{e}_d \big) ds +\int_{1-\epsilon}^{\frac{1}{1+\epsilon}}\frac{s^{\alpha-1-p}-s^{\alpha-1}}{(1-s)^{1+\alpha}}\sB\big(((1-s)\wt{u}, 1), s\mathbf{e}_d \big) ds\, .
\end{align*}
Therefore,
\begin{align*}
I_2
=&\lim_{\epsilon\to 0}\int_0^{1-\epsilon}\frac{(s^p-1)+(s^{\alpha-1-p}-s^{\alpha-1})}{(1-s)^{1+\alpha}}\sB\big(((1-s)\wt{u}, 1), s\mathbf{e}_d \big)\, ds \\
& + \lim_{\epsilon\to 0}\int_{1-\epsilon}^{\frac{1}{1+\epsilon}}\frac{s^{\alpha-1-p}-s^{\alpha-1}}{(1-s)^{1+\alpha}}\sB\big(((1-s)\wt{u}, 1), s\mathbf{e}_d \big)\, ds =: I_{21}+I_{22}\, .
\end{align*}
Next
\begin{align*}
I_{21}
=&\lim_{\epsilon\to 0}\int_0^{1-\epsilon} \frac{(1-s^p)(1-s^{p-(\alpha-1)})}{(1-s)^{1+\alpha}}s^{\alpha-1-p}\sB\big(((1-s)\wt{u}, 1), s\mathbf{e}_d \big)\, ds.
\end{align*}
Since $\sB$ is bounded (uniformly in $\wt{u}$) by \textbf{(B7)} and $\alpha<2$,
the fraction is integrable near 1. 
Further, if $s\in (0,1/2)$, then by using 
\textbf{(B7)}   and \eqref{e:B7_2} 
we obtain the following estimate
$$
\sB\big(((1-s)\wt{u}, 1), s\mathbf{e}_d \big)\le c_1 s^{\beta_1}\log(1/s)^{\beta_3}\, ,
$$
where the constant $c_1>0$ does not depend on $\wt{u}$. 
The fact that $p<\alpha+\beta_1$  implies that the function $s\mapsto s^{\alpha-1-p}s^{\beta_1}\log(1/s)^{\beta_3}$ is integrable near 0. Thus there exists a constant $c_2>0$ independent of $\wt{u}\in \R^{d-1}$ such that
\begin{align}
\label{e:I21n}
I_{21}=\int_0^1 \frac{(s^p-1)(1-s^{\alpha-p-1})}{(1-s)^{1+\alpha}}\sB\big(((1-s)\wt{u}, 1), s\mathbf{e}_d \big)\, ds<c_2 <\infty\,  .
\end{align}
Here we used that the integrand is non-negative since $p>\alpha-1$.
Now we show that $I_{22}=0$. Indeed, since $\sB$ is bounded, we conclude from \cite[p.121]{BBC} that
$$
\left|\int_{1-\epsilon}^{\frac{1}{1+\epsilon}}\frac{s^{\alpha-1-p}-s^{\alpha-1}}{(1-s)^{1+\alpha}} \sB\big((1-s)\wt{u}, 1), s\mathbf{e}_d \big)ds\right| \le c\epsilon^{2-\alpha}\, .
$$

Therefore,
\begin{align*}
I_1&=\int_{\R^{d-1}}\frac{1}{(|\wt{u}|^2+1)^{(d+\alpha)/2}}\left( \int_0^1 \frac{(s^p-1)(1-s^{\alpha-p-1})}{(1-s)^{1+\alpha}}
\sB\big(((1-s)\wt{u}, 1), s\mathbf{e}_d \big)\, ds\right)d\wt{u}\, .
\end{align*}
By combining the uniform (in $\wt{u}$) estimate for $I_{21}$ and the fact that $\int_{\R^{d-1}}(|\wt{u}|^2+1)^{(-d-\alpha)/2}d\wt{u}<\infty$, 
we see that $I_1=C(\alpha, p, \sB)\in (0,\infty)$. The proof is complete.
\qed

\begin{remark}\label{r:C-increasing}{\rm 
(a)
Note that $p\mapsto (s^p-1)(1-s^{\alpha-p-1})$ is increasing for $s\in (0,1)$. This implies that $p\mapsto C(\alpha, p, \sB)$ is also increasing. Further, it is clear
$\lim_{p\downarrow (\alpha-1)_+}C(\alpha, p, \sB)=0$ and, 
by inspecting the integrand in \eqref{e:I21n} near 0 and using Fatou's lemma, we see  $\lim_{p\uparrow \alpha+\beta_1}C(\alpha, p, \sB)=\infty$. 

\noindent (b) Note also that if $\alpha\in [1,2)$, then the above proof shows that $L^\sB_{\alpha} g_p=0$ for $p=\alpha-1$.
}
\end{remark}

Next we set
\begin{equation}\label{e:defn-LB}
L^\sB f(x):=L_{\alpha}^\sB f(x)- C(\alpha, p,\sB)x_d^{-\alpha}f(x)\, , \quad x\in \R^d_+\, ,
\end{equation}
where $L^\sB_\alpha$ is defined in \eqref{d:Lalpha}. Thus, $L^\sB g_p(x)=0$, $x\in \R^d_+$. 
Recall that $L^\sB $ was in a more general setting already introduced in \eqref{e:defn-LBD}.

\bigskip

For $a,b, r>0$ and $\wt{w} \in \R^{d-1}$, we define 
$$
D_{\wt{w}}(a,b):=\{x=(\wt{x}, x_d)\in \R^d:\, |\wt{x}-\wt{w}|<a, 0<x_d<b\}.
$$
Without loss of generality, we will  deal with the case $\wt w=\wt{0}$ only. 
We will write $D(a,b)$ for $D_{\wt{0}}(a,b)$ and
$U(r)=D_{\wt{0}}(\frac{r}2, \frac{r}2).$
Further  we use  $U$ for $U(1)$.

Let $h_q(x):=g_q(x){\bf 1}_{D(1,1)}(x)=x_d^q {\bf 1}_{D(1,1)}(x)$.

\begin{lemma}\label{l:estimate-of-L-hat-B}
\begin{itemize}
\item [(a)] Let $p\in ((\alpha-1)_+, \alpha+\beta_1)$. There exists $C_{22}>0$ such that for every $z\in U$, 
$$
0\ge L^\sB  h_p(z) \ge -C_{22} z_d^{\beta_1}|\log z_d|^{\beta_3}\, .
$$
\item [(b)] Let $(\alpha-1)_+ <p< q<\alpha+\beta_1$. There exist $r_0\in (0,1/2]$ and  $C_{23}>0$ and $C_{24}>0$ such that for every  $z\in D(\frac12, r_0)$,
$$
C_{23} z_d^{q-\alpha}\le L^\sB  h_q(z) \le C_{24} z_d^{q-\alpha}.
$$
\item [(c)] Let $(\alpha-1)_+ <q< p<\alpha+\beta_1$. 
There exist  $C_{25}>0$, $C_{26}>0$ and $C_{27}>0$ such that for every  $z\in U$,
$$
-C_{26} z_d^{q-\alpha}-C_{27} z_d^{\beta_1}|\log z_d|^{\beta_3}\le L^\sB  h_q(z) \le -C_{25} z_d^{q-\alpha}.
$$
\end{itemize}
\end{lemma}
\pf (a) Let $z\in U$. Then by using that $L^\sB  g_p(x)=0$, we see that
\begin{align*}
L^\sB  h_p(z)=-\int_{D(1,1)^c\cap \R^d_+}\frac{y_d^p}{|y-z|^{d+\alpha}}\sB (z,y)\, dy\, ,
\end{align*}
which is negative. 

Further, if $y\in D(1,1)^c$ and $z\in U=D(1/2,1/2)$, then 
$|y|\ge   1>|z|$, $|y-z|\ge z_d$ and $|y-z| \asymp |y|$. 
 Thus it follows from 
Lemma \ref{l:key-estimate-for-B} that
\begin{align*}
&|L^\sB  h_p(z)|\le 
\int_{D(1,1)^c\cap \R^d_+}
\frac{|y|^{p}}{|y-z|^{d+\alpha}}\sB (z,y)\, dy \\
&\le c_1\int_{y \in \R^d_+, |y|\ge  1, |y-z|\ge z_d}
|y|^{-d-\alpha+p}\sB (z,y)\, dy \le c_2 z_d^{\beta_1}|\log z_d|^{\beta_3}
\end{align*}

\noindent (b) Let $z\in U$. Then by using Lemma \ref{l:LB-on-g}, 
\begin{eqnarray*}
L^\sB  h_q(z)
=(C(\alpha, q, \sB )-C(\alpha, p, \sB ))z_d^{q-\alpha}-\int_{D(1,1)^c
\cap \R^d_+}\frac{y_d^q}{|y-z|^{d+\alpha}}\sB (z, y)\, dy\, .
\end{eqnarray*}
Since the second term is non-negative, by removing it we obtain the upper bound.
In the same way as before (this uses $q<\alpha+\beta_1$),
$$
\left| \int_{D(1,1)^c\cap \R^d_+}\frac{y_d^q}{|y-z|^{d+\alpha}}\sB (z, y)\, dy \right| \le c_2 z_d^{\beta_1}|\log z_d|^{\beta_3}\, .
$$
Thus, for any $z\in U$,
$$
L^\sB  h_q(z) \ge  (C(\alpha, q, \sB )-C(\alpha, p, \sB ))z_d^{q-\alpha}-c_2 z_d^{\beta_1}|\log z_d|^{\beta_3}\, .
$$
Since $q-\alpha<\beta_1$ and $C(\alpha, q, \sB )-C(\alpha, p, \sB )>0$, we can find $r_0\in (0,1/2]$ such that the function $t\mapsto (C(\alpha, q, \sB )-C(\alpha, p, \sB ))t^{q-\alpha}- c_2 t^{\beta_1}|\log t|^{\beta_3} $ 
is bounded from below by $c_3 t^{q-\alpha}$ with a positive constant $c_3>0$ for all $t\in (0,r_0)$. 
This concludes the proof of the lower bound. 

\noindent 
(c) This is almost the same as the proof of (b) except that now  $C(\alpha, q,\sB )-C(\alpha,p,\sB )<0$ and we keep the term with $z_d^{\beta_1}|\log z_d|^{\beta_3}$ in the lower bound (hence do not have to shrink). 
\qed

For $r\in (0, 1)$, let 
$C(r):=(D(1,1)\setminus D(\sqrt{r}, \sqrt{r}) )  \cap \{y\in \R^d_+: y_d\ge |\wt{y}|\}$ 
be the part of the cone in  $D(1,1)\setminus D(\sqrt{r}, \sqrt{r})$. Note that for $z\in U(r)$ and $y\in C(r)$ we have $|z-y|\le 2|y|$, $|y|\le \sqrt{2}y_d$, $|z|\le |y|/\sqrt{2}$, and therefore $y_d/4\le |y|/4\le  |z-y| \le 2\sqrt{2}y_d$. 
Now we find a lower bound for $\sB (z,y)$, $z\in U(r)$, $y\in C(r)$. By using \textbf{(B7)}, the lower bound in \eqref{e:B7_2}, 
$z_d<y_d$ and $y_d\asymp |z-y|\asymp |y|$, we get
\begin{eqnarray*}
\sB (z,y)&\ge & c_1 \left(\frac{z_d\wedge y_d}{|z-y|}\wedge 1\right)^{\beta_1}\left(\frac{z_d\vee y_d}{|z-y|}\wedge 1\right)^{\beta_2}\log\left(1+\frac{(z_d\vee y_d)\wedge |z_d-y_d|}{z_d\wedge y_d \wedge |z-y|}\right)^{\beta_3}\\
 &\ge & c_2 \left(\frac{z_d}{|y|}\right)^{\beta_1}\log(1+\frac{y_d}{z_d})^{\beta_3}\, .
\end{eqnarray*}
Since $y_d\ge \sqrt{r}$ and $z_d\le r/2$, we see that $y_d/z_d\ge 2/\sqrt{z_d}$, hence $\log(1+y_d/z_d)\ge \log(y_d/z_d)\ge \log(1/\sqrt{z_d})=|\log z_d|/2$. Thus,
$
\sB (z,y)\ge c_3 \left({z_d}/{|y|}\right)^{\beta_1}|\log z_d|^{\beta_3}\, .
$
By using the  L\'evy system formula we get for  $x\in U(r)$
(with constants $c_4$, $c_5$ \emph{not} depending on $r$),
\begin{align*}
\E_x[h_p(Y_{\tau_{U(r)}})]&\ge  \E_x[h_p(Y_{\tau_{U(r)}}), Y_{\tau_{U(r)}}\in C(r)]\\
&\ge  c_4 \E_x\int_{C(r)} \int_0^{\tau_{U(r)}} |Y_t-y|^{-d-\alpha} \left(\frac{Y_t^d}{|y|}\right)^{\beta_1} |\log Y_t^d|^{\beta_3} y_d^p \, dy \, dt\\
&\ge  c_5 \int_{C(r)}{y_d^p}{|y|^{-d-\alpha-\beta_1}}\, dy \left(\E_x \int_0^{\tau_{U(r)}} (Y_t^d)^{\beta_1} |\log Y_t^d|^{\beta_3} \, dt\right)\, .
\end{align*}
Since $\int_{D(1,1)}y_d^p |y|^{-d-\alpha-\beta_1}\, dy =+\infty$, we can choose $r_0$ small enough so that 
$$c_5 \int_{C(r)}{y_d^p}{|y|^{-d-\alpha-\beta_1}}\, dy\ge 1,
\qquad \mbox{ for all } r\in (0, r_0].
$$

Hence, for all $r\in (0, r_0]$
\begin{equation}\label{e:integral-upper-g} 
\E_x \int_0^{\tau_{U(r)}} (Y_t^d)^{\beta_1} |\log Y_t^d|^{\beta_3} \, dt \le \E_x[h_p(Y_{\tau_{U(r)}})]\, , \quad 
x\in U(r) \, .
\end{equation}

In the remainder of this section, we will, in several places, 
need Lemma \ref{l:dynkin-hp} (which uses the additional assumption \textbf{(B4)}), an extension of Proposition 
\ref{p:Dynkin} to functions which are neither smooth nor of compact support. 
Since its proof is  rather technical, in order not to interrupt the flow of the presentation, we postpone it to Section \ref{s:10}.

\emph{In the rest of this section we assume that 
\textbf{(B4)} also holds.}

\begin{lemma}\label{l:upper-bound-for-integral}
(a) For all $r\in (0, r_0]$ and $x\in U(r)$,
\begin{equation}\label{e:upper-bound-for-integral1}
\E_x \int_0^{\tau_{U(r)}} (Y_t^d)^{\beta_1}  |\log Y_t^d|^{\beta_3}\, dt \le  x_d^p\, . 
\end{equation}

\noindent (b)
There exist a constant $C_{28}>0$ such that for
all $x\in U$,
\begin{equation}\label{e:upper-bound-for-integral}
\E_x \int_0^{\tau_{U}} (Y_t^d)^{\beta_1}  |\log Y_t^d|^{\beta_3}\, dt \le C_{28} x_d^p\, .
\end{equation}
\end{lemma}
\pf
(a) 
By Lemma \ref{l:dynkin-hp} it holds that 
$$
\E_x[h_p(Y_{\tau_{U(r)}})]=x_d^p +\E_x\int_0^{\tau_{U(r)}}L^\sB  h_p (Y_s)\, ds \le x_d^p\, ,
$$
where the inequality follows from Lemma \ref{l:estimate-of-L-hat-B} (a). The required estimate is a consequence of \eqref{e:integral-upper-g}.

\noindent
(b) Let $\tau^{(r_0)}_{r_0U}$ be the first time that the process $Y^{(r_0)}=(Y^{(r_0), 1}, \dots ,Y^{(r_0), d})$ 
exits $r_0U=U(r_0)$. Then it is easy to check by definition that $\tau^{(r_0)}_{r_0U}=r_0^{\alpha}\tau_{U}$. Using this, one can easily see that 
\begin{eqnarray}\label{e:scaling}
\E_{r_0x}\int^{\tau^{(r_0)}_{r_0U}}_0
(Y^{(r_0), d}_s)^{\beta_1} |\log (Y^{(r_0), d}_s)|^{\beta_3}\, ds &= & r_0^{\alpha+\beta_1}\E_x\int^{\tau_U}_0(Y^d_s)^{\beta_1} |\log (r_0 Y_s^{d})|^{\beta_3}\, ds\nonumber \\
&\ge &r_0^{\alpha+\beta_1}\E_x\int^{\tau_U}_0(Y^d_s)^{\beta_1} |\log Y_s^{d}|^{\beta_3}\, ds.
\end{eqnarray}
In the last line we used that since $r_0<1$, for all $t\in (0,1/2) $ it holds that $|\log (r_0 t)|\ge |\log t|$.
Hence, for any  $x\in U$ we have
\begin{eqnarray*}
\E_x \int^{\tau_U}_0(Y^d_s)^{\beta_1} |\log Y_s^{d}|^{\beta_3}\,  ds&\le&r_0^{-\alpha-\beta_1}\E_{r_0x}\int^{\tau^{(r_0)}_{r_0U}}_0(Y^{(r_0), d}_s)^{\beta_1} |\log (Y^{(r_0), d}_s)|^{\beta_3}\, ds \\
 &\le &  r_0^{-\alpha-\beta_1}(r_0 x_d)^p= r_0^{-\alpha-\beta_1+p}x_d^p=C_{28} x_d^p 
\end{eqnarray*}
with $C_{28}=r_0^{-\alpha-\beta_1+p}$. 
\qed

Next we want to prove the opposite inequality. The key step towards this goal is the following lemma.

\begin{lemma}\label{l:key-lemma}
Let $p\in ((\alpha-1)_+, \alpha+\beta_1)$ 
and assume that $ \theta >(\alpha-1)_+$.
\begin{itemize}
\item[(a)] There exist a function $\psi:\R^d\to [0,\infty)$ and 
a constant $C_{29}>0$ such that
$$
L^\sB  \psi(x)\le C_{29} x_d^{\beta_1}\,  
|\log x_d|^{\beta_3}\, ,\quad x\in U
$$
and the following assertions hold: 
\item[(b)] The function $\phi(x):=h_p(x)-\psi(x)$, $x\in \R^d_+$, satisfies the following properties: 
	\begin{itemize}
	\item[(b1)]  $\phi(x)=x_d^p$ for all $x=(\wt{0}, x_d)\in U$ with $0<x_d<1/4$;
	\item[(b2)]  $\phi(x)\le 0$ for all $x\in U^c\cap \R^d_+$;
	\item[(b3)]  There exists $C_{30}>0$ such that $L^\sB \phi(x)\ge -C_{30} x_d^{\beta_1}\,  |\log x_d|^{\beta_3}$ for all $x\in U$.
	\end{itemize}
\item[(c)] Let $(\alpha-1)_+<q<p$. The function $\varphi(x):=h_q(x)-\psi(x)$ satisfies the following properties:
\begin{itemize}
	\item[(c1)]  $\varphi(x)=x_d^q$ for all $x=(\wt{0}, x_d)\in U$ with $0<x_d<1/4$;
	\item[(c2)]  $\varphi(x)\le 0$ for all $x\in U^c\cap \R^d_+$;
	\item[(c3)] There exist $C_{31}>0$ and $C_{32}>0$ such that	
	$L^\sB \varphi(x)\ge -C_{31} x_d^{q-\alpha}-C_{32} x_d^{\beta_1}\,  |\log x_d|^{\beta_3}$ for all $x\in U$.
	\end{itemize}
\end{itemize}
\end{lemma}

The proof of  the lemma is rather involved, we  defer it to the Section \ref{s:8}.

\bigskip

Combining Lemma \ref{l:key-lemma} (b1) with a formal application of Proposition \ref{p:Dynkin}
to $\phi$, we can get the following inequality opposite to \eqref{e:upper-bound-for-integral}.
However, since $\phi$ is not smooth, we can not apply
Proposition \ref{p:Dynkin} to $\phi$ directly. Therefore we postpone its proof to Section \ref{s:10}.

\begin{lemma}\label{l:lower-bound-for-integral}
Let $p\in ((\alpha-1)_+, \alpha+\beta_1)$.
For any $x=(\wt{0}, x_d)$ with $0<x_d<1/4$ it holds that
\begin{equation}\label{e:lower-bound-for-integral}
\E_x \int_0^{\tau_{U}} (Y_t^d)^{\beta_1}\,  |\log Y_t^d|^{\beta_3}\, dt \ge C_{30}^{-1} x_d^p\, ,
\end{equation}
where $C_{30}$ is the constant from Lemma \ref{l:key-lemma}.
\end{lemma}

Note that there exists $c>0$ such that $z^{\beta_1}|\log z|^{\beta_3}\le c$ for all $z\in (0, 1)$. Since $Y_t^d\in (0,1)$ for $t<\tau_U$,
\eqref{e:lower-bound-for-integral} immediately implies the following estimate for the expected exit time from $U$:
\begin{equation}\label{e:exit-time-estimate-U-upper}
\E_x \tau_U\ge c^{-1}C_{30}^{-1} x_d^p\, , \quad
x=(\wt{0}, x_d) \text{ with }  0<x_d<1/4.
\end{equation}

One can easily obtain the following exit probability estimate by combining Lemmas \ref{l:upper-bound-for-integral} and \ref{l:lower-bound-for-integral}.

\begin{lemma}\label{l:exit-probability-estimate}
Let $p\in ((\alpha-1)_+, \alpha+\beta_1)$.
There exists  $C_{33}>1$ such that 
for $x\in D(1/4,1/4)$,
\begin{equation}\label{e:exit-probability-estimateu}
 \P_x(Y_{\tau_{U}}\in D(1/2,1)\setminus D(1/2,3/4))\le C_{33} x_d^p\, ,
\end{equation}
and for $x=(\wt{0}, x_d)\in D(1/8,1/8)$,
\begin{equation}\label{e:exit-probability-estimated}
\P_x(Y_{\tau_{D(1/4,1/4)}}\in D(1/4,1)\setminus D(1/4,3/4)) \ge C_{33}^{-1} x_d^p \, .
\end{equation}
\end{lemma}
\pf  For $y\in U=D(1/2,1/2)$ and $z\in D(1/2,1)\setminus D(1/2,3/4)$, it holds that $y_d<z_d$, $|y-z|\asymp z_d$ and $y_d<2|y-z|$. Hence, 
$
\sB (y,z) \asymp \left({y_d}/{|y-z|}\right)^{\beta_1}
  \log \left(1+({z_d}/{y_d})\right)^{\beta_3}\, .
$
Since $1\ge z_d \ge (3/2)y_d$, we have that $\log(1/y_d)\le \log(z_d/y_d)\le \log(1+z_d/y_d)\le c \log(z_d/y_d)\le c\log(1/y_d)$. We conclude that
$
\sB (y,z)\asymp \left({y_d}/{|y-z|}\right)^{\beta_1} \, |\log y_d|^{\beta_3}.
$
 Hence, by using the L\'evy system formula, we get
\begin{align*}
&\P_x\left(Y_{\tau_{U}}\in D(1/2, 1)\setminus D(1/2, 3/4)\right)\nonumber\\
&\le c_1\E_x\int^{\tau_{U}}_0 (Y^d_t)^{\beta_1}\, |\log Y_t^d|^{\beta_3}\int_{D(1/2, 1)\setminus D(1/2, 3/4)}\frac1{|Y_t-z|^{d+\alpha+\beta_1}}\, dz\, dt\nonumber\\
&\le c_2\E_x\int^{\tau_{U}}_0 (Y^d_t)^{\beta_1}\, |\log Y_t^d|^{\beta_3}\, dt
\end{align*}
and 
\begin{align*}
&\P_x\left(Y_{\tau_{D(1/4,1/4)}}\in D(1/4, 1)\setminus D(1/4, 3/4)\right)\nonumber\\
&\ge c_3\E_x\int^{\tau_{D(1/4,1/4)}}_0 (Y^d_t)^{\beta_1}\, |\log Y_t^d|^{\beta_3}\int_{D(1/4, 1)\setminus D(1/4, 3/4)}\frac1{|Y_t-z|^{d+\alpha+\beta_1}}\, dz\, dt\nonumber\\
&\ge c_4\E_x\int^{\tau_{D(1/4,1/4)}}_0 (Y^d_t)^{\beta_1}\, |\log Y_t^d|^{\beta_3}\, dt\, .
\end{align*}
The claim now follows from Lemmas \ref{l:upper-bound-for-integral} (b) and \ref{l:lower-bound-for-integral}. \qed

\begin{lemma}\label{l:lower-bound-for-integral-q}
Let $(\alpha-1)_+<p<\alpha+\beta_1$ and $\alpha+\beta_2<p$.
There exist $r_1>0$ and $C_{34}>0$
such that for all  $r>0$ and 
$x=(\wt{0}, x_d)$ with $0<x_d<r_1r$,  it holds that
\begin{equation}\label{e:lower-bound-for-integral-q}
\E_x \int_0^{\tau_{U(r)}} (Y_t^d)^{\beta_2} \, dt \ge C_{34} x_d^{\alpha+\beta_2}\, .
\end{equation}
\end{lemma}
\pf  
Let $q:=\alpha+\beta_2$. Note that $q<p$ by assumption. By combining Lemma \ref{l:key-lemma} (c3) with a formal application of Proposition \ref{p:Dynkin} to $\varphi$, we can get the assertion of this lemma. 
However, since $\varphi$ is not smooth, we can not apply
Proposition \ref{p:Dynkin} to $\varphi$ directly. Therefore,
the rigorous proof is postponed to Section \ref{s:10}.
 \qed
 
\begin{lemma}\label{l:exit-time estimate-U_new}
Let $(\alpha-1)_+ <p<\alpha+\beta_1$ and $\alpha+\beta_2<p$.
There exist a constant $C_{35}>0$  such that 
for all $r>0$ and $x\in U(r)$,
\begin{equation}\label{e:exit-time-estimate-U_new}
\E_x\int_{0}^{\tau_{U(r)}}
  (Y_t^d)^{\beta_2}dt  \le C_{35} x_d^{ \alpha+\beta_2}\, .
\end{equation}
\end{lemma}
\pf 
Let $q=\alpha+\beta_2$. Note that $q<p$ by assumption. By combining Lemma \ref{l:estimate-of-L-hat-B} with a formal application of Proposition \ref{p:Dynkin} to $h_q-h_p$, we can get
the assertion of this lemma. 
However, since $h_q-h_p$ is not smooth, we can not apply
Proposition \ref{p:Dynkin} to $h_q-h_p$ directly. Therefore, 
the rigorous proof is postponed to Section \ref{s:10}.
\qed

The following counterpart of \eqref{e:exit-time-estimate-U-upper} is considerably simpler than the previous results, but has an additional restriction on the range of $p$.
\begin{lemma}\label{l:exit-time estimate-U}
Suppose that $p\in ((\alpha-1)_+, \alpha)$. 
There exist a constant $C_{36}>0$  such that
\begin{equation}\label{e:exit-time-estimate-U}
\E_x \tau_U \le C_{36} 
x_d^p\,, \quad x\in U.
\end{equation}
\end{lemma}
\pf Choose $q\in (p, \alpha)$ and let
$
\eta(x):=h_p(x)-h_q(x)$, $x\in \R^d_+.$
By combining Lemma \ref{l:estimate-of-L-hat-B} with a formal application of Proposition \ref{p:Dynkin} to $\eta$, we can get
the assertion of this lemma. 
However, since $\eta$ is not smooth, we can not apply
Proposition \ref{p:Dynkin} to $\eta$ directly. Therefore 
The rigorous proof is postponed to Section \ref{s:10}.
 \qed

\section{Boundary Harnack principle}\label{s:bhp}
In this section we give the proof of Theorem \ref{t:BHP}.

\begin{lemma}\label{e:POTAe7.14}
There exists $C_{37}>0$ such that for all 
$0<4r\le R\le 1$ and $w\in D(r,r)$, 
$$
\P_w\Big(Y_{\tau_{B(w, r)\cap \R^d_+}}\in A(w, R, 4)\cap \R^d_+\Big)\le C_{37}
\frac{r^{\alpha+\beta_1}}{R^{\alpha+\beta_1}}\frac{w_d^p}{r^p}.
$$
\end{lemma}

\pf
Let $y\in B(w,r)\cap \R^d_+$ and $z\in A(w, R, 4)\cap \R^d_+$. 
By Lemma \ref{l:key-estimate-for-B}(b)  we have that
$$
\sB (y,z)\le c_1 y_d^{\beta_1}\, |\log y_d|^{\beta_3} |y-z|^{-\beta_1} \le c_2 y_d^{\beta_1}\, |\log y_d|^{\beta_3} |w-z|^{-\beta_1}.
$$
Thus for all $0<4r\le R\le 1$ and $w\in D(r,r)$ ,
\begin{align*}
&\P_w\Big(Y_{\tau_{B(w, r)\cap \R^d_+}}\in A(w, R, 4)\cap \R^d_+\Big)\\
&\le c_2 \E_w\int^{\tau_{ B(w, r)\cap \R_d^+}}_0(Y_t^d)^{\beta_1}\, |\log Y_t^d|^{\beta_3} dt \int_{B(w, R)^c}\frac1{|w-z|^{d+\alpha+\beta_1}}dz\\
&\le c_3 R^{-\alpha-\beta_1}\E_w\int^{\tau_{ B(w, r)\cap \R_d^+}}_0(Y_t^d)^{\beta_1}\,
|\log Y_t^d|^{\beta_3}dt.
\end{align*}
Since $B(w, r)\cap \R_d^+\subset D(2r,2r)$,
applying
Lemma \ref{l:upper-bound-for-integral} and scaling, we get that there 
exists $c_4>0$ 
such that for all  $0<4r\le R\le 1$ and $w\in D(r,r)$,
\begin{align*}
&\P_w\Big(Y_{\tau_{B(w, r)\cap \R^d_+}}\in A(w, R, 4)\cap \R^d_+\Big)\\
&\le  c_3 R^{-\alpha-\beta_1}\E_w\int^{\tau_{ D(2r,2r)}}_0(Y_t^d)^{\beta_1}\, |\log Y_t^d|^{\beta_3} dt\le c_4 2^{\alpha+\beta_1-p} \frac{r^{\alpha+\beta_1}}{R^{\alpha+\beta_1}}\frac{w_d^p}{r^p}.
\end{align*}
\qed

\begin{lemma}\label{l:POTAl7.4}
There exists $C_{38}>0$  such that 
for any $x \in D(2^{-5}, 2^{-5})$,
$$
\P_x\left(Y_{\tau_{U}}\in D(1, 1)\right)\le C_{38} x_d^p.
$$
\end{lemma}

\pf 
Let 
$$
H_2:=\{Y_{\tau_{U}}\in D(1, 1)\}, \quad H_1:=\{Y_{\tau_{U}}\in D(1/2, 1)\setminus D(1/2, 3/4)\}.
$$
We claim that $\P_x(H_2)\le c\P_x(H_1)$ for all $x\in D(2^{-5}, 2^{-5})$. 
Combining this claim with
\eqref{e:exit-probability-estimateu} 
we arrive at the conclusion of this
lemma.
Now we prove the claim. Note that by \eqref{e:exit-probability-estimated},
\begin{equation}\label{e:POTAe7.15}
\P_w(H_1)\ge\P_w(Y_{\tau_{D_{\wt w}(1/4, 1/4)}}\in D_{\wt w}(1/4, 1)\setminus D_{\wt w}(1/4, 3/4))\ge c_1 w^p_d, \quad 
w\in D(2^{-3}, 2^{-3}).
\end{equation}
For $i\ge 1$, set
$$
s_0=s_1, \quad s_i=\frac18\Big(\frac12-\frac1{50}\sum^i_{j=1}\frac1{j^2}\Big)\quad \text{ and }  \quad J_i=D(s_i, 2^{-i-3})\setminus D(s_i, 2^{-i-4}) .
$$
Note that $1/(20)<s_i<1/(16)$. Define for $i\ge 1$,
\begin{equation}\label{e:POTAe7.16}
d_i=\sup_{z\in J_i}\frac{\P_z(H_2)}{\P_z(H_1)}, 
\quad \widetilde{J}_i=D(s_{i-1}, 2^{-i-3}),
\quad \tau_i=\tau_{\widetilde{J}_i}.
\end{equation}
Repeating the argument leading to \cite[(6.29)]{KSV18-a}, we get that for 
$z\in J_i$ and $i\ge 2$,
\begin{equation}\label{e:POTAe7.17}
\P_z(H_2)\le \Big(\sup_{1\le k\le i-1}d_k\Big)\P_z(H_1)
+\P_z\left(Y_{\tau_i}\in D(1, 1)\setminus \cup^{i-1}_{k=1}J_k\right).
\end{equation}

For  $i\ge 2$, define $\sigma_{i,0}=0, \sigma_{i,1}=\inf\{t>0: |Y_t-Y_0|\geq 2^{-i-3}\} $ and
$\sigma_{i,m+1}=\sigma_{i,m}+\sigma_{i,1}\circ\theta_{\sigma_{i,m}}$
for $m\geq 1$. By Lemma \ref{lower bound},  we have that there exists $k_1 \in (0,1)$ such that 
\begin{align}\label{e:POTAe7.18}
\P_{w}(Y_{\sigma_{i,1}}\in \widetilde{J}_i) \le 1-
\P_{w}( \sigma_{i,1}=\zeta ) \le 1-\P_{w}(
\tau_{B(w,\delta_D(w)/2)}
=\zeta )  <k_1,\ \ \ 
w\in  \widetilde{J}_i.
\end{align}
For the purpose of further estimates, we now choose a
positive integer $l$ 
such that $k_1^l\le 2^{-(\alpha+\beta_1)}$.
Next we choose 
$i_0 \ge 2$ large enough so that
$2^{-i}<1/(200 l i^3)$ for all $i\ge i_0$.
Now we assume $i\ge i_0$.
Using \eqref{e:POTAe7.18} and the strong Markov property we have
that  for   $z\in J_i$,
\begin{align}\label{e:POTAe7.19}    
&\P_z( \tau_{i}>\sigma_{i,li})\leq \P_z(Y_{\sigma_{i,k}}\in \widetilde{J}_i, 1\leq k\leq
li
 )\nonumber\\
 &=
 \E_z \left[ \P_{Y_{\sigma_{i,li-1}}} (Y_{\sigma_{i,1}}\in  \widetilde{J}_i) : Y_{\sigma_{i,li-1}} \in  \widetilde{J}_i,  Y_{\sigma_{i,k}}\in \widetilde{J}_i, 1\leq k\leq
li-2
  \right]
 \nonumber\\      
 &\leq \P_z\left(Y_{\sigma_{i,k}}\in \widetilde{J}_i, 1\leq k\leq
li-1
 \right)k_1\leq k_1^{li}.
\end{align}
Note that if $z\in J_i$ and $y\in D(1, 1) \setminus[ \widetilde{J}_i \cup(\cup_{k=1}^{i-1}J_k)]$, 
then $|y-z|\ge (s_{i-1}-s_i) \wedge (2^{-4}-2^{-i-3})  = 1/(400 i^2)$.
Furthermore, 
since  $2^{-i-3} < 1/(400 i^2)$ (recall that $i\ge  i_0$), 
if $Y_{\tau_i}(\omega)\in D(1, 1)\setminus \cup_{k=1}^{i-1}J_k$ and $\tau_i(\omega)\le \sigma_{i,li}(\omega)$, then $\tau_i(\omega)=\sigma_{i,k}(\omega)$ for some $k=k(\omega)\le li$.  
Dependence of $k$ on $\omega$ will be omitted  in the next few lines. 
Hence on
$\{Y_{\tau_{i}}\in D (1 , 1) \setminus \cup_{k=1}^{i-1}J_k,\ \ \tau_{i}\leq
\sigma_{i,li}\}$ with $Y_0=z\in J_i$,  we have
$|Y_{\sigma_{i,k}}-Y_{\sigma_{i,0}}|=|Y_{\tau_i}-Y_0|>
\frac{1}{400i^2}$ for some $1\leq k\leq li$.
Thus  for some $1\leq k\leq li$,
$
\sum_{j=0}^k|Y_{\sigma_{i,j}}-Y_{\sigma_{i,j-1}}|>
(400i^2)^{-1}
$
which implies for some $1\leq j\leq k\le  li$,
$ |Y_{\sigma_{i,j}}-Y_{\sigma_{i,j-1}}|\geq
({k400i^2})^{-1}\ge (li)^{-1} (400i^2)^{-1}\, .
$
Thus, we have 
\begin{align*}
& \{Y_{\tau_{i}}\in D(1, 1) \setminus
\cup_{k=1}^{i-1}J_k,\ \ \tau_{i}\leq \sigma_{i,li}\}\\
  \subset  &  \cup_{j=1}^{li}\{|Y_{\sigma_{i,j}}- Y_{\sigma_{i,j-1}}|\geq 1/(800li^3),
  Y_{\sigma_{i,j}}\in D(1, 1), 
  Y_{\sigma_{i,j-1}}\in \widetilde{J}_{i}
\}.
\end{align*}
Now, using  Lemma \ref{e:POTAe7.14} (with $r=2^{-i-3}$ and $R=1/(800 l i^3)$) 
 (noting that  $4 \cdot 2^{-i-2}<1/(400 l i^3)$ for all $i\ge i_0$), 
 and repeating the argument leading to \cite[(6.34)]{KSV18-a}, 
 we get that for $z\in J_i$,
\begin{align*}
&\P_z \left(Y_{\tau_{i}}\in D(1, 1) \setminus
\cup_{k=1}^{i-1}J_k,\ \ \tau_{i}\leq \sigma_{i,li} \right) \leq  li \sup_{z\in \widetilde{J}_{i}} 
  \P_z\left(|Y_{\sigma_{i,1}}-z  |\geq (800li^3)^{-1},  Y_{\sigma_{i,1}}\in D(1, 1)\right)\\
&\le \P_z\left(4>|Y_{\sigma_{i,1}}-z  |\geq (800li^3)^{-1} \right)\le c_{12}li  \left(
\frac{800li^3}{2^{i+3}}\right)^{\alpha+\beta_1}.
\end{align*}
By this and (\ref{e:POTAe7.19}),  we have for 
$z\in J_i$,  $i\ge i_0$,
\begin{align}\label{e:POTAe7.22}
&\P_z\left( Y_{\tau_{i}}\in D(1, 1) \setminus \cup_{k=1}^{i-1}J_k \right) \leq k_1^{li} 
+c_{2} li  \left(
\frac{800li^3 }{2^{i+3}}\right)^{\alpha+\beta_1}. 
\end{align} 
By our choice of $l$, we have 
\begin{align}\label{e:POTAe7.23}
&li    \left(
\frac{800li^3}{2^{i+3}}\right)^{\alpha+\beta_1} 
= 100^{\alpha+\beta_1} l^{1+\alpha+\beta_1} i^{1+3(\alpha+\beta_1)} \left(2^{-(\alpha+\beta_1)}\right)^i
\ge \left(2^{-(\alpha+\beta_1)}\right)^i
 \ge (k_1^{l})^{i}. 
\end{align}
Thus combining \eqref{e:POTAe7.23} with \eqref{e:POTAe7.22}, and then 
 using \eqref{e:POTAe7.15},  we get that  for
$z\in J_i$, $i\ge i_0$, 
\begin{align}\label{e:POTAe7.24}
&\frac{\P_z( Y_{\tau_i}\in
D(1, 1) \setminus \cup_{k=1}^{i-1}J_k)}{\P_z(H_1)} \le 
c_{3} li  
2^{ip}
\left(
\frac{800li^3 }{2^{i+3}}\right)^{\alpha+\beta_1}
\le  
c_{4} i ^{1+3(\alpha+\beta_1)}2^{(p-\alpha-\beta_1)i}.
\end{align} 
By this and (\ref{e:POTAe7.17}),  for 
$z\in J_i$,  $i\ge i_0$,for all $i\ge i_0$
\begin{align*}  
&\frac{\P_z( H_2)}{\P_z(H_1)} \leq  \sup_{1\leq k\leq i-1}  d_k   +\frac{\P_z( Y_{\tau_i}\in
D(1 , 1) \setminus \cup_{k=1}^{i-1}J_k)}{\P_z(H_1)}\le \sup_{1\leq k\leq i-1}d_k+
c_{4} \frac{i ^{1+3(\alpha+\beta_1)}}{2^{(\alpha+\beta_1-p)i}}.
\end{align*}
This implies that for all $i\ge1$
\begin{align} 
d_i& \leq  \sup_{1\leq k\leq i_0-1} d_k
+c_{4}\sum_{k=1}^i\frac{i ^{1+3(\alpha+\beta_1)}}{2^{(\alpha+\beta_1-p)i}}
\leq
\sup_{1\leq k\leq i_0-1} d_k
+c_{4}\sum_{k=1}^\infty\frac{i ^{1+3(\alpha+\beta_1)}}{2^{(\alpha+\beta_1-p)i}}
=:c_{5} <\infty.\nonumber
\end{align}
Thus the claim above is valid, since $D(2^{-5}, 2^{-5} ) 
\subset \cup_{k=1}^\infty J_k$. The proof is now complete.
\qed

\medskip
\noindent
\textbf{Proof of Theorem \ref{t:BHP}.}
By scaling, we just need to consider the case $r=1$. Moreover, by the Harnack inequality, the continuity of harmonic functions and a standard chain argument, it suffices to prove \eqref{e:TAMSe1.8} 
for  $x, y\in D_{\wt w}(2^{-8}, 2^{-8})$.

Since $f$ is harmonic in $D_{\wt w}(2, 2)$ and vanishes continuously on $B(\wt w, 2)\cap \partial \R^d_+$,
it is regular harmonic in $D_{\wt w}(7/4, 7/4)$ and vanishes continuously on $B(\wt w, 7/4)\cap \partial \R^d_+$.
Throughout the remainder of this proof, we assume that 
$x\in D_{\wt w}(2^{-8}, 2^{-8})$. 
Without loss of generality we take $\widetilde{w}=0$.

Define $x_0=(\widetilde{x}, 1/(16))$. By the Harnack inequality 
and Lemma \ref{l:exit-probability-estimate}, we have
\begin{align}\label{e:TAMSe6.37}
f(x)&=\E_x[f(Y_{\tau_{D(1/2, 1/2)}})]\ge \E_x[f(Y_{\tau_{D(1/2, 1/2)}}); Y_{\tau_{D(1/2, 1/2)}}\in D(1/2, 1)\setminus D(1/2, 3/4)]\nonumber\\
&\ge c_1f(x_0)\P_x(Y_{\tau_{D_{\wt x}(1/4, 1/4)}}\in D_{\wt x}(1/4, 1)\setminus D(1/4, 3/4))\ge c_2  f(x_0)  x^p_d.
\end{align}

\noindent
(a) {\it The case $\beta_1=\beta_2=\beta\ge0$, $\beta_3= \beta_4  =0$.}  
For $z\in U$ and $y\in \R^d_+\setminus D(1, 1)$, we have $|y-z|\asymp |y|$. Thus, by using \eqref{e:beta1=beta2},
\begin{align*}
&\int_{\R^d_+\setminus D(1, 1)}\frac{f(y)}{|z-y|^{d+\alpha}}\left(\frac{z_d}{|z-y|}\wedge 1\right)^\beta\left(\frac{y_d}{|z-y|}\wedge 1\right)^\beta dy\nonumber\\
&\asymp \int_{\R^d_+\setminus D(1, 1)}\frac{f(y)}{|y|^{d+\alpha}}\left(\frac{z_d}{|y|}\wedge 1\right)^\beta\left(\frac{y_d}{|y|}\wedge 1\right)^\beta dy\asymp z_d^\beta \int_{\R^d_+\setminus D(1, 1)}\frac{f(y)}{|y|^{d+\alpha+\beta}}\left(\frac{y_d}{|y|}\wedge 1\right)^\beta dy.
\end{align*}
Hence
\begin{align*}
\E_x\left[f(Y_{\tau_{U}}); Y_{\tau_{U}}\notin D(1, 1)\right]\asymp \E_x\int^{\tau_{U}}_0(Y^d_t)^\beta dt\int_{\R^d_+\setminus D(1, 1)}\frac{f(y)}{|y|^{d+\alpha+\beta}}\left(\frac{y_d}{|y|}\wedge 1\right)^\beta dy.
\end{align*}
Combining the above with Lemmas \ref{l:upper-bound-for-integral}.(b) and \ref{l:lower-bound-for-integral}  we get
\begin{align}\label{e:TAMSe6.41}
&\E_x\left[f(Y_{\tau_{U}}); Y_{\tau_{U}}\notin D(1, 1)\right]\asymp x_d^p \int_{\R^d_+\setminus D(1, 1)}\frac{f(y)}{|y|^{d+\alpha+\beta}}\left(\frac{y_d}{|y|}\wedge 1\right)^\beta dy.
\end{align}

On the other hand, by Theorem \ref{t:carleson} and Lemma \ref{l:POTAl7.4}, we have
\begin{align}\label{e:TAMSe6.42}
&\E_x\left[f(Y_{\tau_{U}}); Y_{\tau_{U}}\in D(1, 1)\right]\le c_3 f(x_0)\P_x(Y_{\tau_{U}}\in D(1, 1))\le c_4f(x_0)x^p_d.
\end{align}
Combining \eqref{e:TAMSe6.37}, \eqref{e:TAMSe6.41} and \eqref{e:TAMSe6.42}, we get
\begin{align*}
&f(x)=\E_x\left[f(Y_{\tau_{U}}); Y_{\tau_{U}}\notin D(1, 1)\right]+ \E_x\left[f(Y_{\tau_{U}}); Y_{\tau_{U}}\in D(1, 1)\right]\nonumber\\
&\le c_5x^p_d\left(f(x_0)+\int_{\R^d_+\setminus D(1, 1)}\frac{f(y)}{|y|^{d+\alpha+\beta}}\left(\frac{y_d}{|y|}\wedge 1\right)^\beta dy\right)
\end{align*}
and
\begin{align*}
f(x)&\ge \frac12 \E_x\left[f(Y_{\tau_{U}}); Y_{\tau_{U}}\in D(1/2, 1)\setminus D(1/2, 3/4)\right] + \frac12\E_x\left[f(Y_{\tau_{U}}); Y_{\tau_{U}}\notin D(1, 1)\right]\nonumber\\
&\ge c_6x^p_d\left(f(x_0)+\int_{\R^d_+\setminus D(1, 1)}\frac{f(y)}{|y|^{d+\alpha+\beta}}\left(\frac{y_d}{|y|}\wedge 1\right)^\beta dy\right).
\end{align*}
Therefore
\begin{equation}\label{e:TAMSe6.45}
f(x)\asymp x^p_d\left(f(x_0)+\int_{\R^d_+\setminus D(1, 1)}\frac{f(y)}{|y|^{d+\alpha+\beta}}\left(\frac{y_d}{|y|}\wedge 1\right)^\beta dy\right).
\end{equation}
For any $y\in D(2^{-8}, 2^{-8})$, 
we have the same estimate with $f(y_0)$ instead of $f(x_0)$, 
where $y_0=(\widetilde{y}, 1/(16))$. By the Harnack inequality, we have $f(x_0)\asymp f(y_0)$. Therefore it follows from \eqref{e:TAMSe6.45} that for all 
$x, y\in  D(2^{-8}, 2^{-8})$,
$$
\frac{f(x)}{f(y)}\le c_7\frac{x^p_d}{y^p_d},
$$
which is same as the conclusion of the theorem.

\noindent
(b) {\it The case $p<\alpha$. }
Set $w=(\widetilde{0}, 2^{-7})$. Then
\begin{align}\label{e:POTAe7.27}
f(w)&\ge \E_w\left[f(Y_{\tau_{U}}); Y_{\tau_{U}}\notin D(1, 1)\right]\nonumber\\
&\ge \E_w\int^{\tau_{B(w, 2^{-10})}}_0\int_{\R^d_+\setminus D(1, 1)}J^{\R^d_+}(Y_t, y)f(y)dydt\nonumber\\
&\ge c_{10}\E_w\tau_{B(w, 2^{-10})}\int_{\R^d_+\setminus D(1, 1)}J^{\R^d_+}(w, y)f(y)dy=c_{11}\int_{\R^d_+\setminus D(1, 1)}J^{\R^d_+}(w, y)f(y)dy,
\end{align}
where in the last line we used Proposition \ref{p:exit-time-estimate} (a).

We now show that there is a constant $c_{12}>0$ such that for any $z\in U$ and $y\in \R^d_+\setminus D(1, 1)$ it holds that 
\begin{equation}\label{e:new-estimate-for-J}
J^{\R^d_+}(z,y)\le c_{12}J^{\R^d_+}(w,y).
\end{equation} 
To prove this, first note that $|z-y|\asymp |w-y|\asymp |y|$. We consider two cases. The first case is when $y_d\ge 1/2$. Then
\begin{eqnarray*}
J^{\R^d_+}(z,y)&\preceq &\frac{z_d^{\beta_1} y_d^{\beta_2}}{|y|^{d+\alpha+\beta_1+\beta_2}}\log\left(1+1/z_d\right)^{\beta_3}(1+|\log y_d|^{\beta_3})  \log\left(1+\frac{|y|}{y_d}\right)^{\beta_4},  \\
J^{\R^d_+}(w,y)&\asymp &\frac{ y_d^{\beta_2}}{|y|^{d+\alpha+\beta_1+\beta_2}}(1+|\log y_d|^{\beta_3}) \log\left(1+\frac{|y|}{y_d}\right)^{\beta_4}, 
\end{eqnarray*}
implying \eqref{e:new-estimate-for-J} since $t\mapsto t^{\beta_1}\log(1+1/t)^{\beta_3}$ is bounded on $(0,1/2]$. Here we use notation $\preceq$ to denote that the left-hand side is smaller than the right-hand side multiplied by a constant.
The second case is  $y_d\le 1/2$. 
In this case, 
$$
\log\left(1+\frac{|y|}{(w_d\vee y_d)\wedge |y|}\right) \asymp \log |y|,
$$
and hence 
\begin{equation}\label{e:new-equation-w}
J^{\R^d_+}(w,y)\asymp \frac{y_d^{\beta_1}}{|y|^{d+\alpha+\beta_1+\beta_2}}|\log y_d|^{\beta_3}  \left(\log |y|\right)^{\beta_4}.
\end{equation}
If $z_d\le y_d$, then
 $(z_d\vee y_d)\wedge |y|=y_d$, so that
$$
\log\left(1+\frac{|y|}{(z_d\vee y_d)\wedge |y|}\right)^{\beta_4}\asymp \left(\log |y|+|\log y_d|\right)^{\beta_4}\preceq (\log |y|)^{\beta_4}+|\log y_d|^{\beta_4}. 
$$
Therefore
\begin{eqnarray*}
J^{\R^d_+}(z,y) &\preceq  & \frac{z_d^{\beta_1}y_d^{\beta_2}}{|y|^{d+\alpha+\beta_1+\beta_2}}\log\left(1+\frac{y_d}{z_d}\right)^{\beta_3}\left( (\log |y|)^{\beta_4}+|\log y_d|^{\beta_4}\right)\\
&\le &  \frac{z_d^{\beta_1}}{|y|^{d+\alpha+\beta_1+\beta_2}}|\log z_d|^{\beta_3}\left(\log |y|\right)^{\beta_4} + \frac{z_d^{\beta_1}}{|y|^{d+\alpha+\beta_1+\beta_2}}|\log z_d|^{\beta_3}y_d^{\beta_2} |\log y_d|^{\beta_4}\\
& \preceq & \frac{y_d^{\beta_1}}{|y|^{d+\alpha+\beta_1+\beta_2}}|\log y_d|^{\beta_3}\left(\log |y|\right)^{\beta_4},
\end{eqnarray*}
where the last approximate inequality follows from the fact that $t\mapsto t^{\beta_1}|\log t|^{\beta_3}$ is almost increasing  $t\mapsto t^{\beta_2}|\log t|^{\beta_4}$  bounded on $(0,1/2]$ and $y_d\le 1/2$. 
By comparing with \eqref{e:new-equation-w}, we see that \eqref{e:new-estimate-for-J} holds true. 

If $z_d\ge y_d$, then 
 $(z_d\vee y_d)\wedge |y|=z_d$, so that
$$
\log\left(1+\frac{|y|}{(z_d\vee y_d)\wedge |y|}\right)^{\beta_4}\asymp \left(\log |y|+|\log z_d|\right)^{\beta_4}\preceq (\log |y|)^{\beta_4}+|\log z_d|^{\beta_4}. 
$$
Therefore
\begin{eqnarray*}
J^{\R^d_+}(z,y)& \preceq & \frac{y_d^{\beta_1}z_d^{\beta_2}}{|y|^{d+\alpha+\beta_1+\beta_2}}\log\left(1+\frac{y_d}{z_d}\right)^{\beta_3}
\left( (\log |y|)^{\beta_4}+|\log z_d|^{\beta_4}\right)\\
&\le &  \frac{y_d^{\beta_1}}{|y|^{d+\alpha+\beta_1+\beta_2}}|\log y_d|^{\beta_3}\left(\log |y|\right)^{\beta_4} + 
 \frac{y_d^{\beta_1}}
{|y|^{d+\alpha+\beta_1+\beta_2}}|\log y_d|^{\beta_3}z_d^{\beta_2} |\log z_d|^{\beta_4}\\
&\preceq & \frac{y_d^{\beta_1}}{|y|^{d+\alpha+\beta_1+\beta_2}}|\log y_d|^{\beta_3}(\log |y|)^{\beta_4}.
\end{eqnarray*}
Again,  by comparing with \eqref{e:new-equation-w}, we get \eqref{e:new-estimate-for-J}.

Combining \eqref{e:new-estimate-for-J} with \eqref{e:exit-time-estimate-U} and \eqref{e:POTAe7.27} we now have
\begin{align}\label{e:POTAe7.29}
&\E_x\left[f(Y_{\tau_{U}}); Y_{\tau_{U}}\notin D(1, 1)\right]=\E_x\int^{\tau_{U}}_0\int_{\R^d_+\setminus D(1, 1)}J^{\R^d_+}(Y_t, y)f(y)dydt\nonumber\\
&\le c_{13} \E_x\tau_{U}\int_{\R^d_+\setminus D(1, 1)}J^{\R^d_+}(w, y)f(y)dy\le c_{14} x^p_df(w).
\end{align}

On the other hand, by the Harnack inequality, Carleson's estimate and Lemma \ref{l:POTAl7.4}, we have
\begin{align}\label{e:POTAe7.30}
\E_x\left[f(Y_{\tau_{U}}); Y_{\tau_{U}}\in D(1, 1)\right]&\le c_{15}f(x_0)\P_x\left(Y_{\tau_{U}}\in D(1, 1)\right)\le c_{16}f(x_0)x^p_d.
\end{align}
Combining \eqref{e:POTAe7.29} and \eqref{e:POTAe7.30}, and using the Harnack inequality and Carleson's estimate again, we get
\begin{align*}
f(x)&=\E_x\left[f(Y_{\tau_{U}}); Y_{\tau_{U}}\in D(1, 1)\right]+
\E_x\left[f(Y_{\tau_{U}}); Y_{\tau_{U}}\notin D(1, 1)\right]\nonumber\\
&\le c_{17}x^p_d(f(x_0)+f(w))\le c_{18}x^p_d f(x_0).
\end{align*}
Together with \eqref{e:TAMSe6.37} we get that
\begin{align}\label{e:POTAe7.32}
f(x)\asymp x^p_df(x_0).
\end{align}
Using \eqref{e:POTAe7.32} instead of \eqref{e:TAMSe6.45}, we also conclude that the theorem holds by the same argument.
\qed

\section{Proof of Theorem \ref{t:counterexample}}\label{s:counter}
In this section we assume that the assumptions of  Theorem \ref{t:counterexample} hold.
In particular, the parameters satisfy
 $\alpha+\beta_2 <p<\alpha+\beta_1$.

Suppose that the non-scale invariant BHP holds near the boundary of $\R^d_+$
(see the paragraph before Theorem \ref{t:counterexample}).

Note that by taking 
$g(x)= \P_x(Y_{\tau_{U}}\in D(1/2,1)\setminus D(1/2,3/4))$, 
we see from 
 Lemma \ref{l:exit-probability-estimate} that 
there exists $\wh{R} \in (0, 1)$ such that for any $r \in (0, \wh{R}\, ]$ 
there exists a constant $C_{39}=C_{39}(r)>0$ 
such that for any non-negative function $f$ in $\R^d_+$ which is harmonic  in $\R^d_+ \cap B(0, r)$ with respect to $Y$ and vanishes continuously on $ \partial \R^d_+ \cap B(0, r)$, 
\begin{equation}\label{e:bhp_mfx}
\frac{f(x)}{f(y)}\,\le \, C_{39} \,\frac{x_d^p}{y_d^p}, \quad  \hbox{for all } x, y\in  
\R^d_+ \cap B(0, r/2).
\end{equation}

Let $r_0=\wh{R}/4$  and choose a point $z_0 \in \partial \R^d_+$  with $|z_0| = 4$.
For  $n \in \N$, $B(z_0,1/n)$ does not intersect $B(0, 2r_0)$. We define
$$
K_n:=\int_{ \R^d_+   \cap B(z_0, 1/n)} |\log(y_d)|^{\beta_3  +\beta_4 }dy\,  ,\quad  \quad f_n(y):=
K_n^{-1}y_d^{-\beta_1} {\bf 1}_{\R^d_+ \cap B(z_0, 1/n)}(y),
$$
and 
$$
g_n(x):=\E_x\left[f_n(Y_{\tau_{U(r_0)}})\right]=\E_x\int_{0}^{\tau_{U(r_0)}}  \int_{\R^d_+ \cap B(z_0, 1/n)}  
J^{\R^d_+}(Y_t,y)f_n(y) dydt, \quad x \in {U(r_0)}.
$$
\begin{lemma}\label{l:gnlow} 
There exist $C_{40}>0$ and  $a_1 \in (0,1)$  such that 
\begin{align}
\label{e:gnlow}
\liminf_{n \to \infty}g_n(x)  \ge C_{40} x_d^{\beta_2+\alpha}
\end{align}
for all  $x=x^{(s)}=(\wt 0, s)\in \R^d_+$
with $s \in (0,  r_0a_1)$.
\end{lemma}
\pf 
Since 
$$
 6 >   |z-y|  >2 >y_d \wedge z_d \quad 
\text{ for } (y, z)  \in (\R^d_+ \cap B(z_0, 1/n)) \times   {U(r_0)},
$$
using 
\textbf{(B7)} 
we have for $(y, z)  \in (\R^d_+  \cap B(z_0, 1/n)) \times   {U(r_0)}$,
\begin{align*}
J^{\R^d_+}(z,y)&\asymp
(z_d \wedge y_d)^{\beta_1}  (z_d \vee y_d)^{\beta_2}  \log\left(1+\frac{z_d \vee y_d}{z_d \wedge y_d}\right)^{\beta_3}  \log\left( \frac{1}{z_d \vee y_d}\right)^{\beta_4}\\
&\asymp
  z_d^{\beta_1} y_d^{\beta_1}  (z_d \vee y_d)^{-(\beta_1-\beta_2)}  \log\left(1+\frac{z_d \vee y_d}{z_d \wedge y_d}\right)^{\beta_3}  \log\left( \frac{1}{z_d \vee y_d}\right)^{\beta_4} .
\end{align*}
Therefore, for $ x \in {U(r_0)}$,
\begin{align}
\label{e:jnew2}
&g_n(x)\asymp \nonumber\\
& 
K_n^{-1}
\E_x\int_{0}^{\tau_{U(r_0)}}  (Y_t^d)^{\beta_1} \int_{\R^d_+ \cap B(z_0, 1/n)} 
 (Y_t^d \vee y_d)^{-(\beta_1-\beta_2)}   \log\left(1+\frac{Y_t^d \vee y_d}{Y_t^d \wedge y_d}\right)^{\beta_3}  \log\left(\frac{1}{Y_t^d \vee y_d}\right)^{\beta_4}    dydt.
\end{align}
Note that, using 
$\sup_{ t\ge 1}t^{-(\beta_1-\beta_2)} \log(1+t)^{\beta_3}<\infty$, 
 for $z\in U(r_0)$,
\begin{align}
\label{e:jnew12}
&  K_n^{-1}
 \int_{\R^d_+ \cap B(z_0, 1/n)}(z_d \vee y_d)^{-(\beta_1-\beta_2)}   \log\left(1+\frac{z_d \vee y_d}{z_d \wedge y_d}\right)^{\beta_3} 
  \log\left(\frac{1}{z_d \vee y_d}\right)^{\beta_4}
 dy  \nn\\
\le & K_n^{-1}
(z_d)^{-(\beta_1-\beta_2)}  \int_{\R^d_+ \cap B(z_0, 1/n)\cap \{z_d \le y_d\}} (y_d/z_d)^{-(\beta_1-\beta_2)}\log\left(1+\frac{y_d}{z_d}\right)^{\beta_3}  
 \log\left(\frac{1}{ y_d}\right)^{\beta_4}  dy\nn\\
&+ K_n^{-1}(z_d)^{-(\beta_1-\beta_2)}
 \int_{\R^d_+ \cap B(z_0, 1/n)\cap \{z_d > y_d\}  
}\log\left(1+\frac{z_d}{y_d}\right)^{\beta_3}   \log\left(\frac{1}{z_d }\right)^{\beta_4} dy \nn\\
 \le& c K_n^{-1}(z_d)^{-(\beta_1-\beta_2)}
 \int_{\R^d_+ \cap B(z_0, 1/n)}   \log\left(\frac{1}{y_d}\right) ^{\beta_3  +\beta_4 }
 dy   \le c z_d^{-(\beta_1-\beta_2)}   
\end{align} 
and
$$
\lim_{n \to \infty}
 K_n^{-1}
 \int_{\R^d_+  \cap B(z_0, 1/n)} ( z_d \vee y_d)^{-(\beta_1-\beta_2)}   \log\left(1+\frac{z_d \vee y_d}{z_d \wedge y_d}\right)^{\beta_3}    \log\left(\frac{1}{z_d \vee y_d}\right)^{\beta_4}  dy =
  z_d^{-(\beta_1-\beta_2)}\, .
 $$
Moreover, by Lemma \ref{l:exit-time estimate-U_new},
$\E_x\int_{0}^{\tau_{U(r_0)}}  (Y_t^d)^{\beta_2}dt  < \infty$ for all 
$x\in U(r_0)$. 
 Thus we can use  the dominated convergence theorem to get that  for all  
 $x\in U(r_0)$,
 \begin{align}\label{e:notD231}
&\lim_{n \to \infty} K_n^{-1}
\E_x\int_{0}^{\tau_{U(r_0)}}  (Y_t^d)^{\beta_1} \int_{\R^d_+ \cap B(z_0, 1/n)} 
 (Y_t^d \vee z_d)^{-(\beta_1-\beta_2)}   \log\left(1+\frac{Y_t^d \vee y_d}{Y_t^d \wedge y_d}\right)^{\beta_3}   \log\left(\frac{1}{Y_t^d \vee y_d}\right)^{\beta_4}  dydt \nn \\
&=\E_x\int_{0}^{\tau_{U(r_0)}}  (Y_t^d)^{\beta_1} (Y_t^d)^{-(\beta_1-\beta_2)}dt
=  \E_x\int_{0}^{\tau_{U(r_0)}}  (Y_t^d)^{\beta_2}dt. 
\end{align} 
Combining
\eqref{e:notD231} with Lemma \ref{l:lower-bound-for-integral-q}
 we conclude that \eqref{e:gnlow} holds true.
 \qed
 
\noindent
\textbf{Proof of Theorem \ref{t:counterexample}.} 
From \eqref{e:jnew2}, \eqref{e:jnew12} and Lemma \ref{l:exit-time estimate-U_new}, we see 
that for all $x\in U(r_0)$,
\begin{align}
\label{e:jnew3}
&g_n(x)\nonumber\\
 &\le  c
 K_n^{-1}
\E_x\int_{0}^{\tau_{U(r_0)}}  (Y_t^d)^{\beta_1} \int_{\R^d_+ \cap B(z_0, 1/n)} 
 (Y_t^d)^{-(\beta_1-\beta_2)}  \log\left(1+\frac{Y_t^d \vee y_d}{Y_t^d \wedge y_d}\right) ^{\beta_3}   \log\left(\frac{1}{Y_t^d \vee y_d}\right)^{\beta_4}  dydt \nn\\
& 
\le c \E_x\int_{0}^{\tau_{U(r_0)}}  (Y_t^d)^{\beta_2}  dt \le c x_d^{\beta_2+ \alpha}. \end{align}
Thus the $g_n$'s are non-negative functions  in $\R^d_+ $ 
which are harmonic  in $\R^d_+  \cap B(0,2^{-2}  r_0)$ with respect to $Y^{\R^d_+, \kappa} $ and vanish continuously on $ \partial \R^d_+  \cap B(0,2^{-2}  r_0)$. Therefore, by \eqref{e:bhp_mfx},
$$
\frac{g_n(y)}{g_n(w)} \le 
C_{39} 
\frac{(y_d)^p }{(w_d)^p } \quad 
\text{for all } y\in D \cap B(0,2^{-3}  r_0),
$$
where $w=(\wt{0}, 2^{-3} r_0)$ and  
$C_{39}=C_{39}(2^{-2}  r_0)$. 
Thus by \eqref{e:jnew3}, for all 
$y\in \R^d_+  \cap B(0,2^{-3}  r_0)$,
$$
\limsup_{n \to \infty}g_n(y) \le 
 C_{39} 
\limsup_{n \to \infty}g_n(w) \frac{(y_d)^p }{(w_d)^p } \le c_4 (y_d)^p.
$$
This and \eqref{e:gnlow} imply that  for all  
$x=x^{(s)}=(\wt 0, s)\in \R^d_+$
with $s \in (0, 2^{-3} a_1 r_0)$, $x_d^{\beta_2+\alpha} \le c x_d^{p} $, which  gives a contradiction because $\beta_2+\alpha<p$ . \qed

\section{Proof of Lemma \ref{l:key-lemma}}\label{s:proof-of-key-lemma}\label{s:8}
In this section we prove Lemma \ref{l:key-lemma}. 
\emph{Throughout this section, we assume that \textbf{(B1)}, \textbf{(B4)}, \textbf{(B7)}, \textbf{(B8)} and  \eqref{e:kappap} hold. }
Recall that, because of \textbf{(B4)}, we can, without loss of generality, assume that $B(x, x)=1$.
We start with three auxiliary lemmas estimating certain  integrals, 
and then state and 
prove a technical lemma needed later. 
We then define the function $\psi$, show that (b1), (b2), (c1) and (c2) are immediate consequences of the definition of $\psi$, and that (b3) and (c3) follow directly from 
Lemma \ref{l:estimate-of-L-hat-B}
and Lemma \ref{l:key-lemma} (a). The remaining part, namely the proof of (a), is divided into two cases depending on the location of the point $z\in U$.

\subsection{Auxiliary lemmas}
In the next lemma and throughout the section we use the notation $\R^d_- :=\{x=(\wt{x}, x_d)\in \R^d:\,  x_d<0\}$.
\begin{lemma}\label{l:calc}
For $\alpha\in [1, 2)$,
$$
\int_{\R^d_- \cap \{ |y-z|<6\}}\frac{1}{|y-z|^{d+\alpha-1}}dy \le C_{41} \begin{cases}z_d^{1-\alpha}, &\alpha\in (1, 2);\\ 
|\log z_d|, &\alpha=1,
\end{cases}
\quad \quad z\in U.
$$
\end{lemma}
\pf 
Immediate by enlarging the area of integration to $B(z, 6) \setminus B(z, z_d)$.
\qed

\begin{lemma}\label{l:1}
For  $k \in \bN$, let
$I(k):=\int_{D(7, 7)}{|y-z|^{-d-\alpha+k}}\sB (y, z)\, dy$. There exists $C_{42}>0$ such that
\begin{itemize}
\item[(a)] if  $\alpha<k\le \alpha+\beta_1$, then
$I(k)\le C_{42} z_d^{k-\alpha} \,  |\log z_d|^{\beta_3+1}$ for  $z\in D(2^{-1}, 2^{-1})$;
\item[(b)]  if $k$ is the smallest integer such that $k>\alpha+\beta_1$, then
$I(k)\le  C_{42} z_d^{\beta_1}\, |\log z_d|^{\beta_3}$ for  $ z\in D(2^{-1}, 2^{-1})$.
\end{itemize}
\end{lemma}

\pf 
We write
\begin{align*}
\int_{D(7, 7)\cap\{|y-z|<z_d\}}\frac{\sB (y, z)dy}{|y-z|^{d+\alpha-k}}+\int_{D(7, 7)\cap\{|y-z|\ge z_d\}}\frac{\sB (y, z)dy}{|y-z|^{d+\alpha-k}}=:I+II.
\end{align*}
Then,  by using that $\sB (y,z)\le C_1$, 
$$
I\le C_1 \int_{D(7, 7)\cap\{|y-z|<z_d\}}|y-z|^{-d-\alpha+k}dy\le c_1\int^{z_d}_0r^{-1-\alpha+k}dr
\le c_2 z_d^{k-\alpha} \, .
$$

For $II$ we use Lemma  \ref{l:key-estimate-for-B}(b) together with the fact 
that $|y|\le 7\sqrt{2}$ for $y\in D(7,7)$, to get the estimate
$
\sB (y,z)\le c_3 z_d^{\beta_1}|\log z_d|^{\beta_3}|y-z|^{-\beta_1} .
$
This  implies that
\begin{align*}
&II  \le 
\begin{cases}
c_3 z_d^{\beta_1}|\log z_d|^{\beta_3}\int^{c_4}_{z_d}r^{-1-\alpha-\beta_1+k}dr\le c_5 z_{d}^{k-\alpha}|\log z_d|^{\beta_3+1}& \text{if } \alpha<k<\alpha+\beta_1;\\
c_3 z_d^{\beta_1}|\log z_d|^{\beta_3}\int^{c_4}_{z_d}r^{-1}dr
\le c_5 z_d^{k-\alpha}|\log z_d|^{\beta_3+1}&\text{if } \alpha<k=\alpha+\beta_1;\\
c_3 z_d^{\beta_1}|\log z_d|^{\beta_3}\int^{c_4}_{z_d}r^{-1-\alpha-\beta_1+k}dr
\le c_5z_d^{\beta_1} |\log z_d|^{\beta_3}&\text{if } k>\alpha+\beta_1.
\end{cases}
\end{align*}
\qed

\begin{lemma}\label{l:3}
For $\alpha\in [1, 2)$,
$$
\int_{D(7, 7)}\frac{1}{|y-z|^{d+\alpha-1}}|1-\sB (y, z)|dy
\le  C_{43} \begin{cases}z_d^{1-\alpha}, & \alpha \in (1, 2);\\
|\log z_d|, &\alpha=1,
\end{cases} \quad z\in D(2^{-1}, 2^{-1}).
$$
\end{lemma}
\pf 
We write 
$$
\int_{D(7, 7)}\frac{|1-\sB (y, z)|}{|y-z|^{d+\alpha-1}}dy=
\int_{D(7, 7)\cap \{|y-z|<z_d/2\}}+\int_{D(7, 7)\cap \{|y-z|\ge z_d/2\}}=:I+II.
$$
If $y\in B(z, 2^{-1}z_d)$, then $|y-z|\le z_d/2\le y_d$ and $y_d\asymp z_d$ hence by \textbf{(B4)}, 
we have that 
\begin{align*}
I&\le c_1z_d^{-\theta}\int_{|y-z|<z_d/2}|y-z|^{\theta-d-\alpha+1}dy
=c_2z_d^{-\theta}
\int_0^{z_d/2}r^{\theta-\alpha}dr=c_3z_d^{1-\alpha} \, .
\end{align*}
Clearly, $|1-\sB (y,z)|\le c_4$. Therefore,
$$
II\le \int_{D(7, 7)\cap \{|y-z|\ge z_d/2\}}|y-z|^{-d-\alpha+1}dy
=c_5\int_{z_d/2}^Kr^{-\alpha}\le
c_6\begin{cases}z_d^{1-\alpha}, & \alpha \in (1, 2);\\
\log\frac1{z_d}, &\alpha=1.
\end{cases}
$$
\qed

\begin{lemma}\label{l:non-negativity}
Let $N\ge 2$ be an integer and $C>0$. For $t,s>0$, let
$$
F(t, s; C, N):=\left(\frac{t}{s}\right)^N-C\sum_{k=1}^{N-1}\left(\frac{t}{s}\right)^{N-k}|\log s|^{1+\beta_3}.
$$
There exists $\wt{\kappa}=\wt{\kappa}(C,N, \beta_3)\in (0, e^{-1})$ 
such that if $s\le \wt{\kappa}$ and $t\ge s^{1/2}$, 
then $F(t,s; C, N)\ge 0$.
\end{lemma}
\pf 
By defining $c_1=C(N-1)$, we can rewrite $F(t,s; C,N)$ in the form
\begin{align*}
F(t,s;C,N)&=\frac{1}{N-1}\sum_{k=1}^{N-1}\left(\frac{t}{s}\right)^{N-k}\left( \Big(\frac{t}{s}\Big)^k-c_1|\log s|^{1+\beta_3}\right)\\
&\ge \frac{1}{N-1}\sum_{k=1}^{N-1}\left(\frac{t}{s}\right)^{N-k}\left(\frac{t}{s}-c_1|\log s|^{1+\beta_3}\right)
\end{align*}
provided that $t\ge s$. 
Since $\lim_{s\to 0}s^{1/2}|\log s|^{\beta_3+1}=0$, there exists $\wt{\kappa}=\wt{\kappa}(c_1, \beta_3)\in (0, e^{-1})$ such that for $s\le \wt{\kappa}$ we have  
$c_1 s^{1/2}|\log s|^{\beta_3+1}\le 1$. If now $t\ge s^{1/2}$ , then $t/s^{1/2}\ge 1\ge c_1
s^{1/2}|\log s|^{\beta_3+1}$, implying $t/s-c_1 |\log s|^{\beta_3+1}\ge 0$.
This implies that $F(t,s;C,N)\ge 0$ for $s\le \wt{\kappa}$ and $t\ge s^{1/2}$.
\qed

\subsection{Definition of $\psi$}
Let $N:=\lfloor\alpha+\beta_1\rfloor+1$ (the smallest integer larger than $\alpha+\beta_1$) and assume $\gamma>2N+2$. Let $\psi$ be a non-negative 
$C^N$ function in the upper half space with bounded support
and uniformly bounded derivatives up to the order $N$
 such that
$$
\psi(y)=\begin{cases}|\widetilde{y}|^\gamma, &y\in D(2^{-2}, 2^{-2});\\
                                 1,&y\in D(2, 2)\setminus U;\\
                                 0, & y\in D(3, 3)^c,
\end{cases}
$$
and $\psi(y)\ge 4^{-\gamma}$ for $y\in U\setminus D(2^{-2}, 2^{-2})$. 
The function $\psi$ in the half space can be constructed so that, for $y=(\widetilde{y}, y_d)$ with
$y_d\in (0, \frac18)$, $\psi(y)$ depends on $\widetilde{y}$ only. 
We  extend  $\psi$ to be identically zero in the lower half space.

We recall the definition of the  operator $L^\sB $ from \eqref{e:defn-LB}. 
For $\alpha\in (0, 1)$, $L^\sB_ \alpha\psi(z)$ is not really a principal value integral.
For $\alpha\in [1, 2)$, 
it follows from Proposition \ref{p:operator-interpretation} (a)
that, for $z\in \R^d_+$,  the principal value integral $L^\sB_ \alpha\psi(z)$
can be interpreted as the expression after the last equality in
the display below
\begin{align*}
&\int_{\R^d}\frac{\psi(y)-\psi(z)-\nabla \psi(z) {\bf 1}_{\{ |y-z|<6\}}\cdot (y-z)}{|y-z|^{d+\alpha}}dy
+\int_{\R^d_-}\frac{\psi(z)}{|y-z|^{d+\alpha}}dy\nn\\
&\qquad +\int_{\R^d_+}\frac{\psi(y)-\psi(z)}{|y-z|^{d+\alpha}}(\sB (y, z)-1)dy\nn\\
&=\int_{\R^d_+}\frac{\psi(y)-\psi(z)-\nabla \psi(z) {\bf 1}_{\{ |y-z|<6\}}\cdot (y-z)}{|y-z|^{d+\alpha}}dy
+\int_{\R^d_+}\frac{\psi(y)-\psi(z)}{|y-z|^{d+\alpha}}(\sB (y, z)-1)dy\nn\\
&\qquad-\int_{\R^d_-}\frac{\nabla \psi(z) {\bf 1}_{\{ |y-z|<6\}}\cdot(y-z)}{|y-z|^{d+\alpha}}dy\nn\\
&=\int_{\R^d_+}\frac{\psi(y)-\psi(z)-\nabla \psi(z) {\bf 1}_{\{ |y-z|<6\}}\cdot (y-z)}{|y-z|^{d+\alpha}}\sB (y, z)dy
\nn\\
&\qquad+\int_{\R^d_+}\frac{\nabla \psi(z) {\bf 1}_{\{ |y-z|<6\}}\cdot (y-z)}{|y-z|^{d+\alpha}}(\sB (y, z)-1)dy-\int_{\R^d_-}\frac{\nabla \psi(z) {\bf 1}_{\{ |y-z|<6\}}\cdot(y-z)}{|y-z|^{d+\alpha}}dy.
\end{align*}
Thus, since $|y| \ge |y-z| -|z| >5$  for $z \in U$ and $|y-z|>6$, we have that for $z \in U$
\begin{align}\label{e:b}
&L^\sB_ \alpha\psi(z)=\int_{\R^d_+ \cap \{ |y-z|<6\}}\frac{\psi(y)-\psi(z)-\nabla \psi(z)\cdot (y-z)}{|y-z|^{d+\alpha}}\sB (y, z)dy\nn\\
&\qquad \qquad-\psi(z)\int_{\R^d_+ \cap \{ |y-z|>6\}}\frac{1}{|y-z|^{d+\alpha}}\sB (y, z)dy
\nn\\
&\qquad\qquad+\int_{\R^d_+ \cap \{ |y-z|<6\}}\frac{\nabla \psi(z)\cdot (y-z)}{|y-z|^{d+\alpha}}(\sB (y, z)-1)dy-\int_{\R^d_- \cap \{ |y-z|<6\}}\frac{\nabla \psi(z) \cdot(y-z)}{|y-z|^{d+\alpha}}dy\nn\\
& \le \int_{\R^d_+ \cap \{ |y-z|<6\}}\frac{\psi(y)-\psi(z)-\nabla \psi(z)\cdot (y-z)}{|y-z|^{d+\alpha}}\sB (y, z)dy
\nn\\
&\qquad+\int_{\R^d_+ \cap \{ |y-z|<6\}}\frac{\nabla \psi(z) \cdot (y-z)}{|y-z|^{d+\alpha}}(\sB (y, z)-1)dy-\int_{\R^d_- \cap \{ |y-z|<6\}}\frac{\nabla \psi(z)\cdot(y-z)}{|y-z|^{d+\alpha}}dy.
\end{align}

\subsection{Proof of Lemma \ref{l:key-lemma} (b) and (c)}
Let $p\in ((\alpha-1)_+, \alpha+\beta_1)$, and recall that $h_p(x)=x^p_d{\bf 1}_{D(1, 1)}(x)$.
Define $\phi:=h_p-\psi$.
The function $\phi$ is obviously non-positive on $U^c$, hence 
Lemma \ref{l:key-lemma} (b2) holds true. Moreover, since $\psi((\wt{0},x_d))=0$, we have that $\phi((\wt{0}, x_d))=x_d^p$, for $(\wt{0}, x_d)\in U$, which is  Lemma \ref{l:key-lemma} (b1). Similarly, for $(\alpha-1)_+<q<p$,the  function  $\varphi:=h_q-\psi$ satisfies Lemma \ref{l:key-lemma} (c1) and (c2).  
Furthermore  Lemma \ref{l:key-lemma} (b3) and (c3) follow from  Lemma \ref{l:estimate-of-L-hat-B}  and Lemma \ref{l:key-lemma} (a).
In fact, 
\begin{align*}
L^\sB \phi(z)=L^\sB h_p(z)-L^\sB \psi(z)\ge (-C_{22}-C_{29})z_d^{\beta_1}\,  |\log z_d|^{\beta_3} , \qquad 
z\in U(r_0).
\end{align*}
 By Lemma \ref{l:estimate-of-L-hat-B} (c) and part (a) of this lemma, 
for all $x\in U$ we have
\begin{eqnarray*}
L^\sB  \varphi(x)&=&L^\sB  h_q(x)-L^\sB  \psi(x) \ge -C_{26} 
x_d^{q-\alpha}-(C_{27}+C_{29}) x_d^{\beta_1}|\log x_d|^{\beta_3}.
\end{eqnarray*}
\qed

\subsection{Proof of Lemma \ref{l:key-lemma} (a)}
It remains to prove Lemma \ref{l:key-lemma} (a). 
The proof is split into two cases: (i) $z\in D(2^{-2}, 2^{-2})$,  and (ii) $z\in U\setminus D(2^{-2}, 2^{-2})$. 
\subsubsection{Proof when $z\in  D(2^{-2}, 2^{-2})$}
Throughout this subsection we assume that $z\in D(2^{-2}, 2^{-2})$ so that $\psi(z)=\psi(\wt{z})=|\wt{z}|^{\gamma}$. We will use Taylor's theorem for the function $\psi(z)$ so we first briefly go over the notation. 
For a multi-index $\rho=(\rho_1, \dots, \rho_{d-1})$, we define $|\rho|:=\rho_1+\cdots +\rho_{d-1}$ and 
$$
\partial^{\rho}\psi(\wt{z}) 
=\partial_1^{\rho_1}\dots \partial_{d-1}^{\rho_{d-1}} \psi(\wt{z})=\frac{\partial^ {|\rho|}\psi(\wt{z})}{\partial z_1^{\rho_1}\cdots \partial z_{d-1}^{\rho_{d-1}}}\, .
$$
Let further $\rho!=\rho_1!\dots \rho_{d-1}!$ and $\wt{z}^{\rho}=z_1^{\rho_1}\dots z_{d-1}^{\rho_{d-1}}$. 
Then (recall that $N=\lfloor \alpha+\beta_1\rfloor +1$)
\begin{equation}\label{e:taylor}
\psi(y)
=\sum_{|\rho|\le N-1}\frac{\partial^{\rho}\psi(\wt{z})}{\rho!}(\wt{y}-\wt{z})^{\rho}+
\sum_{|\hat \rho|=N}\frac{\partial^{\hat \rho}\psi(c{y}+(1-c){z})}{\hat \rho !}({y}-{z})^{\hat \rho}
\end{equation}
for some $c\in (0,1)$ where $\hat \rho=(\hat \rho_1, \dots,  \hat \rho_{d})$. 
It is easy to see that, for any non-negative integer $k\le N$,
\begin{equation}\label{e:sum-of-partials}
\Big|\sum_{|\rho|=k}\frac{\partial^{\rho} \psi(\wt{z})}{\rho!}\Big| \le c(\gamma, k)|\wt{z}|^{\gamma-k}
\end{equation}
for some $c(\gamma, k)>0$. Similarly, the terms in Taylor's formula \eqref{e:taylor} are estimated as
\begin{align}
&\Big|\sum_{|\rho|= k}\frac{\partial^{\rho}\psi(\wt{z})}{\rho!}(\wt{y}-\wt{z})^{\rho}\Big| 
\le C(k) |\wt{z}|^{\gamma-k} |\wt{y}-\wt{z}|^k\le  C |\wt{z}|^{\gamma-k}|y-z|^k\, ,\quad k\le N-1, \label{e:taylor-term}\\
&
\Big| \sum_{|\hat \rho|=N}\frac{\partial^{\hat \rho}\psi(c{y}+(1-c){z})}{\hat \rho !}({y}-{z})^{\hat \rho}\Big|\le C(N)|y-z|^N\, ,
 \label{e:taylor-remainder}
\end{align}
for some constants $C(k)=C(\gamma, k)>0$, $k=0, 1,\dots, N$.

Recall that $N=\lfloor \alpha+\beta_1\rfloor+1$ and $\gamma>2N+2$. Suppose first that $\alpha\in (0,1)$ and expand  $L^\sB_\alpha\psi(z)$ further:
\begin{align}
L^\sB_\alpha\psi(z)=
&\int_{\R^d_+ \cap \{ |y-z|<6\}} \frac{\psi(y)-\psi(z)}{|y-z|^{d+\alpha}}\sB(y, z)dy- \psi(z)\int_{\R^d_+ \cap \{ |y-z|>6\}} \frac{1}{|y-z|^{d+\alpha}}\sB(y, z)dy\nn\\
 \le
&\int_{\R^d_+ \cap \{ |y-z|<6\}} \frac{\psi(y)-\psi(z)}{|y-z|^{d+\alpha}}\sB(y, z)dy \label{e:case1-numerator}\\
=&\int_{\R^d_+ \cap \{ |y-z|<6\}}\frac{\psi(y)-\psi(z)-\sum_{1\le |\rho|\le N-1}\frac{\partial^{\rho}\psi(\wt{z})}{\rho!} (\wt{y}-\wt{z})^{\rho}}{|y-z|^{d+\alpha}}\sB(y, z)dy \nn
\\
&+\int_{\R_+^d \cap 
\{|y-z|<6\} 
}\sum_{2\le |\rho|\le N-1}\frac{\partial^{\rho}\psi(\wt{z})}{\rho!} \frac{(\wt{y}-\wt{z})^{\rho}}{|y-z|^{d+\alpha}}\sB(y, z)dy\nn
\\
&+\int_{\R_+^d \cap 
\{|y-z|<6\} }\frac{\partial\psi(\wt{z})}{1!} \frac{(\wt{y}-\wt{z})}{|y-z|^{d+\alpha}}\sB(y, z)dy =: I + II +III. \nn
\end{align}
When $\alpha\in [1,2)$, we use \eqref{e:b} to expand $L^{\sB}_{\alpha}$ further to get
\begin{align*}
L^\sB_\alpha\psi(z)\le& \int_{\R^d_+ \cap \{ |y-z|<6\}}\frac{\psi(y)-\psi(z)-\sum_{1\le |\rho|\le N-1}\frac{\partial^{\rho}\psi(\wt{z})}{\rho!} (\wt{y}-\wt{z})^{\rho}}{|y-z|^{d+\alpha}}\sB(y, z)dy 
\\
&+\int_{\R_+^d \cap\{|y-z|<6\}}
\sum_{2\le |\rho|\le N-1}\frac{\partial^{\rho}\psi(\wt{z})}{\rho!} \frac{(\wt{y}-\wt{z})^{\rho}}{|y-z|^{d+\alpha}}
\sB(y, z)dy
\\
& + \int_{\R^d_+ \cap \{ |y-z|<6\}}\frac{\nabla \psi(z)\cdot (y-z)}{|y-z|^{d+\alpha}}(\sB(y,z)-1)\, dy - \int_{\R^d_- \cap \{ |y-z|<6\}}\frac{\nabla \psi(z)\cdot (y-z)}{|y-z|^{d+\alpha}}  dy 
\\
&=:I+II+IV.
\end{align*}

By using \eqref{e:taylor} we estimate the numerator in $I$ by \eqref{e:taylor-remainder}. The obtained integral is then estimated by use of Lemma \ref{l:1} (b) by a constant multiple of $z_d^{\beta_1}\, |\log z_d|^{\beta_3}$ to get that
$$
I\le C(N)z_d^{\beta_1}\, |\log z_d|^{\beta_3}\, .
$$ 
The terms in $II$ are estimated by \eqref{e:taylor-term}, and the obtained integrals by Lemma \ref{l:1} (a) to get that
$$
II\le \sum_{k=2}^{N-1}C(k)|\wt{z}|^{\gamma-k}z_d^{k-\alpha}\, |\log z_d|^{\beta_3+1}.
$$
By the same argument we get the estimate
$$
III\le C(1)|\wt{z}|^{\gamma-1}z_d^{1-\alpha}\, |\log z_d|^{\beta_3+1}.
$$
The terms in $IV$ are estimated by use Lemmas \ref{l:3} and \ref{l:calc} to get that
$$
IV \le c |\wt{z}|^{\gamma-1}z_d^{1-\alpha}|\log z_d|\le  c |\wt{z}|^{\gamma-1}z_d^{1-\alpha}|\log z_d|^{\beta_3+1}.
$$
Combining the estimates for $I$, $II$, $III$ and $IV$, we get that for $\alpha\in (0,2)$ it holds that
\begin{equation}\label{e:case1-estimate}
L^B_{\alpha}(z)\le \wh{C} \left(z_d^{\beta_1} \, |\log z_d|^{\beta_3}+ \sum^{N-1}_{k=1}|\wt{z}|^{\gamma-k}z_d^{k-\alpha}\,  |\log z_d|^{\beta_3+1} \right)
\end{equation}
for some constant $\wh{C}>0$.

 Hence
\begin{align}
&L^\sB\psi(z)\nonumber\\
&\le  \wh{C} z_d^{\beta_1} \, |\log z_d|^{\beta_3}+ \wh{C}\sum^{N-1}_{k=1}|\wt{z}|^{\gamma-k}z_d^{k-\alpha}\,  |\log z_d|^{\beta_3+1} -C(\alpha, p, \sB) |\wt{z}|^{\gamma}z_d^{-\alpha} \label{e:case1-1} \\
&= \wh{C} z_d^{\beta_1}\, |\log z_d|^{\beta_3}-C(\alpha, p, \sB)\frac{|\widetilde{z}|^{\gamma-N}}{z_d^{\alpha-N}}
\left(\left(\frac{|\widetilde{z}|}{z_d}\right)^N-\frac{ \wh{C}}{C(\alpha, p, \sB)}\sum^{N-1}_{k=1}\left(\frac{|\widetilde{z}|}{z_d}\right)^{N-k}\, |\log z_d|^{\beta_3+1}\right)\nonumber \\
&= \wh{C} z_d^{\beta_1}\, |\log z_d|^{\beta_3}-C(\alpha, p, \sB)\frac{|\widetilde{z}|^{\gamma-N}}{z_d^{\alpha-N}}
\left(\left(\frac{|\widetilde{z}|}{z_d}\right)^N-C\sum^{N-1}_{k=1}\left(\frac{|\widetilde{z}|}
{z_d}\right)^{N-k}\, |\log z_d|^{\beta_3+1}\right), \label{e:case1-2}
\end{align}
where $C:= \wh{C}/C(\alpha, p, \sB)$,

Now we take $z_d\le \wt{\kappa}(C,N)$ (where $\wt{\kappa}$ comes from Lemma \ref{l:non-negativity}) and $|\wt{z}|\ge z_d^{1/2}$. Then it follows from \eqref{e:case1-2} and Lemma \ref{l:non-negativity} that 
$$
L^\sB \psi(z)\le \wh{C} z_d^{\beta_1}\, |\log z_d|^{\beta_3}-C(\alpha, p, \sB)
\frac{|\widetilde{z}|^{\gamma-N}}{z_d^{\alpha-N}} F(|\wt{z}|, z_d; C,N)
\le  \wh{C} z_d^{\beta_1}|\log z_d|^{\beta_3}.
$$
In case when $|\wt{z}|\le z_d^{1/2}$ and $z_d\le \wt{\kappa}$, we estimate the middle term in \eqref{e:case1-1} as
\begin{eqnarray*}
\lefteqn{\sum^{N-1}_{k=1}|\wt{z}|^{\gamma-k}z_d^{k-\alpha}\,  |\log z_d|^{\beta_3+1}=\frac{|\wt{z}|^{\gamma-N}}{z_d^{\alpha-N}(z_d^{1/2})^{N}}\sum_{k=1}^{N-1}\left(\frac{|\wt{z}|}{(z_d)^{1/2}}\right)^{N-k}(z_d^{1/2})^{k} |\log z_d|^{\beta_3+1}}\\
&=& \Big(\frac{|\wt{z}|}{z_d^{1/2}}\Big)^{\gamma-N}z_d^{\frac{\gamma}{2}-\frac12 -\alpha}\, |\log z_d|^{\beta_3+1}\left(\sum_{k=1}^{N-1}\big(\frac{|\wt{z}|}{z_d^{1/2}}\Big)^{N-k} (z_d^{1/2})^{k+1}\right)\\
&\le &  z_d^{\frac{\gamma}{2}-\frac12 -\alpha}\, |\log z_d|^{\beta_3+1} (N-1)\le c z_d^{\beta_1}\, .
\end{eqnarray*}
In the last line we first used that $|\wt{z}|/z_d^{1/2}\le 1$ and $z_d\le 1$ and then that $\gamma/2-1/2-\alpha> \beta_1$. 
Thus, in this case we can disregard the last term in \eqref{e:case1-1} and obtain again that $L^\sB\psi(z)\le C_{44}z_d^{\beta_1}|\log z_d|^{\beta_3}$ with $C_{44}=c\wh{C}$.

Finally, it follows from \eqref{e:case1-estimate} that for $z\in U$ with $z_d\ge \kappa:=\wt{\kappa}$ it holds that $L^\sB\psi(z)\le C_{45}$ for some constant $C_{45}=C_{45}(\kappa)>0$.

Set $C_{46}:=C_{45} ((\kappa^{\beta_1}(\log 2)^{\beta_3}))\vee 1)^{-1}$. Then $L^\sB \psi(z)\le (C_{44}\vee C_{46})z_d^{\beta_1} |\log z_d|^{\beta_3}$ for all $z\in D(2^{-2}, 2^{-2})$.

\subsubsection{Proof when $z\in U\setminus D(2^{-2}, 2^{-2})$}
Throughout this subsection we assume that  $z\in U\setminus D(2^{-2}, 2^{-2})$. 
We show that there exist constants $C_{47}>0$ and $\kappa\in (0,1/4]$ such that (i) for $z_d\le \kappa$ and  $|\wt{z}|\in [1/4,1/2)$ it holds that $L^\sB\psi(z)\le 0$; (ii) For $z_d\in [\kappa, 1/2)$  it holds that $L^\sB\psi(z)\le C_{47}$.

When $\alpha\in (0,1)$ we see from \eqref{e:case1-numerator} that
$$
L^\sB_{\alpha}\psi(z)\le c \int_{D(7,7)}\frac{1}{|y-z|^{d+\alpha-1}}\sB(y,z)\, dy\, .
$$
Combining this with Lemma \ref{l:1} we get that
\begin{equation}\label{e:bigger14-1}
L^\sB_{\alpha}\psi(z)\le C_{48}\left(z_d^{\beta_1}+z_d^{1-\alpha} \right)|\log z_d|^{\beta_3+1}\, , \quad z\in U\, .
\end{equation}

When $\alpha\in [1,2)$ we use \eqref{e:b} to estimate
\begin{eqnarray*}
L^\sB_{\alpha}\psi(z)&\le &c\left(\int_{D(7,7)}\frac{1}{|y-z|^{d+\alpha-2}}\sB(y,z)\, dy+\int_{D(7,7)}\frac{1}{|y-z|^{d+\alpha-1}}|1-\sB(y,z)|\, dy\right. \nonumber \\
& & \left. +\int_{\R_-^d\cap\{|y-z|<6\}}\frac{1}{|y-z|^{d+\alpha-1}}dy \right)\, .
\end{eqnarray*}
Combining this with Lemmas \ref{l:calc}--\ref{l:3} we get that
there exists a positive constant, which by slight abuse of notation we denote by $C_{48}>0$, such that
\begin{equation}\label{e:bigger14-2}
L^\sB_{\alpha}\psi(z)\le C_{48}\left(z_d^{2-\alpha}+z_d^{\beta_1}+z_d^{1-\alpha}\right)|\log z_d|^{\beta_3+1}\, , \quad z\in U\, .
\end{equation}

Let
$$
f(z_d):=\left(z_d^2+z_d^{\beta_1+\alpha}+z_d\right)|\log z_d|^{\beta_3+1}\, .
$$
Assume that $|\wt{z}|\in (4^{-1}, 1)$. By the assumption on $\psi$, we have that $\psi(z)=\psi(\wt{z}, z_d)\ge 4^{-\gamma}$. Since $\lim_{z_d\to 0}f(z_d)=0$,  we can choose $\kappa >0$ so that
$$
f(z_d)-\frac{C(\alpha, p, \sB)4^{-\gamma}}{C_{48}}\le 0
$$
for all $z_d\in (0,\kappa)$. Then,
\begin{eqnarray*}
L^\sB\psi(z)&=&L^\sB_{\alpha}\psi(z)-C(\alpha,p,\sB)z_d^{-\alpha}\psi(z)\le L^\sB_{\alpha}\psi(z)-C(\alpha,p,\sB)4^{-\gamma}z_d^{-\alpha}\\
&\le &C_{48}z_d^{-\alpha}\left(f(z_d)- \frac{C(\alpha, p, \sB)4^{-\gamma}}{C_{48}}\right)\le 0
\end{eqnarray*}
for all $z\in U\setminus D(2^{-2},2^{-2})$ with $|\wt{z}|\in  (4^{-1}, 1)$ and $z_d\in (0, \kappa)$. 

Finally, it follows from \eqref{e:bigger14-1} and \eqref{e:bigger14-2} that there exists $C_{47}=C_{47}(\kappa)$ such that $L^\sB\psi(z)\le L^\sB_{\alpha}\psi(z)\le C_{47}$ for all $z\in U$ with $z_d\ge \kappa$. \qed

\section{Proofs of Lemmas 5.9 and 5.11--5.13}\label{s:10}
Before we give the proofs of these lemmas, we first do some preparation.

Let $v\in C^{\infty}_c(\R^d)$ be a non-negative smooth radial function such that $v(y)=0$ for $|y|\ge 1$ and $\int_{\R^d}v(y)\, dy=1$. For any $a>0$, set $D(a):=\{x=(\widetilde{x}, x_d)\in \R^d: x_d>a\}$ and $U_a(r):= \{y\in U(r): \delta_{U(r)}>a\}$. We write $U_a$ for $U_a(1)$ and recall that $U=D(1/2,1/2)$.
For $b\ge 10$ and $k\ge 1$,  set $v_k(y):=b^{kd} v(b^ky)$. Next we define
$g_k:=v_k\ast(g{\bf 1}_{D(5^{-k})})$ for a bounded, compactly supported function $g$ vanishing on $\R^d\setminus  \R^d_+$.
Since $b^{-k}<5^{-k}$, we have $g_k\in C^\infty_c(\R^d_+)$ and hence $L^\sB g_k$ is defined everywhere. Also note that $v_k\ast g \in C^{\infty}_c(\R^d_+; \R^d)$ and thus $L^\sB (v_k \ast g)$ is well defined (cf.~Subsection 3.2). 

Let $(a_k)_{k\ge 1}$ be a decreasing sequence of positive numbers such that $\lim_{k\to \infty}a_k=0$ and  
$$
a_k\ge 2^{-k(\beta_1/2+1)/(1+\alpha+3\beta_1/2)}\ge 2^{-k}.
$$

\begin{lemma}\label{l:Lvk}
Let $g:\R^d \to [0,1]$ be a compactly supported Borel function vanishing on $\R^d\setminus  \R^d_+$. For any $z\in U$ it holds that
\begin{equation}\label{e:Lvk-1}
\lim_{k\to \infty}L^\sB(v_k\ast g -g_k)(z)=0\, .
\end{equation}
Moreover,  there exists $C_{49}>0$ independent of $g$ such that  for all $k \ge 2$ and $z\in U_{a_k}$
\begin{equation}\label{e:Lvk-2}
0\le L^\sB(v_k\ast g -g_k)(z) \le C_{49} (2/3)^{k(\beta_1/2+1)} z_d^{\beta_1}\, .
\end{equation}
\end{lemma}
\pf 
Let  $z\in U_{a_k}$. We first estimate the difference 
\begin{eqnarray*}
L^\sB(v_k\ast g -g_k)(z)&=&\lim_{\epsilon \to 0}\int_{\R^d_+, |y-z|>\epsilon}\frac{((v_k\ast g)(y)-g_k(y))-((v_k\ast g)(z)-g_k(z))}{|y-z|^{d+\alpha}}\sB(y,z)\, dy \\
& &  -C(\alpha, p, \sB)z_d^{-\alpha}(v_k\ast g -g_k)(z).
\end{eqnarray*}
Note that for $k\ge 2$, $u\in B(0, b^{-k})$ and $y\in \R^d_+$ with $y_d>3^{-k}$, $y_d-u_d>3^{-k}-10^{-k}>5^{-k}$. Therefore
\begin{align}\label{e:3.11}
1-{\bf 1}_{D(5^{-k})}(y-u)=0.
\end{align}
Since $v_k$ is supported in $B(0, b^{-k})$, for all $k\ge 2$ and $z\in \R^d_+$ with $z_d> a_k>2^{-k}$,
$$
\int_{\R^d}(1-{\bf 1}_{D(5^{-k})}(z-u))g (z-u)v_k(u)du=0.
$$
Thus $(v_k\ast g -g_k)(z)=0$. 
Because of the same reason we have that for $z\in U_{a_k}$,
\begin{align*}
&\int_{\R^d_+, |y-z|>\epsilon}\frac{\big((v_k\ast g)(y)-g_k(y)\big)-\big((v_k\ast g)(z)-g_k(z)\big)}{|y-z|^{d+\alpha}}\sB(y,z)\, dy\\
& =\int_{\R^d_+, |y-z|>\epsilon, y_d\le 3^{-k}}\int_{\R^d} v_k(u)
\frac{(1-{\bf 1}_{D(5^{-k})})(y-u)g(y-u)}{|y-z|^{d+\alpha}}\, du\sB(y,z)\, dy \\
&\le \int_{\R^d}v_k(u)\, du \int_{\R^d, y_d\le 3^{-k}}\frac{\sB(y,z)}{|y-z|^{d+\alpha}}\, dy\\
&\le c_1\int_{y_d\le 3^{-k}}\frac{1}{|y-z|^{d+\alpha}}\left(\frac{y_d}{|y-z|}\right)^{\beta_1/2}\, dy\\
&\le c_2 (3^{-k})^{\beta_1/2+1}\int_0^{\infty}\frac{t^{d-2}}{(t^2+cz_d^2)^{(d+\alpha+\beta_1/2)/2}}\, dt \\
&= c_3 (3^{-k})^{\beta_1/2+1} z_d^{-1-\alpha-\beta_1/2} \int_0^{\infty}\frac{s^{d-2}}{(s^2+1)^{(d+\alpha+\beta_1/2)/2}}\, ds \\
&\le c_4 (2/3)^{k(\beta_1/2+1)} z_d^{\beta_1}\, .
\end{align*}
In the third line we used that $0\le g\le 1$,  in the fourth the fact that (together with \eqref{e:B7_2}) 
$$
\left(\frac{y_d\wedge z_d}{|y-z|}\wedge 1\right)^{\beta_1}\log\left(1+\frac{(y_d\vee z_d)\wedge|y-z|}{y_d\wedge z_d \wedge|y-z|}\right)^{\beta_3}
\le c \left(\frac{y_d}{|y-z|}\wedge 1\right)^{\beta_1/2},
$$
in the fifth integration in polar coordinates in $\R^{d-1}$, in the sixth the change of variables $t=c^{1/2}z_d s$, and in the last line the fact that $2^{-k(\beta_1/2+1)}z_d^{-1-\alpha-3\beta_1/2}\le 1$ which follows from $z_d\ge a_k$ and the choice of $a_k$.  Note also that it is clear from the second line that the first line is non-negative. Thus by letting $\epsilon \to 0$ we get for $z\in U_{a_k}$,
$$
0\le L^\sB(v_k\ast g -g_k)(z) \le c_4 (2/3)^{k(\beta_1/2+1)} z_d^{\beta_1}\, .
$$
Now take $z\in U$. Then there exists $k_0\ge 1$ such that $z\in U_{a_k}$ for all $k\ge k_0$, and it follows from above that
$$
\lim_{k\to \infty}L^\sB(v_k\ast g -g_k)(z)=0\, .
$$
\qed

\begin{lemma}\label{l:vLg}
Assume that $g:\R^d_+\to [0,1]$ is a function which is $C^2$  on $D(1,1)$. 
For any $k\ge 2$, $z\in U_{a_k}$ and $|u|<b^{-k}$,
\begin{equation}\label{e:vLg-0}
\mathrm{p.v.} \int_{\R^d_+} \frac{g(y-u)-g(z-u)}{|y-z|^{d+\alpha}}\sB(y,z)\, dy
\end{equation}
is well defined. Moreover, for $z\in U_{a_k}$,
\begin{equation}\label{e:vLg-1}
L^\sB_{\alpha}(v_k\ast g)(z)=\int_{\R^d}v_k(u)\left( \mathrm{p.v.}  \int_{\R^d_+} \frac{g(y-u)-g(z-u)}{|y-z|^{d+\alpha}}\sB(y,z)\, dy\right)\, du\, ,
\end{equation}
and there exists $C_{50}(z)>0$ such that $|L^\sB_{\alpha} (v_k\ast g)(z)|\le C_{50}(z)$ for all $k\ge 2$.
\end{lemma}
\pf Let $z\in U_{a_k}$ and $|u|<b^{-k}$. 
Let $G(y, z, u):=(g(y-u)-g(z-u))|y-z|^{-d-\alpha}$.
For $0<\epsilon <\eta <z_d/10$, consider
\begin{align*}
&\int_{\R^d_+, \epsilon < |y-z| }G(y, z, u)\sB(y,z)\, dy-\int_{\R^d_+, \eta < |y-z| }G(y, z, u)\sB(y,z)\, dy\\
&=\int_{\R^d_+, \epsilon < |y-z| < \eta}G(y, z, u)\sB(y,z)\, dy\\
&= \int_{ \epsilon < |y-z| <\eta}G(y, z, u)\, dy+ \int_{ \epsilon < |y-z| <\eta}G(y, z, u)(\sB(y,z)-1)\, dy\\
&=: I+II\, .
\end{align*}
Since $g$ is $C^2$ on $D(1,1)$ and $y-u, z-u\in D(1,1)$, we see that
\begin{align*}
&| I |\le \int_{\epsilon < |(y-u)-(z-u)|<\eta} \frac{|g(y-u)-g(z-u)-\nabla g(z-u){\bf 1}_{(|(y-u)-(z-u)|<1)}\cdot (y-z)|}{|(y-u)-(z-u)|^{d+\alpha}}\, dy\\
&\le c_1 \sup_{w\in B(z, z_d/5)}|\partial^2 g (w)| \int_{ \epsilon < |y-z| <\eta} |y-z|^{-d-\alpha+2}\, dy=c_2(z) (\eta^{2-\alpha}-\epsilon^{2-\alpha})\, .
\end{align*}
Further, by using the mean value theorem in the first line and \textbf{(B4)}  in the second, we get 
\begin{align*}
&| II |\le  \sup_{w\in B(z, z_d/5)}|\nabla g (w)| \int_{\epsilon < |y-z| <\eta} \frac{|\sB(y,z)-1|}{|y-z|^{d+\alpha-1}}\, dy\\
&\le c_3(z)  \int_{\epsilon < |y-z| <\eta} |y-z|^{-d-\alpha}\left(\frac{|y-z|}{y_d\wedge z_d}\right)^{\theta}\, dy\\
&\le c_4 c_3(z) z_d^{-\theta} \int_{\epsilon < |y-z| <\eta} |y-z|^{-d-\alpha+\theta}\, dy=c_5(z) (\eta^{\theta-\alpha+1}-\epsilon^{\theta-\alpha+1})\, .
\end{align*}
The estimates  for $I$ and $II$ imply that the principal value integral in \eqref{e:vLg-0} is well defined.

Let $z\in U_{a_k}$. For $\epsilon < z_d/10$ and $|u|<b^{-k}$, we have
\begin{align*}
&\left| \int_{\R^d_+, |y-z|>\epsilon} G(y, z, u)\sB(y,z)\, dy\right|\\
&\le \left|\int_{\R^d_+,  |y-z|\ge z_d/10} G(y, z, u)\sB(y,z)\, dy\right| +\left|\int_{\R^d_+, \epsilon < |y-z| < z_d/10} G(y, z, u)\sB(y,z)\, dy\right| \\
&=:III+IV\, .
\end{align*}
By estimating $g$ by 1, we get that
$$
III\le 2\int_{|y-z|\ge z_d/10}|y-z|^{-d-\alpha}dy \le c_6 z_d^{-\alpha}=c_7(z) \, .
$$
The integral in $IV$ is estimated in $I$ and $II$ with $\eta=z_d/10$, so we have 
$$
IV\le c_2(z)(z_d/10)^{2-\alpha}+ c_5(z)(z_d/10)^{\theta-\alpha +1}= c_8(z)\, .
$$
Thus we have that
\begin{equation}\label{e:vLg-2}
\left| \int_{\R^d_+, |y-z|>\epsilon} \frac{g(y-u)-g(z-u)}{|y-z|^{d+\alpha}}\sB(y,z)\, dy\right| \le c_9(z)\, .
\end{equation}
Hence we can use the dominated convergence theorem to conclude that
\begin{align*}
&L^\sB_{\alpha} (v_k\ast g)(z) =\lim_{\epsilon \to 0}\int_{\R^d_+, |y-z|>\epsilon}\frac{(v_k\ast g)(y)-(v_k\ast g)(z)}{|y-z|^{d+\alpha}}\sB(y,z)\, dy\\
&=\lim_{\epsilon\to 0} \int_{|u|<b^{-k}}v_k(u) \int_{\R^d_+, |y-z|>\epsilon} \frac{g(y-u)-g(z-u)}{|y-z|^{d+\alpha}}\sB(y,z)\, dy \, du\\
&=\int_{|u|<b^{-k}} v_k(u)\left(\lim_{\epsilon \to 0} \int_{\R^d_+, |y-z|>\epsilon} \frac{g(y-u)-g(z-u)}{|y-z|^{d+\alpha}}\sB(y,z)\, dy\right)\, du\, ,
\end{align*}
which is \eqref{e:vLg-1}. The last statement follows from \eqref{e:vLg-2} . 
\qed

We note also that since $g$ is continuous in $D(1, 1)$, it holds that $\lim_{k\to \infty}(v_k\ast g)(z)=g(z) $ for all $z\in D(1, 1)$.

Recall that the function $h_q$ was defined in Section \ref{s:estimates}: $h_q(x)=x_d^q {\bf 1}_{D(1,1)}$. 
\begin{lemma}\label{l:hq}
Let $q\in ( (\alpha-1)_+, \alpha+\beta_1)$ 
and set $b:=10\vee 2^{4(q-2)_-+3}$ and 
$$
a_k:=2^{-k(q+1+\frac12\beta_1)/(\alpha+1+\frac32\beta_1)}\vee 2^{-k(2+\beta_1)/(1+\alpha+\frac32\beta_1-q)}.
$$
There exists a constant $C_{51}>0$ such that for any $k\ge 1$ and $z\in U_{a_k}$,
\begin{equation}\label{e:hq-1}
|L^\sB(v_k \ast h_q)(z)-L^\sB h_q(z)|\le C_{51} \left(\frac{4}{5}\right)^k \, .
\end{equation}
In particular, the functions $z\to |L^\sB(v_k \ast h_q)(z)-L^\sB h_q(z)|$ are all bounded
by the constant $C_{51}$on $U$, and for any $z\in U$, 
$
\lim_{k\to \infty} \left| L^\sB(v_k \ast h_q)(z)-L^\sB h_q(z)\right| =0.
$
\end{lemma}
\pf First note that $a_k\ge 2^{-k}$ since the first term in its definition is larger than $2^{-k}$. Next, by using Lemma \ref{l:vLg} with $g=h_q$ in the third line below we see that
\begin{align*}
&L^\sB(v_k\ast h_q)(z)\\
&=L^\sB_{\alpha}(v_k\ast h_q)(z)-C(\alpha, p, \sB)z_d^{-\alpha}(v_k\ast h_q)(z)\\
&=\int_{\R^d}v_k(u)\left(\mathrm{p.v.}\int_{\R^d_+}\frac{h_q(y-u)- h_q(z-u)}{|y-z|^{d+\alpha}}\sB(y, z)dy\right)du
-C(\alpha, p, \sB)z_d^{-\alpha}(v_k\ast h_q)(z)\\
&=\int_{\R^d}v_k(u)\left(\mathrm{p.v.}\int_{\R^d_+}\frac{h_q(y-u)- h_q(z-u)-(h_q(y)-h_q(z))}{|y-z|^{d+\alpha}}\sB(y, z)dy\right)du\\
&\quad +C(\alpha, p, \sB)z_d^{-\alpha}(h_q(z)-(v_k\ast h_q)(z)) + L^\sB h_q(z).
\end{align*}
Note that for $z\in U_{a_k}$, 
\begin{align}
|h_q(z)-(v_k\ast h_q)(z)|=\left|\int_{B(0, b^{-k})}v_k(u)(h_q(z)-h_q(z-u))du\right|\le c_0(1/5)^k z_d^q.\label{e:new1-1}
\end{align}

 Set $b=10\vee 2^{4(q-2)_-+3}$. Now we write, for $u\in B(0, b^{-k})$,
\begin{align*}
&
\mathrm{p.v.} \int_{\R^d_+}\frac{h_q(y-u)- h_q(z-u)-(h_q(y)-h_q(z))}{|y-z|^{d+\alpha}}\sB(y, z)dy\\
&=\int_{D(1+b^{-k}, 1+b^{-k})\setminus U, y_d>5^{-k}}+\int_{D(1+b^{-k}, 1+b^{-k}), y_d<5^{-k}}\\
&\quad + 
 \int_{U, y_d>5^{-k}, |y-z|>2^{-1}z_d}+\  \mathrm{p.v.}\int_{U \cap B(z, z_d/2)} =:I+II+III+IV.
\end{align*}

We deal with $I$ first. For $u\in B(0, b^{-k})$,
\begin{align*}
I&=\int_{D(1+b^{-k}, 1+b^{-k})\setminus D(1-b^{-k}, 1-b^{-k}), y_d>5^{-k}}+
\int_{D(1-b^{-k}, 1-b^{-k})\setminus U, y_d>5^{-k}}=:I_1+I_2.
\end{align*}
Obviously, we have
$
|I_1|\le c_1 b^{-k}.
$

Let $A_k:=(D(1-b^{-k}, 1-b^{-k})\setminus U)\cap\{y:y_d>5^{-k}\}$. 
Then, we have
\begin{align*}
&|I_2|=\big|\int_{A_k}
\frac{(y_d-u_d)^q-y_d^q-((z_d-u_d)^q-z_d^q)}{|y-z|^{d+\alpha}}
\sB(y, z)dy\big|\\
&=\big|\int_{A_k}\frac{q(y_d-z_d)\cdot \int^1_0\left((z_d-u_d+ t (y_d-z_d))^{q-1}-(z_d+ t (y_d-z_d ))^{q-1}  \right) d t }{|y-z|^{d+\alpha}}\sB(y, z)dy\big|\\
&
\le c|u_d| \int_{A_k}\frac{(z_d \wedge y_d)^{-(q-2)_-} |y_d-z_d|}{|y-z|^{d+\alpha}}dy
\le c b^{-k}5^{k(q-2)_-}\int_{A_k}\frac{|y_d-z_d|
}
{|y-z|^{d+\alpha}} dy\\
&
\le c b^{-k}2^{3k(q-2)_-}\int_{2>|y-z|>a_k}\frac{1}{|y-z|^{ d+\alpha-1 }} dy 
\le c b^{-k}2^{3k(q-2)_-} a_k^{-1}
\le c_2 b^{-k}2^{3k(q-2)_-+1}\le c_2 2^{-k((q-2)_{-}+2)},
\end{align*}
where in the first inequality we use mean value theorem  for the difference  inside the integral in the numerator, in the second last inequality, the fact $a_k\ge 2^{-k}$ and  in the last inequality, the fact $b\ge 2^{4(q-2)_- +3}$.

For $II$, we have
$$
|II| \le c_3 \int_{D(1+b^{-k}, 1+b^{-k}), y_d<5^{-k}}\frac{ y_d^q+z_d^q+b^{-kq}}{|y-z|^{d+\alpha}}\sB(y, z)dy.
$$
Similarly as in the proof of Lemma \ref{l:Lvk}  we first estimate 
\begin{align*}
&\int_{D(1+b^{-k}, 1+b^{-k}), y_d<5^{-k}}\frac{y_d^q+b^{-kq}}{|y-z|^{d+\alpha}}\sB(y, z)dy\\
&\le c_4\int_{\R^{d-1}}\int^{5^{-k}}_0\frac{y_d^{q+\beta_1/2}+10^{-kq}y_d^{\beta_1/2}}{|y-z|^{d+\alpha+\frac12\beta_1}}dy_dd{\widetilde y}\\
&\le c_5 5^{-k(q+1+\frac12\beta_1)}\int^\infty_0\frac{t^{d-2}}{(t^{2}+cz_d^2)^{(d+\alpha+\frac12\beta_1)/2}}dt \\
&= c_6 5^{-k(q+1+\frac12\beta_1)}z_d^{-1-\alpha-\frac12\beta_1}\int^\infty_0\frac{s^{d-2}}{(s^{2}+1)^{(d+\alpha+\frac12\beta_1)/2}}ds\\
& \le c_7 (2/5)^{k(q+1+\frac12\beta_1)}z_d^{\beta_1}.
\end{align*}
For the remaining part, we use a similar argument:
\begin{align*}
&\int_{D(1+b^{-k}, 1+b^{-k}), y_d<5^{-k}}\frac{z_d^q}{|y-z|^{d+\alpha}}\sB(y, z)dy\\
&\le c_8 z_d^q\int_{\R^{d-1}}\int^{5^{-k}}_0\frac{y_d^{\beta_1/2}}{|y-z|^{d+\alpha+\frac12\beta_1}}dy_dd\widetilde{y}\\
&\le c_9 z_d^q 5^{-k(1+\frac12\beta_1)}\int^\infty_0\frac{t^{d-2}}{(t^2+cz^2_d)^{(d+\alpha+\frac12\beta_1)/2}}dt\\
&=c_{10} z_d^q 5^{-k(1+\frac12\beta_1)}z_d^{-1-\alpha-\frac12\beta_1}\int^\infty_0\frac{s^{d-2}}{(s^2+1)^{(d+\alpha+\frac12\beta_1)/2}}ds \\
&\le c_{11} (4/5)^{k(1+\frac12\beta_1)}2^{-k(2+\beta_1)}z_d^{q-1-\alpha-\frac12\beta_1}\le c_{12} (4/5)^{k(1+\frac12\beta_1)}z_d^{\beta_1}.
\end{align*}
Thus
$$
|II|\le c_{13} \big((2/5)^{k(q+1+\frac12 \beta_1)}+(4/5)^{k(1+\frac12\beta_1)}\big)z_d^{\beta_1}.
$$

Let $B_k:=U\cap \{y_d>5^{-k}\}\cap \{y:|y-z|>2^{-1}z_d\}$. Then, we have
\begin{align*}
&|III|=\left|\int_{B_k}\frac{(y_d-u_d)^q-y_d^q-((z_d-u_d)^q-z_d^q)}{|y-z|^{d+\alpha}}
\sB(y, z)dy\right|\\
&=\left|\int_{B_k}\frac{q(y_d-z_d)\cdot\int^1_0\left((z_d-u_d+ t (y_d-z_s))^{q-1}-(z_d+t (y_d-z_d))^{q-1}\right)d t}{|y-z|^{d+\alpha}}
\sB(y, z)dy\right|\\
&\le c_{14}b^{-k}2^{k3(q-2)_-}\int_{U, |y-z|>2^{-1}z_d}\frac{1}{|y-z|^{d+\alpha-1}}dy\\
&\le c_{15}(4/5)^{k}2^{-3k}(z_d^{1-\alpha}\vee \log \frac1{z_d})
\le c_{16} (4/5)^{k}z_d^3(z_d^{1-\alpha}\vee \log \frac1{z_d}),
\end{align*}
where in the first inequality we use mean value theorem inside the integral in the numerator
and the fact the derivative of the integrand is bounded above by $c(5^{-k}-b^{-k})^{-(q-2)_-}\le c2^{3k(p-2)_-}$.

Let $F(y_d, z_d, u_d):=q(q-1)\int^1_0\left((z_d-u_d+ t (y_d-z_d))^{q-2}-(z_d+ t (y_d-z_d))^{q-2}\right)(1- t )d t $. 
For $t \in [0, 1]$, $u\in B(0, b^{-k})$
and $y\in B(z, \frac12 z_d)$, $z_d-u_d+ t (y_d-z_d)$ and $z_d+ t (y_d-z_d)$ are both comparable with $z_d$. Thus,
for $IV$, we have that, when $\alpha \ge 1$,  for large $k$,
\begin{align*}
&|IV| \le \left|
 \mathrm{p.v.} 
\int_{U\cap B(z, 2^{-1}z_d)}\frac{(y_d-u_d)^q-y_d^q-((z_d-u_d)^q-z_d^q)}{|y-z|^{d+\alpha}}
(\sB(y, z)-1)dy\right|\\
&+\left|\int_{U\cap B(z, 2^{-1}z_d)}
\frac{(y_d-u_d)^q-y_d^q-((z_d-u_d)^q-z_d^q) -q((z_d-u_d)^{q-1}-z_d^{q-1})
{\bf 1}_{B(z, a_k)} (y) (y_d-z_d)}{|y-z|^{d+\alpha}}
dy\right|\\
&\le c_{17}
\int_{U\cap B(z, 2^{-1}z_d)}
\frac{\left|\int^1_0((z_d-u_d+t(y_d-z_s))^{q-1}-(z_d+t(y_d-z_d))^{q-1}dt\right|}{z_d^{\theta} |y-z|^{d+\alpha-1- \theta }}
dy\\
& \quad + c_{18}
\int_{U\cap (B(z, 2^{-1}z_d) \setminus B(z, a_k) )}
\frac{\left|\int^1_0((z_d-u_d+t(y_d-z_s))^{q-1}-(z_d+t(y_d-z_d))^{q-1}dt\right|}{ |y-z|^{d+\alpha-1}}
dy\\
& \quad +\int_{B(z, a_k)}\frac{\left| F(y_d, z_d, u_d) \right|}{|y-z|^{d+\alpha-2}}
dy\\
&\le c_{19} \left(
\int_{B(z, 2^{-1}z_d)}
\frac{ |u_d|z_d^{q-2}
dy}
{z_d^\theta |y-z|^{d+\alpha-1-\theta}}
+ 
\int_{B(z, 2^{-1}z_d) \setminus B(z, a_k) }
\frac{
 |u_d|z_d^{q-2}dy
}{ |y-z|^{d+\alpha-1}}
+\int_{B(z, a_k)}\frac{ |u_d|z_d^{q-3}dy}{|y-z|^{d+\alpha-2}}
 \right)\\
&\le c_{ 20 }b^{-k} \left ( z_d^{q-1-\alpha} +z_d^{q-2 }
a_k^{1-\alpha}
+
z_d^{q-1-\alpha}  \right)   
\le c_{21} (2/b)^k  z_d^{q-1-\alpha }
\le c_{ 22 } (4/5)^k z_d^{q+1-\alpha },
\end{align*}
where  in the second inequality we used \textbf{(B4)}  and   in the third inequality 
the mean value theorem
for the difference inside the integral in the definition of $F$. 
 When $\alpha < 1$, the argument is much simpler and we skip it. 

By putting together \eqref{e:new1-1} and the estimates for $I$, $II$ , $III$ and $IV$ we see that \eqref{e:hq-1} is true for some constant $C$ (independent of $z$ and $k\ge 1$). \qed

Recall that the function $\psi$ in the next lemma was introduced in Section \ref{s:proof-of-key-lemma}.
\begin{lemma}\label{l:psi}
Let $(a_k)_{k\ge 1}$ be a decreasing sequence of positive numbers such that $\lim_{k\to \infty}a_k=0$ and  
$a_k\ge 2^{-k(2+\beta_1)/(2+\alpha+\frac32\beta_1)}$. 
There exists a constant $C_{52}>0$ such that for every $k\ge 1$ and $z\in U_{a_k}$
\begin{equation}\label{e:psi-1}
|L^{\sB}(v_k \ast \psi)(z)-L^{\sB} \psi(z)|\le C_{52} \left(\frac{4}{5}\right)^k\, .
\end{equation}
In particular, for any $z\in U$, 
$
\lim_{k\to \infty} \left| L^{\sB}(v_k \ast \psi)(z)-L^{\sB} \psi(z)\right| =0\, .
$
\end{lemma}
\pf The proof is similar to the proof of Lemma \ref{l:hq}. 
By using Lemma \ref{l:vLg} with $g=\psi$ in the third line below we see that
\begin{align*}
&L^{\sB}(v_k\ast \psi)(z)\\
&=L^{\sB}_{\alpha}(v_k\ast \psi)(z)-C(\alpha, p, {\sB})z_d^{-\alpha}(v_k\ast \psi)(z)\\
&=\int_{\R^d}v_k(u)\left(\mathrm{p.v.}\int_{\R^d_+}\frac{\psi(y-u)- \psi(z-u)}{|y-z|^{d+\alpha}}{\sB}(y, z)dy\right)du
-C(\alpha, p, {\sB})z_d^{-\alpha}(v_k\ast \psi)(z)\\
&=\int_{\R^d}v_k(u)\left(\mathrm{p.v.}\int_{\R^d_+}\frac{\psi(y-u)- \psi(z-u)-(\psi(y)-\psi(z))}{|y-z|^{d+\alpha}}{\sB}(y, z)dy\right)du\\
&\quad +C(\alpha, p, {\sB})z_d^{-\alpha}(\psi(z)-(v_k\ast \psi)(z)) + L^{\sB} \psi(z).
\end{align*}
 Set $b=10\vee 2^{4(q-2)_-+3}$. Note that for $z\in U_{a_k}$, 
\begin{align}
&|\psi(z)-(v_k\ast \psi)(z)|=\big|\int_{B(0, b^{-k})}v_k(u)(\psi(z)-\psi(z-u))du\big|
\le cb^{-k}.\label{e:new1-2}
\end{align}

Now we write, for $u\in B(0, b^{-k})$,
\begin{align*}
&
\mathrm{p.v.}\int_{\R^d_+}\frac{\psi(y-u)- \psi(z-u)-(\psi(y)-\psi(z))}{|y-z|^{d+\alpha}}{\sB}(y, z)dy\\
&=\int_{D(3+b^{-k}, 3+b^{-k})\setminus U, y_d>5^{-k}}+\int_{D(3+b^{-k}, 3+b^{-k}), y_d<5^{-k}}\\
&
\quad + 
 \int_{U, y_d>5^{-k}, |y-z|>2^{-1}z_d}+\ \mathrm{p.v.}\int_{U \cap B(z, z_d/2)} 
=:I+II+III+IV.
\end{align*}

We deal with $I$ first. For $u\in B(0, b^{-k})$,
\begin{align*}
I&=\int_{D(3+b^{-k}, 3+b^{-k})\setminus D(3-b^{-k}, 3-b^{-k}), y_d>5^{-k}}+
\int_{D(3-b^{-k}, 3-b^{-k})\setminus U, y_d>5^{-k}}=:I_1+I_2.
\end{align*}

Obviously, we have
$
|I_1|\le c b^{-k}.
$
Let $A_k:=(D(3-b^{-k}, 3-b^{-k})\setminus U)\cap\{y:y_d>5^{-k}\}$. 
Then, we have
\begin{align*}
&|I_2|=\big|\int_{A_k}
\frac{\psi(y-u)-\psi(y)-(\psi(z-u)-\psi(z))}{|y-z|^{d+\alpha}}
{\sB}(y, z)dy\big|\\
&=\big|\int_{A_k}\frac{(y-z)\cdot \int^1_0\left(\nabla \psi(z-u+ t (y-z))-\nabla \psi(z+ t (y-z))\right)d t}{|y-z|^{d+\alpha}}
{\sB}(y, z)dy\big|\le cb^{-k},
\end{align*}
where in the first inequality we use mean value theorem inside the integral in the numerator.

For $II$, we have
\begin{align*}
|II|&\le c\int_{D(3+b^{-k}, 3+b^{-k}), y_d<5^{-k}}\frac{1}{|y-z|^{d+\alpha}}
{\sB}(y, z)dy\\
&\le c\int_{\R^{d-1}}\int^{5^{-k}}_0\frac{y_d^{\beta_1/2}}{|y-z|^{d+\alpha+\frac12\beta_1}}dy_dd\widetilde{y}\\
&\le c5^{-k(1+\frac12\beta_1)}\int^\infty_0\frac{t^{d-2}}{(t^2+cz^2_d)^{(d+\alpha+\frac12\beta_1)/2}}dt\\
&=c5^{-k(1+\frac12\beta_1)}z_d^{-1-\alpha-\frac12\beta_1}\int^\infty_0\frac{s^{d-2}}
{(s^2+1)^{(d+\alpha+\frac12\beta_1)/2}}ds\\
&\le c(4/5)^{k(1+\frac12\beta_1)}2^{-k(2+\beta_1)}z_d^{-1-\alpha-\frac12\beta_1}
\le c(4/5)^{k(1+\frac12\beta_1)}z_d^{1+\beta_1},
\end{align*}
since $a_k\ge 2^{-k(2+\beta_1)/(2+\alpha+\frac32\beta_1)}$.

Let $B_k:=U\cap \{y_d>5^{-k}\}\cap \{y:|y-z|>2^{-1}z_d\}$. Then, we have
\begin{align*}
&|III|=\big|\int_{B_k}\frac{\psi(y-u)-\psi(y)-(\psi(z-u)-\psi(z)  )}{|y-z|^{d+\alpha}}
B(y, z)dy\big|\\
&=\big|\int_{B_k}\frac{(y-z)\cdot\int^1_0\left(\nabla \psi(z-u+t(y-z))-\nabla \psi(z+ t (y-z))\right)d t}{|y-z|^{d+\alpha}}
{\sB}(y, z)dy\big|\\
&\le cb^{-k}\int_{U, |y-z|>2^{-1}z_d}\frac{1}{|y-z|^{d+\alpha-1}}dy\\
&\le c(4/5)^{k}2^{-3k}(z_d^{1-\alpha}\vee \log \frac1{z_d})
\le c (4/5)^{k}z_d^3(z_d^{1-\alpha}\vee \log \frac1{z_d}),
\end{align*}
where in the first inequality we use mean value theorem for the difference inside the integral in the numerator
and the fact the derivative of the integrand is bounded.

Let $F(y, z,  u  ):=\int^1_0\left(\nabla^2\psi(z-u+ t (y-z))-\nabla^2\psi(z+ t (y-z))\right)(1- t )d t $. For $IV$, we have
that, when $\alpha \ge 1$,  
for large $k$, 
\begin{align*}
&|IV| \le \left|
 \mathrm{p.v.} \int_{U\cap B(z, 2^{-1}z_d)}\frac{\psi(y-u)-\psi(y)-(\psi(z-u)-\psi(z)  )}{|y-z|^{d+\alpha}}
(\sB(y, z)-1)dy\right|\\
&+\left|\int_{U\cap B(z, 2^{-1}z_d)}\!\!\!\!\!\!\!\!\!\!\!\!\!\!\!
\frac{\psi(y-u)-\psi(y)-(\psi(z-u)-\psi(z)  )-(\nabla \psi(z-u)-\nabla \psi(z))
{\bf 1}_{B(z, a_k)} (y) (y_d-z_d)}{|y-z|^{d+\alpha}}
dy\right|\\
&\le c_{17}
\int_{U\cap B(z, 2^{-1}z_d)} 
\frac{\left|\int^1_0\left(\nabla \psi(z-u+t(y-z))-\nabla \psi(z+t(y-z))\right)dt
\right|}{z_d^\theta |y-z|^{d+\alpha-1-\theta}}
dy\\
&\quad + c_{18}
\int_{U\cap (B(z, 2^{-1}z_d) \setminus B(z, a_k) )}
\frac{\left|\int^1_0\left(\nabla \psi(z-u+t(y-z))-\nabla \psi(z+t(y-z))\right)dt\right|}{ |y-z|^{d+\alpha-1}}
dy\\
&\quad +\int_{B(z, a_k)}\frac{\left|(y-z)\otimes (y-z)\cdot  F(y, z, u) \right|}{|y-z|^{d+\alpha}}
dy\\
&\le c_{19} \left(
\int_{B(z, 2^{-1}z_d)}
\frac{ |u_d|
dy}
{z_d^\theta |y-z|^{d+\alpha-1-\theta}}
+ 
\int_{B(z, 2^{-1}z_d) \setminus B(z, a_k) }
\frac{
 |u_d|dy
}{ |y-z|^{d+\alpha-1}}
+\int_{B(z, a_k)}\frac{ |u_d|dy}{|y-z|^{d+\alpha-2}}
 \right)\\
&\le c_{20 }b^{-k} \left ( z_d^{1-\alpha} + a_k^{1-\alpha}+
a_k^{2-\alpha}  \right)   
\le  c_{21}   (2/b)^k \le c_{ 22 } (4/5)^k,
\end{align*}
where  in the second inequality we used \textbf{(B4)} and  in the third inequality 
mean value theorem for the difference inside the integral in the definition of $F$. 
When $\alpha < 1$, the argument is much simpler and we skip it.

By putting together \eqref{e:new1-2} and estimates for $I$, $II$ , $III$ and $IV$ we see that \eqref{e:psi-1} is true for some constant $C$ (independent of $z$ and $k\ge 1$).
\qed

\begin{lemma}\label{l:dynkin-hp}
Let $r\le 1$. For every $x\in U(r)$ it holds that
\begin{equation}\label{e:dynkin-hp}
\E_x [h_p(Y_{ \tau_{U(r)}})]=h_p(x)+\E_x \int_0^{\tau_{U(r)}}  L^\sB h_p(Y_s)\, ds \, .
\end{equation}
\end{lemma}
\pf 
Set $g_k:=v_k\ast (h_p {\bf 1}_{D(5^{-k})})$. By combining Lemmas \ref{l:Lvk} and \ref{l:hq}  (with $q=p$), we see that for every $z\in U(r)$ with $r\le 1/2$,
$$
\lim_{k\to \infty}|L^\sB g_k (z)- L^\sB h_p(z)|=0
$$
and $|L^\sB g_k (z)- L^\sB h_p(z)|$
is bounded by the constant $C_{51}>0$. Let $x\in U(r)$, $r\le 1$. There is $k_0\ge 1$ such that $x\in U_{a_k}(r)$ for all $k\ge k_0$. Note that since $g_k\in C_c^{\infty}(\R^d_+)$, it follows from Proposition \ref{p:Dynkin} that for all $t\ge 0$,
$$
\E_x [g_k(Y_{t\wedge \tau_{U_{a_k}(r)}})]=g_k(x)+\E_x\int_0^t {\bf 1}_{s<\tau_{U_{a_k}(r)}} L^\sB g_k (Y_s)\, ds \, .
$$
Clearly, $\lim_{k\to \infty}\tau_{U_{a_k}(r)}=\tau_{U(r)}$. Since $g_k\to h_p$ as $k\to \infty$, we get that the left-hand side above converges to $\E_x [h_p(Y_{t\wedge \tau_{U(r)}})]$. Further, by Lemma \ref{l:hq},
$$
\lim_{k\to \infty}\left({\bf 1}_{s<\tau_{U_{a_k}(r)}} L^\sB g_k (Y_s)-{\bf 1}_{s<\tau_{U(r)}} L^\sB h_p(Y_s)\right)=0,
$$
and $\left| {\bf 1}_{s<\tau_{U_{a_k}(r)}} L^\sB g_k (Y_s)-{\bf 1}_{s<\tau_{U(r)}} L^\sB h_p(Y_s)\right|$
is bounded by the constant $C_{51}>0$. By Lemma \ref{l:estimate-of-L-hat-B} (a), $0\ge L^\sB h_p(z)\ge -C_{22}z_d^{\beta_1}|\log z_d|^{\beta_3}$, and hence that 
$$
\E_x \int_0^t |{\bf 1}_{s<\tau_{U(r)}} L^\sB h_p(Y_s)|\, ds <\infty\, .
$$
Thus we can use the dominated convergence theorem to conclude that 
$$
\lim_{k\to \infty} \E_x\int_0^t {\bf 1}_{s<\tau_{U_{a_k}(r)}} L^\sB g_k (Y_s)\, ds = \E_x \int_0^t {\bf 1}_{s<\tau_{U(r)}} L^\sB h_p(Y_s)\, ds \, .
$$
By letting $t\to \infty$ we obtain \eqref{e:dynkin-hp}. \qed

\noindent
\textbf{Proof of Lemma \ref{l:lower-bound-for-integral}.}
From Lemma \ref{l:key-lemma} (b) we know that $\phi(x)=x_d^p$, 
$\phi(w)\le 0$ for $w\in U^c\cap \R^d_+$ and $L^\sB \phi(z)\ge -C_{30}z_d^{\beta_1}|\log z_d|^{\beta_3}$ for all $z\in U$. In particular, $L^\sB \phi(z)\ge -C_{30}$ for all $z\in U$. 

Let $g_k=v_k\ast(\phi {\bf 1}_{D(5^{-k})})$. It follows from Lemmas \ref{l:Lvk}, \ref{l:hq} and \ref{l:psi} that $\lim_{k\to \infty}|L^\sB g_k(z)-L^\sB \phi(z)|=0$ for every $z\in U$, and the sequence of functions $|L^\sB g_k(z)-L^\sB \phi(z)|$ on $U$ is bounded by some constant $C$. Then $L^\sB g_k(z)\ge L^\sB \phi(z)-C\ge -C_{30}-C$ for every $z\in U$.

Note that since $g_k\in C_c^{\infty}(\R^d_+)$, it follows from Proposition \ref{p:Dynkin} that for all $t\ge 0$,
$$
\E_x [g_k(Y_{t\wedge \tau_{U_{a_k}}})]=g_k(x)+\E_x\int_0^t {\bf 1}_{s<\tau_{U_{a_k}}} L^\sB g_k (Y_s)\, ds \, .
$$
Clearly, $\lim_{k\to \infty}\tau_{U_{a_k}}=\tau_{U}$. Since $g_k\to \phi$ as $k\to \infty$, we get that the left-hand side above converges to $\E_x [\phi(Y_{t\wedge \tau_{U}})]$. Further
$$
\lim_{k\to \infty}\left({\bf 1}_{s<\tau_{U_{a_k}}} L^\sB g_k (Y_s)\right)={\bf 1}_{s<\tau_{U}} L^\sB\phi(Y_s)
$$
and ${\bf 1}_{s<\tau_{U_{a_k}}} L^\sB g_k (Y_s)\ge -C_{30}-C$. Thus we may use Fatou's lemma to conclude that
$$
\E_x\int_0^t {\bf 1}_{s<\tau_U} L^\sB \phi(Y_s)\, ds \le \liminf_{k\to \infty} \E_x\int_0^t {\bf 1}_{s<\tau_{U_{a_k}}} L^\sB g_k (Y_s)\, ds.
$$
By using that the right-hand side above is equal to $\E_x [\phi(Y_{t\wedge \tau_{U}})]-\phi(x)$ and then letting $t\to \infty$, we arrive at
$$
\E_x [\phi(Y_{\tau_{U}})]-\phi(x)\ge \E_x\int_0^{\tau_U} L^\sB \phi(Y_s)\, ds  \ge -C_{30}\E_x\int_0^{\tau_U} (Y_s^d)^{\beta_1}|\log Y_s^d|^{\beta_3}\, ds\, .
$$
Since $\E_x [\phi(Y_{\tau_{U}})]\le 0$ and $\phi(x)=x_d^p$, this concludes the proof. 
\qed

\noindent
\textbf{Proof of Lemma \ref{l:lower-bound-for-integral-q}.}
Let $q:=\alpha+\beta_2$. Note that $q<p$ by assumption. By repeating the proof of Lemma \ref{l:lower-bound-for-integral} (with $h_p$ replaced by $h_q$, and 
Lemma \ref{l:key-lemma}.(b) replaced by Lemma \ref{l:key-lemma}.(c)), we get
that $x=(\wt{0}, x_d)\in U(r)$,
$$
 C_{31} \E_x \int_0^{\tau_U(r)} (Y_t^d)^{\beta_2}\, dt +C_{32}  \E_x \int_0^{\tau_U(r)}
(Y_t^d)^{\beta_1} |\log Y_t^d|^{\beta_3}\, dt \ge x_d^q\, .
$$
By Lemma \ref{l:upper-bound-for-integral}, we have $\E_x \int_0^{\tau_U} (Y_t^d)^{\beta_1} |\log Y_t^d|^{\beta_3} \, dt \le C_{28} x_d^p$ for $x\in U$.
Thus for $x\in U$,
$$
\E_x \int_0^{\tau_U} (Y_t^d)^{\beta_2}\, dt \ge C_{31}^{-1}(x_d^q-C_{32} C_{28} x_d^p)\, .
$$
There exist $a_1>0$ and $C_{34}>0$ so that the last term is greater than $C_{34} x_d^q$ whenever $0<x_d<a_1$. 
By scaling, we can replace $U$ by $U(r)$. \qed

\noindent
\textbf{Proof of Lemma \ref{l:exit-time estimate-U_new}.}
Let $q:=\alpha+\beta_2$ and $\eta(x):=h_q(x)-h_p(x)$. Note that $q<p$ by assumption.
For $x\notin D(1,1)$, $\eta(x)=0$, while if $x\in D(1,1)$ we have $\eta(x)=x_d^q-x_d^p > 0$. 
Since
 $(\alpha -1)_+ <q < p$, 
by Lemma \ref{l:estimate-of-L-hat-B} (a) and (c), for all $x\in U$ we have that 
$L^\sB h_p(x)\ge -C_{22} x^{\beta_1} |\log x_d|^{\beta_3} $ 
and $L^\sB h_q(x)\le -C_{25} x^{q-\alpha}=-C_{25} x^{\beta_2}$. Thus, since $\beta_2 < \beta_1$,  we can find $r_1\in (0, 1/2]$ such that
$$
L^\sB \eta(x)=L^\sB h_q(x)-L^\sB h_p(x)\le -C_{25} x_d^{\beta_2} +
C_{22}
x_d^{\beta_1}|\log x_d|^{\beta_3}\le -2^{-1}C_{25} x_d^{\beta_2} , \qquad x\in U(r_1)
$$
and, also by Lemma \ref{l:estimate-of-L-hat-B} (c),
$$
L^\sB \eta(x)=L^\sB h_q(x)-L^\sB h_p(x)\ge -C_{26} z_d^{q-\alpha}-C_{27} z_d^{\beta_1}|\log z_d|^{\beta_3},  \qquad x\in U(r_1).
$$
Repeating the proof of Lemma \ref{l:dynkin-hp}, we get that for $x\in U(r_1)$,
$$
\E_x [\eta(Y_{ \tau_{U(r_1)}})]=\eta(x)+\E_x \int_0^{\tau_{U(r_1)}}  L^\sB \eta(Y_s)\, ds \, .
$$
Thus for $x\in U(r_1)$,
$$
\E_x [\eta(Y_{\tau_{U(r_1)}})]\le \eta(x)  -2^{-1}C_{25} \E_x\int_{0}^{\tau_U}  (Y_t^d)^{\beta_2}dt\, .
$$
Since $\eta\ge 0$ everywhere, we get $0\le x_d^q  -2^{-1}C_{25} \E_x\int_{0}^{\tau_U}  (Y_t^d)^{\beta_2}dt$ for all
$x\in U(r_1)$. 
Now \eqref{e:exit-time-estimate-U_new} follows from this and the scaling argument in \eqref{e:scaling}. \qed

\noindent
\textbf{Proof of Lemma  \ref{l:exit-time estimate-U}.}
Choose $q\in (p, \alpha)$ and let $\eta(x):=h_p(x)-h_q(x)$, $x\in \R^d_+$.
For $x\notin D(1,1)$, $\eta(x)=0$, while if $x\in D(1,1)$ we have $\eta(x)=x_d^p-x_d^q>0$. By Lemma \ref{l:estimate-of-L-hat-B}, for all $x\in U(r_0)$ we have that 
$L^\sB h_p(x)\le 0$ and $L^\sB h_q(x)\ge C_{23} x^{\alpha-q}$. Thus we can find $r_1\in (0, r_0]$ such that
\begin{equation}\label{e:exit-time-estimate-U-3}
L^\sB \eta(x)=L^\sB h_p(x)-L^\sB h_q(x)\le -C_{23} x_d^{\alpha-q}\le -1, \qquad x\in U(r_1).
\end{equation}
Let $g_k=v_k\ast (\eta {\bf 1}_{D(5^{-k}})$. It follows from Lemmas \ref{l:Lvk} and \ref{l:hq} applied to $h_g$ and $h_p$ that  $L^\sB g_k\to L^\sB \eta$ on $U$ and the sequence of functions
$|L^\sB g_k-L^\sB \eta|$  is bounded by some constant $C>0$. In particular, 
\begin{equation}\label{e:exit-time-estimate-U-2}
-L^\sB g_k(z) \ge -L^\sB \eta(z)-C \ge 1-C\, , \qquad z\in U(r), r\le r_1\, .
\end{equation}
It follows from Proposition \ref{p:Dynkin} that for all $t\ge 0$,
$$
\E_x [g_k(Y_{t\wedge \tau_{U_{a_k}(r)}})]-g_k(x)=\E_x\int_0^t {\bf 1}_{s<\tau_{U_{a_k}(r)}} L^\sB g_k (Y_s)\, ds \, .
$$
As $k\to \infty$, the left-hand side converges to $\E_x [\eta(Y_{t\wedge \tau_{U(r)}})]-\eta(x)$. For the right-hand side we can use Fatou's lemma (justified because of \eqref{e:exit-time-estimate-U-2}) to conclude that
\begin{align*}
\limsup_{k\to \infty}\E_x\int_0^t {\bf 1}_{s<\tau_{U_{a_k}(r)}} L^\sB g_k(Y_s)\, ds \le \E_x\int_0^t {\bf 1}_{s<\tau_{U(r)}} L^\sB \eta(Y_s)\, ds \le -\E_x (t\wedge \tau_{U(r)})\, .
\end{align*}
Thus we get that $\E_x [\eta(Y_{t\wedge \tau_{U(r)}})]-\eta(x) \le -\E_x (t\wedge \tau_{U(r)})$, and by letting $t\to \infty$,
$$
-\eta(x)\le \E_x [\eta(Y_{\tau_{U(r)}})]-\eta(x) \le -\E_x \tau_{U(r)}\, .
$$
Thus we get $\E_x \tau_{U(r)}\le \eta(x)\le x_d^p$. By using that $U(r_1)=r_1 U$ and \eqref{e:exit-time-scaling},
 for any $x\in U$,
$$
\E_x \tau_U=r_1^{-\alpha} \E_{r_1x}\tau_{r_1U}\le r_1^{-\alpha}(r_1x_d)^p= r_1^{p-\alpha} x_d^p\, .
$$
This proves the claim with $C_{36}=r_1^{p-\alpha} $. \qed

\vspace{.1in}
\textbf{Acknowledgment}:
We thank the referee for helpful comments on the first version of this paper.
Part of the research for this paper was done while the second-named author was visiting
Jiangsu Normal University, where he was partially supported by  a grant from
the National Natural Science Foundation of China (11931004) and by
the Priority Academic Program Development of Jiangsu Higher Education Institutions

 \vspace{.1in}
\begin{singlespace}

\end{singlespace}
\vskip 0.1truein

\parindent=0em

{\bf Panki Kim}

Department of Mathematical Sciences and Research Institute of Mathematics,

Seoul National University, Seoul 08826, Republic of Korea

E-mail: \texttt{pkim@snu.ac.kr}

\bigskip

{\bf Renming Song}

Department of Mathematics, University of Illinois, Urbana, IL 61801,
USA

E-mail: \texttt{rsong@math.uiuc.edu}

\bigskip

{\bf Zoran Vondra\v{c}ek}

Department of Mathematics, Faculty of Science, University of Zagreb, Zagreb, Croatia,

Email: \texttt{vondra@math.hr}

\end{document}